\RequirePackage[l2tabu, orthodox]{nag}
\documentclass[10pt]{article}%

\usepackage{amssymb}
\usepackage{amsfonts}
\usepackage{amsmath}
\usepackage{mathtools}
\usepackage{dsfont}
\usepackage{esint} 
\usepackage{upgreek}

\usepackage{amsthm}
\usepackage{hyperref}
\usepackage{cleveref}

\usepackage{float}
\usepackage{graphicx}
\usepackage{wrapfig}
\usepackage{subfig}
\usepackage{fancybox}
\usepackage{framed}
\usepackage[usenames,dvipsnames,svgnames,table]{xcolor}
\usepackage{esint}
\usepackage{caption}
\usepackage{epstopdf} 
\usepackage[all]{xy} 
\usepackage{tikz} 
\usepackage{tabularx}
\usepackage{afterpage}
\usepackage{arydshln}

\usepackage[titletoc, title]{appendix} 
\usepackage{titling} 
\usepackage{url} 
\usepackage{color}
\usepackage{enumerate, multirow, longtable} 

\usepackage[latin1]{inputenc}
\usepackage[a4paper]{geometry}
\usepackage{microtype}

\providecommand{\U}[1]{\protect\rule{.1in}{.1in}}
\newtheorem{theorem}{Theorem}[section]
\newtheorem{proposition}[theorem]{Proposition}
\newtheorem{lemma}[theorem]{Lemma}
\newtheorem{corollary}[theorem]{Corollary}

\newtheorem{remark}[theorem]{Remark}
\newtheorem{definition}[theorem]{Definition}

\newtheorem{conjecture}[theorem]{Conjecture}
\numberwithin{equation}{section}

\newcommand{\CC}{\mathbb{C}}
\newcommand{\DD}{\mathbb{D}}
\newcommand{\EE}{\mathbb{E}}

\newcommand{\NN}{\mathbb{N}}

\newcommand{\PP}{\mathbb{P}}

\newcommand{\RR}{\mathbb{R}}

\newcommand{\ZZ}{\mathbb{Z}}

\newcommand{\Ba}{ {\mathcal B }}

\newcommand{\Na}{ {\mathcal N }}

\newcommand{\Ga}{ {\mathcal G }}

\newcommand{\Oa}{ {\mathcal O }}

\newcommand{\Ma}{ {\mathcal M }}

\renewcommand{\Re}{\mathrm{Re}}

\DeclareMathOperator{\mom}{\mathrm{MoM}}

\title{On the critical-subcritical moments of moments of random characteristic polynomials: a GMC perspective}
\author{Jonathan P. Keating and Mo Dick Wong}
\date{\today}
\begin{document}
\maketitle

\abstract{
We study the `critical moments' of subcritical Gaussian multiplicative chaos (GMCs) in dimensions $d \le 2$. In particular, we establish a fully explicit formula for the leading order asymptotics, which is closely related to large deviation results for GMCs and demonstrates a similar universality feature. We conjecture that our result correctly describes the behaviour of analogous moments of moments of random matrices, or more generally structures which are asymptotically Gaussian and log-correlated in the entire mesoscopic scale. This is verified for an integer case in the setting of circular unitary ensemble, extending and strengthening the results of Claeys et al. and Fahs to higher-order moments.
}

\section{Introduction}
\subsection{Moments of moments of random characteristic polynomials}
There has been a good deal of interest in recent years in understanding asymptotically log-correlated Gaussian fields and their extremal processes. In the context of random matrix theory, these structures arise from the study of the logarithms of characteristic polynomials of different ensembles, such as Haar-distributed random matrices from the classical compact groups, the unitary ensemble of random Hermitian matrices and the complex Ginibre ensemble. This has led to a line of investigation into the so-called moments of moments, which capture multifractal properties of the underlying characteristic polynomial as well as non-trivial information about its maxima.

To motivate our discussion, let us consider the circular unitary ensemble (CUE). Suppose $U_N$ is an $N \times N$ Haar-distributed unitary matrix and $P_N(\theta) = \det(I - U_N e^{-i\theta})$ is the associated characteristic polynomial. For $k, s > 0$, we denote by
\begin{align*}
\mom_{\mathrm{U}(N)}(k, s) := \EE_{\mathrm{U}(N)} \left[ \left(\frac{1}{2\pi}\int_0^{2\pi} |P_N(\theta)|^{2s}d\theta\right)^{k}\right]
\end{align*}

\noindent the (absolute) moments of moments (MoMs) of $P_N(\theta)$.  We shall be interested in the behaviour of $\mom_{\mathrm{U}(N)}(k, s)$ as $N \to \infty$. Due to the multifractality of the random characteristic polynomials, this quantity exhibits a phase transition as one varies the values of $(k, s)$, and it was conjectured in \cite{FK2014} that
\begin{align}\label{eq:CUEmom}
\mom_{\mathrm{U}(N)}(k, s) \overset{N \to \infty}{\sim}
\begin{dcases}
\left(\frac{G(1+s)^2}{G(1+2s) \Gamma(1-s^2)}\right)^k \Gamma(1-ks^2) N^{ks^2}, & \text{if } k < 1/s^2, \\
c(k, s)  N^{k^2s^2 - k+1}, & \text{if } k > 1/s^2,
\end{dcases}
\end{align}

\noindent where $G(s)$ is the Barnes $G$-function, and $c(k, s)$ is some non-explicit constant depending on $(k, s)$.

Much progress has been made in establishing the asymptotics \eqref{eq:CUEmom}, especially when $k, s$ are both positive integers, in which case combinatorial approaches are possible. In this setting, the pair $(k, s)$ always\footnote{Except the relatively trivial case $k = s = 1$, where the moment of moment is simply reduced to $\EE_{\mathrm{U}(N)}\left[|P_N(\theta)|^2\right]$.} lies in the supercritical regime $k > 1/s^2$, and the corresponding asymptotics was first verified in \cite{BK2019}. An alternative combinatorial proof was also given in \cite{AK2020}, where a connection to the Painlev\'e V equation was established for $c(2, s)$, $s \in \NN$.

For positive integers $k$ and general $s$ satisfying $s > -1/4$, most results have been established with the help of the asymptotics of Toeplitz determinants obtained by  the Riemann-Hilbert approach\footnote{Some of the results hold for complex-valued $s$,  but this is beyond the scope of the discussion here.}. The case where $k=2$ was treated in \cite{CK2015}, where the asymptotic formulae in both regimes in \eqref{eq:CUEmom} were established alongside the connection of $c(2, s)$ to the Painlev\'e V equation. As for general $k \in \NN$, the subcritical case $s^2 < 1/k$ was fully established in \cite{Fa2019}, but only a weaker result $\mom_{U(N)}(k, s) = \Omega(N^{k^2 s^2 - k+1})$ was obtained in the supercritical regime.

Common to all the aforementioned works is the crucial assumption that $k$ is a positive integer, so that one can expand the moments and obtain
\begin{align}\label{eq:expandedmom}
\mom_{\mathrm{U}(N)}(k, s) 
= \frac{1}{(2\pi)^k} \int_0^{2\pi} \dots \int_0^{2\pi}\EE_{\mathrm{U(N)}} \left[\prod_{j=1}^k |P_N(\theta_j)|^{2s}\right] d\theta_1 \dots d\theta_k.
\end{align}

\noindent From there, the problem becomes  the analysis of the asymptotics for the cross moment in the integrand as well as its behaviour as one integrates it over $[0, 2\pi]^k$. It is not clear, however, how one may analyse the cross moment when $\theta_j$'s are arbitrarily close to each other (except when $k = 2$ \cite{CK2015}), let alone extend the investigation of moments of moments to non-integer values of $k$. Our goal here is to address this situation through the lens of Gaussian multiplicative chaos (GMC).

Before we continue, let us briefly mention the works \cite{ABK2020, BK2020a, BK2020b} where analogous moments of moments are studied in the context of classical compact groups, an approximate model for Riemann zeta function, and the branching random walk via the combinatorial approach, as well as the article \cite{CGMY2020} where the results in \cite{Fa2019} are translated to Toeplitz$+$Hankel determinants, yielding similar weak asymptotics for classical compact groups.

\subsection{Perspective from the theory of multiplicative chaos}\label{subsec:per}
The theory of GMCs studies the construction and properties of a random multifractal measure formally defined as the exponentiation of log-correlated Gaussian fields.
First proposed as a model of turbulence in the 60's, GMCs have attracted a lot of attention in the last two decades due to their central role in the mathematical formulation of Polyakov's path integral theory of random surfaces, also known as Liouville quantum gravity \cite{DKRV2016, DS2011}. The significance of this theory lies in the fact that it describes (at least conjecturally) the scaling limits of a wide class of models in statistical mechanics, and is closely related to other `canonical' geometric objects in the continuum with conformal symmetries, see e.g. \cite{DMS2014, Sh2016} and the references therein. Thanks to techniques from conformal field theory and complex analysis, finer results have been obtained in recent years including various exact formulae for Liouville correlation functions and related quantities \cite{KRV2020, Re2020, RZ2020a, RZ2020b}.

In a different direction, GMCs have been used as a tool to study the geometry of the underlying fields, and new applications in other areas of mathematics have also been studied over the last couple of years. For example, it is now known that the characteristic polynomials of many different random matrix ensembles behave asymptotically like GMC measures as the size of the matrix goes to infinity \cite{We2015, BWW2018, LOS2018, NSW2018, CFLW2019, FK2020, CN2019}. We would like to exploit such connections to GMCs, and explore the implications of seemingly unrelated exact integrability results from Liouville quantum gravity in the context of random matrix theory.

For illustration, let us revisit our CUE example. It has been shown in \cite{NSW2018} that when $s \in (0, 1)$, the sequence of random measures
\begin{align} \label{eq:cue2gmc}
\frac{|P_N(\theta)|^{2s}}{\EE_{\mathrm{U}(N)}|P_N(\theta)|^{2s}} d\theta
\end{align}

\noindent converges in distribution (and in the weak$^*$ topology) to the GMC measure $M_{\gamma}(d\theta) = e^{\gamma X(e^{i\theta}) - \frac{\gamma^2}{2} \EE[X(e^{i\theta})^2]}d\theta$, where $\gamma = \sqrt{2}s$ and $X$ is the Gaussian free field
\begin{align}\label{eq:cuegff}
\EE\left[X(x) X(y) \right] = -\log|x-y|, \qquad x, y \in \partial \DD = \text{unit circle},
\end{align}

\noindent as the size of the matrix $N$ goes to infinity. When $k < 1/s^2 = 2/\gamma^2$, the $k$-th moment of GMC exists and is given by the Fyodorov-Bouchaud formula \cite{FB2008, Re2020}
\begin{align*}
\EE \left[\left(\frac{1}{2\pi}M_{\gamma}([0, 2\pi])\right)^k\right] = \frac{\Gamma(1 - k \frac{\gamma^2}{2})}{\Gamma(1 - \frac{\gamma^2}{2})^k}.
\end{align*}

It follows from the work of Keating and Snaith \cite{KS2000} that
\begin{align}\label{eq:KScue}
\EE_{\mathrm{U}(N)}|P_N(\theta)|^{2s} \overset{N \to \infty}{\sim} \frac{G(1+s)^2}{G(1+2s)} N^{s^2}.
\end{align}

\noindent  Combining both results, we arrive at
\begin{align}
\notag
\mom_{\mathrm{U}(N)}(k, s) 
& =  \EE_{\mathrm{U}(N)} \left[ \left(\frac{1}{2\pi}\int_0^{2\pi} |P_N(\theta)|^{2s}d\theta\right)^{k}\right]\\
\notag
& = \left(\EE_{\mathrm{U}(N)} |P_N(\theta)|^{2s}\right)^k \EE_{\mathrm{U}(N)} \left[ \left(\frac{1}{2\pi} \int_0^{2\pi} \frac{|P_N(\theta)|^{2s}}{\EE_{\mathrm{U}(N)}|P_N(\theta)|^{2s}} d\theta\right)^{k}\right]\\
\label{eq:cue_subsub}
& \overset{N \to \infty}{\sim} \frac{\Gamma(1 - k s^2)}{\Gamma(1 - s^2)^k} \left[ \frac{G(1+s)^2}{G(1+2s)}\right]^k N^{ks^2}
\end{align}

\noindent which matches the first asymptotics in \eqref{eq:CUEmom}, 
 i.e. the regime $k < 1/s^2$, $0 < s < 1$ is essentially resolved\footnote{See \Cref{app:subsub} for further discussion.}.

The focus of this article is to extend the analysis beyond the regime where GMC moments exist, particularly when $k$ takes the critical value. In the CUE example this would correspond to the situation where $k\gamma^2/2 = ks^2 = 1$ with $s^2 < 1$. A naive substitution of $k = 1/s^2$ into \eqref{eq:cue_subsub} results in the blow-up of the factor $\Gamma(1-ks^2)$, justifying the need for a different asymptotic formula in this setting. It is widely conjectured that
\begin{align}\label{eq:cs_conjecture}
\mom_{\mathrm{U}(N)}\left(k, s = \frac{1}{\sqrt{k}}\right) \sim \alpha(k, \frac{1}{\sqrt{k}}) N \log N,
\end{align}

\noindent but nothing is known rigorously except for $k=1/s^2 = 2$ where it was shown in \cite[Theorem 1.15]{CK2015} that
\begin{align*}
\alpha\left(2, \frac{1}{\sqrt{2}}\right) = \frac{1}{\pi} \frac{G(1+\frac{1}{\sqrt{2}})^2}{G(1+\sqrt{2})}.
\end{align*}

\noindent In general, the best known result to date is
\begin{align*}
\log \mom_{\mathrm{U}(N)}\left(k, \frac{1}{\sqrt{k}}\right) = \log N + \log \log N + \mathcal{O}(1)
\end{align*} 

\noindent which was established for positive integers $k$ in \cite{Fa2019}, and to our knowledge there have been no conjectures giving an explicit formula for the $\Oa(1)$-term.

While one can no longer directly rely on the convergence of GMC in \cite{NSW2018} to derive an asymptotic formula for $\mom_{\mathrm{U}(N)}(k, s)$ in the regime where $ks^2 \ge 1$, one may still consider the following GMC heuristic: if $(X_{\epsilon}(\cdot))_{\epsilon > 0}$ is a collection of continuous Gaussian fields such that
\begin{align*}
M_{\gamma, \epsilon}(dx) := e^{\gamma X_{\epsilon}(x) - \frac{\gamma^2}{2}\EE[X_\epsilon(x)^2]}dx
\xrightarrow{\epsilon \to 0^+} M_{\gamma}(dx)
\end{align*}

\noindent and the covariance profiles of the log-characteristic polynomials are comparable to those of $X_{\epsilon}$ for suitably chosen $\epsilon = \epsilon(N)$, then perhaps
\begin{align}\label{eq:GMCheuristic}
\EE_{\mathrm{U}(N)}\left[ \left(\int_0^{2\pi} \frac{|P_N(\theta)|^{2s}}{\EE_{\mathrm{U}(N)}|P_N(\theta)|^{2s}} d\theta\right)^k\right]
\overset{\epsilon \to 0^+}{\sim} \EE\left[ \left(\int_0^{2\pi} M_{\gamma, \epsilon}(dx)\right)^k\right].
\end{align}

\noindent We shall therefore consider `critical moments' of (regularised) subcritical GMCs \eqref{eq:GMCheuristic} in this article, and study their asymptotics as $\epsilon \to 0^+$.

\subsection{Setting and main result}
For generality we consider $d \le 2$ so that our result is potentially relevant to  not only one-dimensional models like the classical compact groups, but also two-dimensional ones such as complex Ginibre ensemble \cite{WW2019} as well as other asymptotically log-correlated structures in other probability problems. 

Let $X$ be a log-correlated Gaussian field on some bounded open domain $D \subset \RR^d$ with covariance
\begin{align}
\label{eq:LGF}
\EE[X(x) X(y)] = -\log |x-y| + f(x, y) \qquad \forall x, y \in D
\end{align}

\noindent where $f$ is some continuous function on $\overline{D} \times \overline{D}$. 

For $\gamma^2 < 2d$, the Gaussian multiplicative chaos associated to $X$ is defined as the weak$^*$ limit
\begin{align}\label{eq:GMCdef}
M_{\gamma}(dx) 
:= \lim_{\epsilon \to 0^+} M_{\gamma, \epsilon}(x) 
= \lim_{\epsilon \to 0^+} e^{\gamma X_\epsilon(x) - \frac{\gamma^2}{2} \EE[X_\epsilon(x)^2]}dx
\end{align}

\noindent where $\left(X_{\epsilon}\right)_{\epsilon > 0}$ is some sequence of continuous Gaussian fields on $D$ that converges to $X$ as $\epsilon \to 0^+$ in a suitable sense. The covariance structure of such approximate fields are typically of the form
\begin{align*}
\EE[X_\epsilon(x) X_\epsilon(y)] = -\log \left(|x-y| \vee \epsilon \right) + f_{\epsilon}(x, y)
\xrightarrow{\epsilon \to 0} \EE[X(x) X(y)] \qquad \forall x, y \in D
\end{align*}

\noindent and it is an established fact in the GMC literature that the limiting measure $M_{\gamma}$ is independent of the choice of approximation (see \Cref{theo:GMCexist} in \Cref{subsec:GMC} for a version for the mollification construction, or \cite{JS2017} for a generalisation).

It is also well-known that the total mass of $M_\gamma$ possesses moments of all orders strictly less than $p_c := 2d/\gamma^2$, and for convenience let us call the pair $(p_c, \gamma)$ critical--subcritical. In particular, if $g \ge 0$ is a non-trivial\footnote{By which we mean $\int_D g > 0$.} continuous function on $\overline{D}$, then the sequence
\begin{align}\label{eq:cmom}
\EE \left[\left(\int_D g(x) M_{\gamma, \epsilon}(dx)\right)^{p_c}\right]
\end{align}

\noindent blows up as $\epsilon \to 0^+$ by Fatou's lemma. Our first theorem provides a precise description of this behaviour when the underlying approximate fields $X_\epsilon$ are constructed via mollification.

\begin{theorem}\label{theo:main}
For $d \le 2$ and $\gamma \in (0, \sqrt{2d})$, consider
\begin{align*}
M_{\gamma, \epsilon}(dx) = e^{\gamma X_\epsilon(x) - \frac{\gamma^2}{2}\EE[X_\epsilon(x)^2]}dx
\end{align*}  
\noindent with $X_\epsilon := X \ast \nu_\epsilon$, where $X$ is the log-correlated Gaussian field in \eqref{eq:LGF} and $\nu_\epsilon(\cdot):= \epsilon^{-d} \nu(\cdot/\epsilon)$ for some mollifier $\nu \in C_c^\infty(\RR^d)$.\footnote{By mollifiers we mean functions $\nu \in C_c^\infty(\RR^d)$ such that $\nu \ge 0$ and $\int \nu= 1$.} Suppose $g \ge 0$ is a non-trivial continuous function on $\overline{D}$. Writing $p_c = 2d/\gamma^2$, we have
\begin{equation}\label{eq:main}
\begin{split}
&\EE\left[\left(\int_D g(x) M_{\gamma, \epsilon}(dx)\right)^{p_c}\right]\\
& \qquad \overset{\epsilon \to 0^+}{\sim}
\left(\int_D e^{d(p_c-1) f(u, u)} g(u)^{p_c} du\right)\frac{\gamma^2}{2}(p_c-1)^2 \overline{C}_{\gamma, d}\log \frac{1}{\epsilon}.
\end{split}
\end{equation}

\noindent The constant $\overline{C}_{\gamma, d}$ is the reflection coefficient of GMC, which is explicitly given by
\begin{align}\label{eq:reflection0}
\overline{C}_{\gamma, d} 
= \begin{dcases}
\frac{(2\pi)^{\frac{2}{\gamma^2}-1}}{(1-\frac{\gamma^2}{2})\Gamma(1-\frac{\gamma^2}{2})^{\frac{2}{\gamma^2}}}, & d = 1,\\
-\frac{\left(\pi \Gamma(\frac{\gamma^2}{4}) / \Gamma(1- \frac{\gamma^2}{4})\right)^{\frac{4}{\gamma^2}- 1}}{\frac{4}{\gamma^2}-1}
 \frac{\Gamma(\frac{\gamma^2}{4} - 1)}{\Gamma(1-\frac{\gamma^2}{4}) \Gamma(\frac{4}{\gamma^2}-1)}, & d = 2.
\end{dcases}
\end{align}
\end{theorem}

The expression \eqref{eq:main} is reminiscent of the tail expansion of subcritical GMCs, as it was shown in \cite{Won2020} (under extra regularity condition on $f$) that
\begin{align}\label{eq:LD_GMC}
\PP\left(\int_D g(x) M_{\gamma}(dx) > t\right)
 \overset{t \to \infty}{\sim}
\left(\int_D e^{d(p_c-1) f(u, u)} g(u)^{p_c} du\right)
\left(1 -\frac{1}{p_c}\right)\frac{\overline{C}_{\gamma, d}}{t^{p_c}}.
\end{align}

\noindent This should not be entirely surprising given that the blow-up of the $p_c$-th moment is intrinsically related to the power law asymptotics of GMCs, and these results together seem to suggest that \eqref{eq:LD_GMC} may also hold for the approximate GMC $M_{\gamma, \epsilon}$ for $t = t(\epsilon)$ tending to infinity not too quickly with respect to $\epsilon \to 0^+$. We also see the recurring theme of universality for these large-deviation type problems, i.e. the asymptotics only depends on the diagonal of the covariance $f(u, u)$ but not any further specification of the underlying log-correlated field.

The asymptotic formula \eqref{eq:main} also looks relatively robust and should hold beyond the mollification setting, e.g. when $X_\epsilon$ is constructed by truncating the series expansion of $X$. We have the following corollary.
\begin{corollary}\label{cor:main}
Under the same setting as \Cref{theo:main}, except that $(X_{\epsilon})_{\epsilon > 0}$ is any collection of continuous Gaussian fields satisfying the following two conditions: 
\begin{itemize}
\item There exists some $C>0$ such that
\begin{align*}
\limsup_{\epsilon \to 0^+} \sup_{x, y \in D} \left|\EE[X_\epsilon(x) X_\epsilon(y)] + \log (|x-y| \vee \epsilon) \right| \le C.
\end{align*}

\item For any $\delta > 0$,
\begin{align*}
\lim_{\epsilon \to 0^+} \sup_{|x-y| \ge \epsilon^{1-\delta}} \left|\EE[X_\epsilon(x) X_\epsilon(y)] - \EE[X(x) X(y)]\right| & = 0.
\end{align*}
\end{itemize}

\noindent Then the asymptotics in \eqref{eq:main} remains true.
\end{corollary}

\paragraph{On circular unitary ensemble.}
We are able to establish the following asymptotics for Haar-distributed unitary matrices.

\begin{theorem}\label{theo:CUEinteger}
Let $k\ge 2$ be any integer. As $N \to \infty$, we have
\begin{align}\label{eq:cue_conj}
\mom_{\mathrm{U}(N)}\left(k, \frac{1}{\sqrt{k}}\right)
\overset{N \to \infty}{\sim}
\frac{k-1}{\Gamma(1-\frac{1}{k})^k}
\left[\frac{G(1+\frac{1}{\sqrt{k}})^{2}}{G(1+\frac{2}{\sqrt{k}})}\right]^k
N \log N.
\end{align}
\end{theorem}

Let us quickly explain how \eqref{eq:cue_conj} follows from the GMC heuristics. Since the log-characteristic polynomial behaves asymptotically like a log-correlated Gaussian field in the entire mesoscopic regime,
we anticipate that the GMC approximation would be accurate, and hence expect
\begin{align*}
\mom_{\mathrm{U}(N)}(k, s) 
& \overset{N \to \infty}{\sim} \frac{1}{(2\pi)^k}
\left(\EE_{\mathrm{U}(N)} |P_N(\theta)|^{2s}\right)^k \EE \left[ \left( \int_0^{2\pi} M_{\gamma, \epsilon(N)}(d\theta)\right)^{k}\right]
\end{align*}

\noindent to hold for any $k = 1/s^2 > 1$. Based on the discussion in \Cref{subsec:per}, it is reasonable to believe that the suitable log-correlated field for approximation is given by the Gaussian free field
\begin{align}\label{eq:hcue1}
\EE\left[X(e^{i\theta})X(e^{i\phi})\right]
= -\log |e^{i(\theta - \phi)} - 1|
\overset{\theta \to \phi}{=} -\log|\theta - \phi| + \mathcal{O}(\theta - \phi),
\end{align}

\noindent  in which case the $f$-function appearing in \eqref{eq:LGF} is simply zero on the diagonal. Meanwhile, the comparison
\begin{align*}
\Omega(\epsilon^{-\frac{k\gamma^2}{2}})
= \EE\left[e^{\frac{\gamma^2}{2}\EE[X_\epsilon(x)^2]}\right]^{k}
\approx \left(\EE |P_N(\theta)|^{2s}\right)^k
= \Omega(N^{ks^2})
\end{align*}

\noindent with $\gamma = \sqrt{2}s$ seems to suggest the choice $\epsilon = \epsilon(N) =\Omega(N^{-1})$.The asymptotics \eqref{eq:cue_conj} then follows by combining all these considerations with \Cref{theo:main}. When $k \ge 2$ is an integer, this heuristic argument can be made rigorous, since we may consider the expanded moments \eqref{eq:expandedmom} and show that
\begin{align*}
\EE_{\mathrm{U(N)}} \left[\prod_{j=1}^k |P_N(\theta_j)|^{2s}\right]
\qquad \text{and} \qquad
\EE \left[\prod_{j=1}^k e^{\gamma X_{\epsilon}(e^{i\theta_j})}\right]
\end{align*}

\noindent behave similarly as we integrate over $[0, 2\pi)^k$. This will be verified carefully in \Cref{subsec:CUEproof}. This strategy does not work for non-integers $k$ unfortunately, but we nevertheless conjecture that
\begin{conjecture}
The asymptotics \eqref{eq:cue_conj} hold for any $k > 1$.
\end{conjecture}

\subsection{Interpretation of results: mesoscopic statistics}\label{subsec:interpret}
Our main result of course does not cover the asymptotics for
$\EE\left[\left(\int_D g(x) M_{\gamma, \epsilon}(dx)\right)^{p}\right]$

\noindent in the `supercritical' cases where
\begin{itemize}
\item $\gamma^2 >2d$ and $p = \frac{2d}{\gamma^2}$; or
\item $\gamma \in \RR$ but $p > \frac{2d}{\gamma^2}$.
\end{itemize}

\noindent We would like to emphasise that this restriction is not due to any technical reason, but rather is largely due to our motivation for \eqref{eq:cmom} as an accurate toy model that provides the correct leading order.

As we shall see in the proofs, the GMC heuristic \eqref{eq:GMCheuristic} is highly plausible in the critical-subcritical setting as long as the original model behaves asymptotically like a log-correlated Gaussian field in the entire mesoscopic scale. Going beyond this phase, the leading order will depend on the microscopic behaviour, and there is no hope to identify the correct leading coefficients even if our GMC toy model still successfully captures the correct order. 
In the context of CUE, the microscopic statistics is described by the Sine kernel, and it should not be surprising that the present Gaussian analysis does not provide the right answer to the leading order coefficients. 
The interested readers may find more discussion in \Cref{app:supercritical} where we sketch some of the calculations in the supercritical case and explain the connections between the leading order asymptotics and the microscopic regime.

The same GMC heuristic may be applied to other models in probability, and we would like to mention two other examples in random matrix theory.

\paragraph{Circular beta ensemble (C$\beta$E).} For $\beta > 0$ and $N \in \NN$, the C$\beta$E$_N$ is a point process of $N$ particles $\theta_1, \dots, \theta_N$ on the unit circle with joint density
\begin{align*}
d\PP_{\mathrm{C\beta E}(N)}(\theta_1, \dots, \theta_N) \propto \prod_{j < k} |e^{i\theta_j} - e^{i \theta_k}|^\beta \prod_{1 \le k \le N} d\theta_k.
\end{align*}

\noindent When $\beta = 1, 2, 4$ this exactly corresponds to the eigenvalue distributions of random matrices sampled from the circular orthogonal, unitary and symplectic ensembles respectively, but a general construction is also available for arbitrary $\beta > 0$ in \cite{KN2004}. If we define the `characteristic polynomial' $P_N(\theta) = \prod_{j=1}^N (1 - e^{i(\theta_j-\theta)})$ as before, then we know from \cite{La2019} that
\begin{align*}
\log |P_N(\theta)| \to \frac{1}{\sqrt{\beta}} X(e^{i\theta})
\end{align*}

\noindent in the sense of finite dimensional distribution as $N\to\infty$, where $X(\cdot)$ is again the Gaussian free field on the unit circle in \eqref{eq:cuegff}. Partial progress has been made towards proving the convergence to GMCs in the same article \cite{La2019}; see also \cite{CN2019} where a different connection between C$\beta$E-characteristic polynomials and GMCs was established.

Not much is known about the asymptotics for moments of moments of C$\beta$E models in general, since the underlying point process is no longer determinantal for $\beta \ne 2$ and powerful techniques such as Riemann-Hilbert analysis no longer apply. The only work in this direction is the recent article \cite{As2020} where the asymptotics are established for part of the supercritical regime $\beta < 2ks^2$ with integral parameters $k$ and $s$, based on symmetric function theory.

We consider the `critical-subcritical' regime based on the GMC heuristics and our main theorem.
Let us recall from the work of Keating and Snaith \cite{KS2000} that
\begin{align*}
\EE_{\mathrm{C\beta E}(N)}\left[|P_N(\theta)|^{2s}\right] 
& = \prod_{k=0}^{N-1} \frac{\Gamma(1+\frac{\beta}{2}k) \Gamma(1+2s + \frac{\beta}{2}k)}{\Gamma(1+s+\frac{\beta}{2}k)^2}.
\end{align*}

\noindent Using the leading order asymptotics derived in \cite[Lemma 4.17]{BHR2019}, we arrive at the following conjecture.
\begin{conjecture}
For any $s> 0$ and $\beta > 0$ satisfying $2s^2 < \beta$ and $k = \frac{\beta}{2s^2}$,
\begin{align}
\notag
\mathrm{MoM}_{\mathrm{C\beta E}(N)}(k, s)
&:= \EE_{\mathrm{C\beta E}(N)}\left[\left(\frac{1}{2\pi} \int_0^{2\pi} |P_N(\theta)|^{2s}d\theta\right)^{k}\right]\\
\label{eq:cbe_conj}
& \overset{N \to \infty}{\sim} \frac{k-1}{\Gamma(1-\frac{1}{k})^k} \left[\frac{e^{\Upupsilon^\beta(1- \frac{\beta}{2}) + \Upupsilon^\beta(1+2s - \frac{\beta}{2})}}{e^{2\Upupsilon^\beta(1 + s - \frac{\beta}{2})}}\right]^k N \log N
\end{align}

\noindent where
\begin{align*}
\Upupsilon^\beta(x) 
&= \frac{\beta}{2} \log G\left(1+\frac{2x}{\beta}\right) - \left(x - \frac{1}{2}\right) \log \Gamma\left(1+\frac{2x}{\beta}\right)\\
& \qquad + \int_0^\infty \left[\frac{1}{2t} - \frac{1}{t^2} + \frac{1}{t(e^t-1)}\right] \frac{e^{-xt - 1}}{e^{\frac{\beta}{2}t}-1} dt + \frac{x^2}{\beta} + \frac{x}{2}
\end{align*}

\noindent and $G$ is the Barnes G-function.
\end{conjecture}

\begin{remark}
At first glance it is not obvious how \eqref{eq:cbe_conj} coincides with our earlier conjecture regarding the CUE when $\beta = 2$, but this may be deduced from the equivalent and simpler representation of $\Upupsilon^\beta(\cdot)$ for rational values of $\beta$, see \cite[Lemma 7.1]{BHR2019}.
\end{remark}

\paragraph{Classical compact groups.} Following the spirit of previous works on unitary matrices, the convergence of characteristic polynomials to GMC measures has recently been established in \cite{FK2020} for other classical compact groups, namely the orthogonal group $\mathrm{O}(N)$ and symplectic group $\mathrm{Sp}(2N)$\footnote{By $\mathrm{Sp}(2N)$ we mean the collection of unitary symplectic matrices of size $2N \times 2N$.}. More precisely, if $\Ga(N) = \mathrm{O}(N)$ or $\mathrm{Sp}(N)$ (the latter case only makes sense when $N$ is even), and $P_N(\theta) = \det(I - U_N e^{-i\theta})$ where $U_N \in \Ga(N)$ is Haar-distributed, then
\begin{align*}
\frac{|P_N(\theta)|^{2s}}{\EE_{\Ga(N)} |P_N(\theta)|^{2s}} ds
\xrightarrow[N\to\infty]{d} e^{\gamma X(e^{i\theta}) - \frac{\gamma^2}{2} \EE[X(e^{i\theta})^2]}d\theta
\end{align*}

\noindent where $\gamma = \sqrt{2}s$ as before but now the Gaussian field $X$ has covariance
\begin{align}\label{eq:GFF_OS}
\EE[X(e^{i\theta}) X(e^{i\phi})] = -\log \left|\frac{e^{i\theta} - e^{i\phi}}{e^{i\theta} - e^{-i\phi}} \right|, \qquad \forall \theta, \phi \in [0, 2\pi).
\end{align}

\noindent While this log-correlated field has an extra singularity at the complex conjugate and does not exactly satisfy \eqref{eq:LGF}, we believe that the GMC heuristic would work equally well here\footnote{The main reasoning is that blow-up of the `critical moment' should be attributed to very positive covariance (i.e. near the diagonal), which is opposite to the `anti-diagonal' behaviour, i.e. that the covariance tends to minus infinity when $\theta \to -\phi$.}, provided that the leading order is not dominated by the contribution from the neighbourhood of $\theta = 0$ or $\theta = \pi$ where the behaviour of the field $X$ is different.

To apply the GMC heuristic, we again extract the value of the $f$-function on the diagonal: as $\theta \to \phi$,
\begin{align*}
 -\log \left|\frac{e^{i\theta} - e^{i\phi}}{e^{i\theta} - e^{-i\phi}} \right|
\to -\log|\theta - \phi| + \log |1-e^{2i\phi}| + \Oa(\theta - \phi)
\end{align*}

\noindent everywhere except for $\phi \in \{0, \pi\}$, which is a set of measure zero and hence may be neglected, i.e. we may just take $f(\theta, \theta) = \log |1- e^{2i\theta}|$.  Unlike the unitary case, we also need to pick a suitable $g$-function for \eqref{eq:main}, since $\log|P_N(\theta)|$ converges to a non-stationary Gaussian field with non-zero mean (see \cite[Theorem 3]{FK2020}). For orthogonal matrices, it is known (e.g. \cite[equations (4.32)]{FK2020}) that
\begin{align*}
\EE_{\mathrm{O}(N)}|P_N(\theta)|^{2s} \overset{N \to \infty}
= (1+o(1)) N^{s^2} \frac{G(1+s)^2}{G(1+2s)} |1-e^{2i\theta}|^{-s^2 +s},
\qquad \text{for } \theta \not \in \{0, \pi\}
\end{align*}

\noindent and so we may want to pick $g(\theta) = N^{s^2} \frac{G(1+s)^2}{G(1+2s)} 
|1-e^{2i\theta}|^{-s^2 +s}$. Note that this is a continuous function without any singularity when $s \le 1$. Therefore,
\begin{conjecture}
For any $k>1$ and $s = \frac{1}{\sqrt{k}}$,
\begin{align*}
&\mathrm{MoM}_{\mathrm{O}(N)}(k, s)
:= \EE_{\mathrm{O}(N)} \left[\left(\frac{1}{2\pi}\int_0^{2\pi} |P_N(\theta)|^{2s}\right)^k\right]\\
& \qquad \overset{N \to \infty}{\sim} \frac{k-1}{\Gamma(1-\frac{1}{k})^k}
\left(\frac{1}{2\pi}\int_0^{2\pi} |1-e^{2i\theta}|^{k+\sqrt{k}-2}d\theta\right)
\left[\frac{G(1+\frac{1}{\sqrt{k}})^{2}}{G(1+\frac{2}{\sqrt{k}})}\right]^k
N \log N.
\end{align*}
\end{conjecture}

For the symplectic case, we have (see e.g. \cite[equations (4.18)]{FK2020}) 
\begin{align}\label{eq:sym_abs}
\EE_{\mathrm{Sp}(2N)}|P_{2N}(\theta)|^{2s} \overset{N \to \infty}
= (1+o(1)) (2N)^{s^2} \frac{G(1+s)^2}{G(1+2s)}  
|1-e^{2i\theta}|^{-s^2 -s}
\end{align}

\noindent for $\theta \not \in \{0, \pi\}$. Naturally one would like to pick $g(\theta) = (2N)^{s^2} \frac{G(1+s)^2}{G(1+2s)}  |1-e^{2i\theta}|^{-s^2 -s}$, but unlike the orthogonal case the `density function' now has singularities at $\theta = 0$ or $\pi$. This violates our assumption on $g$ in \Cref{theo:main}, but perhaps the GMC heuristic may be extended to the current case provided that the integral on the RHS of \eqref{eq:main} remains finite, i.e. the leading order asymptotics is not dominated by local contributions near $\theta \in \{0, \pi\}$, so that the non-uniformity of \eqref{eq:sym_abs} or the different behaviour of the covariance kernel \eqref{eq:GFF_OS} near these singularities are irrelevant. Interestingly this gives rise to a `singularity threshold' that is determined by the golden ratio:
\begin{conjecture}
For any $k>(\frac{1+\sqrt{5}}{2})^2$ and $s = \frac{1}{\sqrt{k}}$,
\begin{align*}
&\mathrm{MoM}_{\mathrm{Sp}(2N)}(k, s)
:= \EE_{\mathrm{Sp}(2N)} \left[\left(\frac{1}{2\pi}\int_0^{2\pi} |P_{2N}(\theta)|^{2s}\right)^k\right]\\
& \quad \overset{N \to \infty}{\sim} \frac{k-1}{\Gamma(1-\frac{1}{k})^k}
\left(\frac{1}{2\pi}\int_0^{2\pi} |1-e^{2i\theta}|^{k-\sqrt{k}-2}d\theta\right)
\left[\frac{G(1+\frac{1}{\sqrt{k}})^{2}}{G(1+\frac{2}{\sqrt{k}})}\right]^k
2N \log N.
\end{align*}
\end{conjecture}

\begin{remark}
The leading order of the critical-subcritical moments of moments is conjectured to be $\Omega(N\log N)$ for the three classical compact groups $\mathrm{U}(N)$, $\mathrm{O}(N)$ or $\mathrm{Sp}(N)$. This is in stark contrast to the results in \cite{ABK2020} where the exponent of $N$ depends on the underlying group, but is not surprising because of different singularity behaviour of the `density function' $g$, which plays a role in determining the size of the leading order in these supercritical cases.
\end{remark}

\paragraph{A note on Gaussian branching random walk.}
In \cite{BK2020b}, the asymptotics for moments of moments were considered for a toy model based on Gaussian branching random walk.  For instance, if we take $\gamma = \sqrt{2} \beta \ne  0$, $k= 1/\beta^2 \in \NN$ and $\epsilon = 2^{-n}$, then \cite[Theorem 2.1]{BK2020b} says that
\begin{align*}
\EE\left[\left(\int_0^1 e^{\gamma \widetilde{X}_\epsilon(x) - \frac{\gamma^2}{2}\EE[\widetilde{X}_\epsilon(x)^2]} dx\right)^k \right]
\sim \frac{\sigma(k)}{\log 2} \log \frac{1}{\epsilon}
\qquad \text{as } n \to \infty
\end{align*}

\noindent where $\widetilde{X}_\epsilon$ is a centred Gaussian field on $D = (0,1)$ with covariance
\begin{align*}
\EE[\widetilde{X}_{2^{-n}}(x) \widetilde{X}_{2^{-n}}(y)] = -\log \left(d(x, y) \vee 2^{-n}\right)
\end{align*}

\noindent and  $d(x, y) = 2^{-m}$ if the binary expansions of $x$ and $y$ only agree up to the first $m$ digits.\footnote{For definiteness one may assume that dyadic points are represented by finite-length binary expansions.}  Neither \Cref{theo:main} nor \Cref{cor:main} would apply to this field\footnote{It is expected that the leading order coefficients obtained in this model do not correspond to the analogous constants in GMC/random matrix theory, regardless of the regime under consideration.} since $d(x, y)$ is not continuous with respect to the Euclidean distance (e.g. consider $x < 1/2 < y$ with $x, y$ arbitrarily close to $1/2$) and $-\log d(x, y)$ is not of the form \eqref{eq:LGF}.  A one-sided bound may still be obtained, however, if we use \Cref{theo:main} with a different Gaussian field $X_\epsilon = X \ast \nu_\epsilon$ where
\begin{align*}
\EE[X(x)X(y)] = -\log |x-y| \qquad \forall x, y \in D.
\end{align*}

\noindent Indeed, as $|x-y| \le d(x, y)$, the comparison principle from \Cref{lem:Kahane} suggests that
\begin{align*}
\frac{\sigma(k)}{\log 2} \le \left(\int_D  du\right)\frac{\gamma^2}{2}(p_c-1)^2 \overline{C}_{\gamma, 1} = (2\pi)^{k-1} \frac{k-1}{\Gamma(1-\frac{1}{k})^k}
\end{align*}

\noindent and this inequality is most likely to be a strict inequality.

\subsection{Main idea}\label{subsec:idea}
Despite the connections to large deviation of GMC, our proof is orthogonal to that in \cite{Won2020}. As we shall see later, the analysis here is closely related to the asymptotics for
\begin{align}\label{eq:BMexpf}
\EE\left[\left(\int_0^T e^{\gamma (B_t - \mu t)} dt\right)^{p_c-1}\right]
\qquad \text{as} \qquad T \to \infty,
\end{align}

\noindent where $(B_t)_{t \ge 0}$ is a standard Brownian motion and $\mu = \frac{\gamma}{2}(p_c-1)> 0$. This guiding toy model has been used previously in the literature, such as the renormalisability of Liouville CFT at the Seiberg bound \cite{DKRV2017} and fusion estimates for GMC \cite{BW2018}, but these papers deal with negative moments of exponential functionals of Brownian motion with non-negative drift, the behaviours of which are drastically different from those in the current article.

Going back to \eqref{eq:BMexpf}, we note that $\int_0^\infty e^{\gamma (B_t - \mu t)}dt$ is actually almost surely finite. To understand what causes the blow-up when we compute the $(p_c-1)$-th moment, we make use of Williams' path decomposition theorem (\Cref{theo:path_dec}), which says that $(B_t - \mu t)_{t \ge 0}$ 
\begin{itemize}
\item behaves like a Brownian motion with positive drift $\mu$ (which is denoted by $(B_s^\mu)_{s \ge 0}$) until it hits a height $\Ma$ that is independently distributed according to $\mathrm{Exp}(2\mu)$;
\item runs like $\Ma - \Ba_{t - \tau_\Ma}^\mu$ when one goes beyond the first hitting time $\tau_{\Ma}$, where $(\Ba_s^\mu)_{s \ge 0}$ is a Brownian motion with drift $\mu$ conditioned to stay positive.
\end{itemize}

With these ideas in mind, we can rewrite our exponential functional as
\begin{align*}
\EE\left[e^{\gamma(p_c-1)\Ma} \left( \int_{-\log r}^{T \wedge \tau_\Ma} e^{\gamma (B_t^\mu -\Ma) } dt + \int_0^{(T-\tau_\Ma)_+} e^{- \gamma \Ba_{t }^\mu}  dt\right)^{p_c-1} \right].
\end{align*}

\noindent and attribute the blow-up behaviour to the presence of the term $e^{\gamma(p_c-1)\Ma} = e^{2\mu \Ma}$ in the expectation. 

The use of Williams' result is inspired by an earlier paper \cite{RV2019} on the tail expansions of 2-dimensional GMC measures, but the analysis there concerns negative moments and is closer in spirit to our ``continuity from below" discussion in \Cref{lem:lower} or \Cref{app:lowercont}. Our main theorem, however, concerns positive moments which require a different treatment, and it is crucial to identify the range of $\Ma$ that contributes to the leading order of the exponential functional. Not surprisingly, this is given by $\{\Ma \le \mu T\}$ as supported by the intuition of law of large numbers.

\paragraph{Organisation of the paper.} The remainder of this article is organised as follows.
\begin{itemize}
\item In \Cref{sec:prelim}, we compile a list of results that will be used in our proofs. This includes a collection of methods for the analysis of Gaussian processes, as well as various  facts about path distribution of Brownian motions with drift.

\item Our main proofs are presented in \Cref{sec:mainproof}. We first study a slightly unrelated \Cref{lem:lower} to illustrate how the path decomposition results may be used in our analysis. We then go back to \Cref{theo:main}, showing how the general problem may be reduced to one concerning the exactly log-correlated field from which we obtain our main result. After establishing a technical estimate, we then extend our result to \Cref{cor:main}. The section concludes with the proof of \Cref{theo:CUEinteger} on CUE moments of moments.

\item For completeness, we explain in the appendices the current status of the subcritical-subcritical asymptotics \eqref{eq:cue_subsub}, and how the leading order coefficient in \Cref{theo:main} may be seen as the renormalised limit of GMC moments from below. We also discuss the reason why GMC heuristics \eqref{eq:GMCheuristic} for random matrix models should fail completely as one goes beyond the critical-subcritical regime, and comment briefly on the validity of the strong asymptotics for Toeplitz determinants with slowly merging singularities which is used in \Cref{subsec:CUEproof}.
\end{itemize}

\paragraph{Acknowledgements}
We would like to thank Benjamin Fahs for helpful discussions on Toeplitz determinants with merging singularities, and Joseph Najnudel for useful comments on an earlier draft of this article. We are also most grateful to an anonymous referee for pointing out some of the typos in a previous version of this article, as well as for providing us with many remarks that helped improve the presentation.

\section{Preliminaries}\label{sec:prelim}
\subsection{Gaussian processes}
We start with our Gaussian toolbox. The first result we need is a standard change-of-measure lemma.
\begin{lemma}[Cameron-Martin]\label{lem:CMG}
Let $N, X(\cdot)$ be centred and jointly Gaussian. Then
\begin{align}\label{eq:CMG}
\EE \left[e^{N - \frac{1}{2}\EE[N^2]} F(X(\cdot))\right]
=\EE \left[F(X(\cdot) + \EE[N X(\cdot)])\right].
\end{align}
\end{lemma}

Next, we state Kahane's interpolation formula.
\begin{lemma}\label{lem:Kahane}
Let $\rho$ be a Radon measure on $D$ and $F: \RR_+ \to \RR$ some smooth function with at most polynomial growth at infinity, and suppose $G_0, G_1$ are two centred continuous Gaussian fields on $D$. For $t \in [0,1]$, define $G_t(x) := \sqrt{t} G_1(x) + \sqrt{1-t} G_0(x)$ and
\begin{align*}
\varphi(t) := \EE\left[F(W_t)\right] 
\qquad \text{where}\qquad W_t := \int_D e^{G_t(x) - \frac{1}{2} \EE[G_t(x)^2]} \rho(dx).
\end{align*}

\noindent Then the derivative of $\varphi$ is given by
\begin{align}
\notag 
\varphi'(t) 
& = \frac{1}{2} \int_D \int_D\left(\EE[G_0(x) G_0(y)] - \EE[G_1(x) G_1(y)]\right)\\
\label{eq:interpolate}
& \qquad \qquad \times \EE\left[e^{G_t(x) + G_t(y) - \frac{1}{2} \EE[G_t(x)^2] - \frac{1}{2} \EE[ G_t(y)^2] }F''(W_t)\right] \rho(dx) \rho(dy). 
\end{align}

\noindent In particular, if
\begin{align*}
\EE[G_0(x) G_0(y)] \le \EE[G_1(x) G_1(y)] \qquad \forall x, y \in D,
\end{align*}

\noindent then for any (integrable) convex function $F$ and (non-negative) measure $\mu$ on $D$, we have
\begin{align}\label{eq:Kahane}
\EE \left[ F\left(\int_D e^{G_0(x) - \frac{1}{2}\EE[G_0(x)^2]} \rho(dx) \right)\right]
\le \EE \left[ F\left(\int_D e^{G_1(x) - \frac{1}{2}\EE[G_1(x)^2]} \rho(dx)\right)\right].
\end{align}
\end{lemma}

A simple corollary to \Cref{lem:Kahane} is as follows:
\begin{corollary}\label{cor:GP}
Suppose $G_0, G_1$ are two continuous centred Gaussian fields on $D$ such that
\begin{align*}
\EE[G_0(x) G_0(y)] - \EE[G_1(x) G_1(y)] \le C \qquad \forall x, y \in D
\end{align*}

\noindent for some $C > 0$, then for any (non-negative) measure $\mu$, we have
\begin{align*}
& \EE \left[ \left(\int_D e^{\gamma G_0(x) - \frac{\gamma^2}{2}\EE[G_0(x)^2]} \mu(dx) \right)^k\right]\\
& \qquad \ge e^{\frac{\gamma^2}{2}k(k-1)C} 
\EE \left[ \left(\int_D e^{\gamma G_1(x) - \frac{\gamma^2}{2}\EE[G_1(x)^2]} \mu(dx) \right)^k\right]
\end{align*}

\noindent for $k \in [0, 1]$, and the inequality is reversed for $k \not \in [0,1]$. In particular, if the two-sided inequality
\begin{align*}
|\EE[G_0(x) G_0(y)] - \EE[G_1(x) G_1(y)]| \le C \qquad \forall x, y \in D
\end{align*}

\noindent is satisfied instead, then for any $k \in\RR$ we have
\begin{align*}
\frac{\EE \left[ \left(\int_D e^{\gamma G_0(x) - \frac{\gamma^2}{2}\EE[G_0(x)^2]} \mu(dx) \right)^k\right]}
{\EE \left[ \left(\int_D e^{\gamma G_1(x) - \frac{\gamma^2}{2}\EE[G_1(x)^2]} \mu(dx) \right)^k\right]}
\in [ e^{-\frac{\gamma^2}{2}|k(k-1)| C}, e^{\frac{\gamma^2}{2}|k(k-1)| C}].
\end{align*}

\end{corollary}

\begin{proof}
Suppose $k \in [0,1]$ such that $x \mapsto x^k$ is concave. Let $\Na_C$ be an independent Gaussian random variable with variance $C$, then
\begin{align*}
\EE\left[(G_1(x) + \Na_C)(G_1(y) + \Na_C)\right]
= \EE\left[G_1(x) G_1(y)\right] + C \ge \EE[G_0(x) G_0(y)].
\end{align*}

\noindent By Kahane's inequality, we have
\begin{align*}
&\EE\left[\left(\int_D e^{\gamma G_0(x) - \frac{\gamma^2}{2} \EE[G_0(x)^2]}\mu(dx)\right)^k\right]\\
&\qquad \ge \EE\left[\left(\int_D e^{\gamma (G_1(x) + \Na_C) - \frac{\gamma^2}{2} \EE[(G_1(x)+\Na_C)^2]}\mu(dx)\right)^k\right]\\
& \qquad = e^{\frac{\gamma^2}{2}k(k-1)C} \EE\left[\left(\int_D e^{\gamma G_1(x) - \frac{\gamma^2}{2} \EE[G_1(x)^2]}\mu(dx)\right)^k\right].
\end{align*}

\noindent The case where $k \not\in[0,1]$ is similar.
\end{proof}

\subsection{On Brownian motions and $3$-dimensional Bessel processes}\label{subsec:path_dec}
This subsection collects a few important results regarding the path distributions of Brownian motions as well as $3$-dimensional Bessel processes starting from $0$ (or $\mathrm{BES}_0(3)$-processes in short). We shall be using these notations throughout the entire paper.
\begin{itemize}
\item $(B_t)_{t \ge 0}$ is a standard Brownian motion.
\item $(\beta_t)_{t \ge 0}$ is a $\mathrm{BES}_0(3)$-process.
\item For $\mu \ge 0$:
\begin{itemize}
\item $(B_t^\mu)_{t \ge 0}$ is a Brownian motion with drift $\mu$.
\item $(\Ba_t^\mu)_{t \ge 0}$ is a Brownian motion with drift $\mu$ conditioned to stay non-negative, and $(\Ba_t^\mu)_{t \in \RR}$ a two-sided version of it, i.e. $(\Ba_t^\mu)_{t \ge 0}$ and $(\Ba_{-t}^\mu)_{t \ge 0}$ are two independent copies.
\end{itemize}
\end{itemize}

To be more precise:
\begin{definition}
$(\Ba_t^\mu)_{t \ge 0}$ is a Markov process starting from the origin with infinitesimal generator
\begin{align*}
\frac{1}{2}\frac{d^2}{dx^2} + \mu \coth(\mu x) dx, 
\end{align*}

\noindent or equivalently solution to the singular stochastic differential equation
\begin{align} \label{eq:BMcon}
d\Ba_t^\mu = \mu\coth(\mu \Ba_t^\mu) dt + dW_t, \qquad \Ba_0^\mu = 0
\end{align}

\noindent where $(W_t)_{t \ge 0}$ is the Wiener process. These definitions extend to the case $\mu = 0$, where $(\Ba_t^0)_{t \ge 0}$ should be interpreted as a $\mathrm{BES}_0(3)$-process $(\beta_t)_{t \ge 0}$ which solves
\begin{align*}
d\beta_t = \frac{1}{\beta_t} dt + dW_t, \qquad \beta_0 = 0.
\end{align*}

\noindent (In particular, despite its colloquial name, the process $(\Ba_t^\mu)_{t \ge 0}$ and similarly its two-sided analogue have discontinuous drift at $t = 0$.)
\end{definition}

Note that $(\Ba_t^\mu)_{t \ge 0}$ is known under different names in the literature, such as \emph{the radial part of a $3$-dimensional Brownian motion with drift $\mu \ge 0$}. We prefer to call it a conditioned drifting Brownian motion because this is the more commonly known name in the GMC--related literature. Indeed $(\Ba_t^\mu)_{t \ge 0}$ can be constructed via such conditioning procedure in the sense of Doob's theory of $h$-transform, see e.g. \cite[Section 2.4]{Wil1974} and \cite[Section 3]{RP1981} and the references therein.

Let us begin with an elementary lemma on stochastic dominance.
\begin{lemma}\label{lem:stocdom}
$(\Ba_t^{\mu_1})_{t \ge 0}$ is stochastically dominated by $(\Ba_t^{\mu_2})_{t \ge 0}$ whenever $\mu_1 \le \mu_2$. This comparison may be extended to $\mu_1 = 0$, in which case we interpret $(\Ba_t^{0})_{t \ge 0}$ as a $\mathrm{BES}_0(3)$-process $(\beta_t)_{t \ge 0}$.
\end{lemma}

\begin{proof}
The stochastic order should immediately follow from a suitable conditioning construction, but we nevertheless check it by showing that the drift in \eqref{eq:BMcon} is monotonically increasing in $\mu$: if we write
\begin{align*}
m_x(\mu) := \mu \coth(\mu x) = \mu \frac{e^{\mu x} + e^{- \mu x}}{e^{\mu x} - e^{-\mu x}}
\end{align*}

\noindent for any fixed $x > 0$, then
\begin{align*}
\frac{\partial}{\partial \mu} m_x(\mu)
& = \frac{e^{\mu x} + e^{- \mu x}}{e^{\mu x} - e^{-\mu x}}
+  \mu x \frac{(e^{\mu x} - e^{- \mu x})^2 - (e^{\mu x} + e^{- \mu x})^2}{(e^{\mu x} - e^{-\mu x})^2}\\
& = \frac{e^{2\mu x} - e^{-2\mu x} - 4\mu x}{(e^{\mu x} - e^{- \mu x})^2} \ge 0
\end{align*}

\noindent since $\partial_\mu (e^{2\mu x} - e^{-2\mu x} - 4\mu x) \ge 0$ for $\mu x \ge 0$.
\end{proof}

We now mention the generalised ``2M-B" theorem due to Rogers and Pitman, which provides a duality between $(B_t^\mu)_{t \ge 0}$ and $(\Ba_t^\mu)_{t \ge 0}$.
\begin{theorem}[{cf. \cite[Theorem 1]{RP1981}}]\label{theo:2m-b}
Let $\mu \ge 0$, and write
\begin{align*}
S_t^\mu := \max_{s \le t} B_s^\mu,
\qquad J_t^\mu := \inf_{s \ge t} \Ba_s^\mu.
\end{align*}

\noindent Then we have
\begin{align}
(2S_t^\mu - B_t^\mu, S_t^\mu)_{t \ge 0}
\overset{d}{=}
(\Ba_t^\mu, J_t^\mu)_{t \ge 0}
\end{align}

\noindent where the $\mu = 0$ case also holds if we interpret $(\Ba_t^{\mu=0})_{t \ge 0}$ as a $\mathrm{BES}_0(3)$-process $(\beta_t)_{t \ge 0}$ accordingly.
\end{theorem}

A simple but useful corollary of the above result is the following.
\begin{corollary}\label{cor:indep}
Let $\mu \ge 0$ and $x > 0$. If we define
\begin{align*}
T_x := \sup \{t > 0: \Ba_t^\mu = x\},
\end{align*}

\noindent then the process $(\Ba_{T_x + t}^\mu - \Ba_{T_x}^\mu)_{t \ge 0}$ is independent of $(\Ba_t^\mu)_{t \le T_x}$ and again distributed as $(\Ba_t^\mu)_{t \ge 0}$.
\end{corollary}

\begin{proof}
Given $(\Ba_t^\mu)_{t \ge 0}$, we define $J_t^\mu = \inf_{s \ge t} \Ba_s^\mu$ and set
\begin{align*}
B_t^\mu := 2J_t^\mu - \Ba_t^\mu,
\qquad S_t^\mu = \max_{s \le t} B_s^\mu.
\end{align*}

\noindent Then according to \Cref{theo:2m-b}, $(B_t^\mu)_{t \ge 0}$ is a Brownian motion with drift $\mu$, and $J_t^\mu = S_{t}^\mu$ for any $t>0$.
We then observe that
\begin{align*}
T_x = \inf \{t > 0: J_t^\mu = x\}
= \inf \{t > 0: S_t^\mu = x\}
\end{align*}

\noindent is a stopping time for $(B_t^\mu)_{t \ge 0}$. By the strong Markov property, the process 
\begin{align*}
(\widetilde{B}_t^\mu)_{t \ge 0}
= (B_{T_x + t}^\mu - B_{T_x}^\mu)_{t \ge 0}
\end{align*}

\noindent is independent of $(B_t^\mu)_{t \le T_x}$, or equivalently $(\Ba_t^\mu)_{t \le T_x}$, and it is clear that $(\widetilde{B}_t^\mu)_{t \ge 0} \overset{d}{=} (B_t^\mu)_{t \ge 0}$.  Combining this with the fact that $\Ba_{T_x}^\mu = x = J_{T_x}^\mu = B_{T_x}^\mu$, we obtain
\begin{align*}
\Ba_{T_x+t}^\mu - \Ba_{T_x}^\mu
&= (2J_{T_x + t}^\mu - B_{T_x + t}^\mu) - (2J_{T_x}^\mu - B_{T_x}^\mu)\\
& = 2(S_{T_x + t}^\mu - S_{T_x}^\mu) - (B_{T_x + t}^\mu-B_{T_x}^\mu)
= 2\widetilde{S}_{t}^\mu - \widetilde{B}_t^\mu
\end{align*}

\noindent where $\widetilde{S}_t^\mu := S_{T_x + t}^\mu - S_{T_x}^\mu = \max_{s \le t} \widetilde{B}_{s}^\mu$. The result follows from one final application of \Cref{theo:2m-b}.
\end{proof}

Note that \Cref{theo:2m-b} is intrinsically related to the work of Williams \cite{Wil1974} (see Remarks in \cite[Section 5]{RP1981} for a discussion of their connections), and we would like to highlight two important consequences that will play an indispensable role in our main proof. The first one concerns the time reversal of Brownian motion from its first hitting time.
\begin{lemma}[{cf. \cite[Corollary 1]{RP1981}}]\label{lem:timerev}
Let $\mu > 0$ and suppose $\tau_x= \inf\{t > 0: B_t^\mu = x\}$ is the first hitting time of $x>0$ by $(B_t^\mu)_{t \ge 0}$. Then
\begin{align}
(B_{\tau_x - t}^\mu)_{t \in [0, \tau_x]} 
\overset{d}{=} 
(x - \Ba_t^\mu)_{t \in [0, L_{x}]}
\end{align}

\noindent where  $L_{x} = \sup \{t > 0: \Ba_t^{\mu} = x\}$ is the last hitting time  of $(\Ba_t^\mu)_{t \ge 0}$.
\end{lemma}

The second one is the celebrated result of Williams on the path decomposition of Brownian motions with drift.\footnote{The case with unit drift $\mu = 1$ was stated in \cite[Theorem 2.2]{Wil1974} but this can be extended to arbitrary $\mu > 0$ by a suitable space--time rescaling
\begin{align*}
(B_t^1)_{t \ge 0} \overset{d}{=} \left( \mu B_{t / \mu^2}^\mu \right)_{t \ge 0}
\qquad \text{and} \qquad 
(\Ba_t^1)_{t \ge 0} \overset{d}{=} \left(\mu \Ba_{t / \mu^2}^\mu\right)_{t \ge 0},
\end{align*}

\noindent which can be verified by e.g. checking the infinitesimal generators of these processes.}
\begin{theorem}[cf. {\cite[Theorem 2.2]{Wil1974}}, {\cite[Corollary 2]{RP1981}}]\label{theo:path_dec}
Let $\mu > 0$ and consider the following independent objects:
\begin{itemize}
\item $(B_t^\mu)_{t \ge 0}$ is a Brownian motion with drift $\mu$.
\item $\Ma$ is an $\mathrm{Exp}(2\mu)$ variable.
\item $(\Ba_t^\mu)_{t \ge 0}$ is a Brownian motion with drift $\mu$ conditioned to stay non-negative.
\end{itemize}

Then the process $(R_t)_{t \ge 0}$ defined by
\begin{align}
R_t = \begin{cases}
B_t^\mu & t \le \tau_\Ma \\
\Ma - \Ba_{t - \tau_\Ma}^{\mu} & t \ge \tau_\Ma,
\end{cases}
\end{align}

\noindent where $\tau_x = \inf \{t > 0: B_t^\mu = x\}$ as usual, is a Brownian motion with drift $-\mu$.
\end{theorem}

\subsection{Gaussian multiplicative chaos} \label{subsec:GMC}
In this subsection we collect a few facts and fix some notations about multiplicative chaos. Before we begin, let us recall that $X(\cdot)$ is said to be a log-correlated Gaussian field on a bounded open domain $D \subset \RR^d$ if it is a (centred) Gaussian field on $D$ with covariance of the form
\begin{align*}
\EE[X(x) X(y)] = - \log |x-y| + f(x, y) \qquad \forall x, y \in D
\end{align*}

\noindent where (for the purpose of this paper) the function $f$ is continuous on $\overline{D} \times \overline{D}$. The logarithmic singularity on the diagonal of the covariance kernel means that $X(\cdot)$ is not defined pointwise and has to be interpreted as a generalised function: for any two test functions $\phi_1, \phi_2 \in C_c^\infty(D)$, the pair of random variables (formally written as)
\begin{align*}
\left( \int_D X(x) \phi_1(x)dx, \int_D X(x) \phi_2(x)dx\right)
\end{align*}

\noindent is jointly Gaussian with mean $0$ and covariance given by
\begin{align*}
\int_{D \times D} \phi_1(x) \underbrace{\Bigg[ - \log |x-y| + f(x, y)\Bigg]}_{=\EE[X(x)X(y)]} \phi_2(y) dxdy,
\end{align*}

\noindent and the set of test ``functions" may be extended to the set of finite (signed) measures 
\begin{align*}
\left\{\phi(d\cdot): \mathrm{supp}(\phi) \subset D \text{ and }~ \int_{D \times D} \EE[X(x)X(y)] \phi(dx) \phi(dy) < \infty \right\}.
\end{align*}

More concretely, the log-correlated field $X(\cdot)$ may be constructed via Karhunen--Lo\`{e}ve expansion and viewed as an abstract random variable in the negative Sobolev space $H^{-\epsilon}(\RR^d)$ for any $\epsilon > 0$; we refer the interested readers to \cite[Section 2.1]{JSW2020} for a self-contained proof of the Sobolev regularity.

\medskip
Given a log-correlated field $X(\cdot)$ on $D$ and $\gamma \in (0, \sqrt{2d})$, the associated Gaussian multiplicative chaos (GMC) is formally defined as the random measure
\begin{align*}
M_\gamma(dx) = e^{\gamma X(x) - \frac{\gamma^2}{2}\EE[X(x)^2]}dx.
\end{align*}

\noindent While the exponentiation of $X(\cdot)$ is not well-defined, the GMC measure may be constructed via various methods, such as martingale limits \cite{Kah1985}, smoothing procedures \cite{RV2010, Be2017}, and a generic approach based on the abstract language of Cameron-Martin space is also available \cite{Sha2016}. The main takeaway is that all these constructions lead to the same random measure and we are free to choose the most convenient one to work with. The following formulation is based on \cite{Be2017}.
\begin{theorem}\label{theo:GMCexist}
Let $X(\cdot)$ be a log-correlated Gaussian field on $D \subset \RR^d$ and $X_\epsilon := X \ast \nu_\epsilon$ where $\nu \in C_c^\infty(D)$ is a mollifier. For any fixed $\gamma \in (0, \sqrt{2d})$, the sequence of random measures
\begin{align*}
M_{\gamma, \epsilon}(dx) := e^{\gamma X_\epsilon(x) - \frac{\gamma^2}{2}\EE[X_\epsilon(x)^2]}dx, \qquad x \in D
\end{align*}

\noindent converges in probability, as $\epsilon \to 0^+$, to a non-trivial and non-atomic random measure $M_{\gamma}(dx)$ in the space of Radon measures on $D$ equipped with the weak$^*$ topology. Moreover, the limit $M_{\gamma}(dx)$ is independent of the choice of $\nu \in C_c^\infty(D)$.
\end{theorem}

We will need to deal with GMC measures associated to \emph{exact fields} and related objects; these are summarised in the definition below.
\begin{definition}\label{defi:exactppt}
Let $d \le 2$. By exact field, we mean the log-correlated Gaussian field $\overline{X}$ defined on the (closed) unit ball with covariance
\begin{align}\label{eq:exact}
\EE\left[\overline{X}(x) \overline{X}(y)\right] = -\log|x-y|
\qquad \forall x, y \in  \overline{B}(0,1).
\end{align}

\noindent The field $\overline{X}$ can be decomposed into two independent Gaussian components
\begin{align}
\EE\left[\overline{X}(x) \overline{X}(y)\right]
\label{eq:exact_cov}
& = -\log |x| \vee |y| + \log \frac{|x| \vee |y|}{|x-y|}\\
\notag
& =: \EE[B_{-\log |x|} B_{-\log|y|}] + \EE\left[\widehat{X}(x) \widehat{X}(y)\right]
\end{align}

\noindent where $(B_t)_{t \ge 0}$ is a standard Brownian motion. We also denote by\footnote{We have abused the notation in $d=2$ when we identify $\RR^2 \cong \CC$ and represent a point in the unit ball as a complex number.}
\begin{align}\label{eq:lateral}
Z_\gamma(dt) :=\begin{dcases}
e^{\gamma \widehat{X}(e^{-t}) - \frac{\gamma^2}{2} \EE[\widehat{X}(e^{-t})^2]} dt
+ e^{\gamma \widehat{X}(-e^{-t}) - \frac{\gamma^2}{2} \EE[\widehat{X}(-e^{-t})^2]} dt
& d = 1 \\
\left(\int_0^{2\pi} e^{\gamma \widehat{X}(e^{-t} e^{i\theta}) - \frac{\gamma^2}{2} \EE[\widehat{X}(e^{-t} e^{i\theta})^2]}d\theta\right)dt & d = 2
\end{dcases}
\end{align}

\noindent the GMC measure associated to $\widehat{X}$.
\end{definition}

\begin{remark}
We elaborate on a few things in \Cref{defi:exactppt} for those who may not be familiar with the theory of log-correlated fields and multiplicative chaos.
\begin{itemize}
\item The fact that the covariance kernel \eqref{eq:exact_cov} of $\overline{X}$ can be split into two independent parts follows from the radial-lateral decomposition of Gaussian free field in $d=2$ (the case $d=1$ may be seen as a one-dimensional restriction). Indeed one can construct the Brownian motion $B_{-\log |x|}$ by considering the circle averaging
\begin{align}
B_{-\log|x|} := \frac{1}{2\pi} \int_0^{2\pi} X(|x| e^{i\theta}) d\theta.
\end{align}

\noindent One can then define $\widehat{X}(x) := X(x) - B_{-\log|x|}$ and verify (using rotational symmetry) that
\begin{align*}
\EE\left[\widehat{X}(x) B_{-\log|y|}\right] = 0 \qquad \forall x, y \in \overline{B}(0, 1)
\end{align*}

\noindent which implies independence due to Gaussianity.

\item A priori, the field $\widehat{X}$ (often known as the lateral part of $\overline{X}$) is defined on $\overline{B}(0,1)$ in the sense of a generalised Gaussian field, but its scale invariance (see \Cref{lem:lateralinvariant} below) implies that it may be extended to the entire $\RR^2$ (on a suitable probability space).\footnote{This may also be obtained directly by performing a similar radial--lateral decomposition on a log-correlated field with covariance $-\log|x-y| +1_{\{|x| \ge 1\}} \log|x|+1_{\{|y| \ge 1\}} \log|y|$ on $\RR^2$, which coincides with that of the exact field when restricted to the unit ball. } Since $B_{-\log|\cdot|}$ has a continuous version (as a Brownian motion), $\widehat{X}$ inherits the regularity of $\overline{X}$ and may be seen as an abstract random variable living in a (local) Sobolev space of negative exponent.

\item $Z_{\gamma}(\cdot)$ (and its formal expression) should be interpreted in the sense of \Cref{theo:GMCexist}, i.e. the weak$^*$ limit (in probability) of
\begin{align*}
Z_{\gamma, \epsilon}(dt)
:=\begin{dcases}
e^{\gamma \widehat{X}_\epsilon(e^{-t}) - \frac{\gamma^2}{2} \EE[\widehat{X}_\epsilon(e^{-t})^2]} dt
+ e^{\gamma \widehat{X}_\epsilon(-e^{-t}) - \frac{\gamma^2}{2} \EE[\widehat{X}_\epsilon(-e^{-t})^2]} dt
& d = 1 \\
\left(\int_0^{2\pi} e^{\gamma \widehat{X}_\epsilon(e^{-t} e^{i\theta}) - \frac{\gamma^2}{2} \EE[\widehat{X}_\epsilon(e^{-t} e^{i\theta})^2]}d\theta\right)dt & d = 2
\end{dcases}
\end{align*}

\noindent where $\widehat{X}_{\epsilon} (\cdot):= (\widehat{X} \ast \nu_\epsilon)(\cdot)$, and again the limit is independent of the choice of mollifiers $\nu \in C_c^\infty(\RR^d)$. Similar to $\widehat{X}$, the domain of $Z_{\gamma}(d\cdot)$ can be extended to the entire real line $\RR$.
\end{itemize}
\end{remark}

For future reference, we record the invariance property of $\widehat{X}$ and $Z_{\gamma}$ in the following lemma.
\begin{lemma}\label{lem:lateralinvariant}
For any $c > 0$, we have
\begin{align*}
\left(\widehat{X}(cx)\right)_{x \in \RR^2}
\overset{d}{=} \left(\widehat{X}(x)\right)_{x \in \RR^2}
\qquad \text{and} \qquad
\left(Z_{\gamma} (dx)\right)_{x \in \RR}
\overset{d}{=}\left(Z_{\gamma}\circ \phi_c (dx)\right)_{x \in \RR}
\end{align*}

\noindent where $\phi_c: x \mapsto x+c$ is the shift operator.
\end{lemma}

\begin{proof}
The scale invariance of $\widehat{X}$ follows from a simple covariance check that\footnote{This can be made rigorous by integrating $\widehat{X}$ against test functions which we omit here.}
\begin{align*}
\EE\left[\widehat{X}(cx) \widehat{X}(cy)\right]
 = \log \frac{|cx| \vee |cy|}{|cx-cy|}
=\log \frac{|x| \vee |y|}{|x-y|}= \EE\left[\widehat{X}(x) \widehat{X}(y)\right].
\end{align*}

As for the shift invariance of $Z_{\gamma}$, let us start by considering $Z_{\gamma, \epsilon}$ associated to $\widehat{X}_{\epsilon} := \widehat{X} \ast \nu_\epsilon$ for some mollifier $\nu \in C_c^\infty(\RR^d)$. Since
\begin{align}
\notag
\left(\widehat{X}_{\epsilon}(e^{-c}x) \right)_{x \in \RR^d}
&= \left(\int \widehat{X}(e^{-c}x - \epsilon u) \nu(u) du\right)_{x \in \RR^d}\\
\label{eq:regularised_scale}
&\overset{d}{=} \left(\int \widehat{X}(x - e^c \epsilon u) \nu(u) du\right)_{x \in \RR^d} 
= \left(\widehat{X}_{e^c\epsilon}(x) \right)_{x \in \RR^d}
\end{align}

\noindent by the scale invariance of $\widehat{X}$, we see that
\begin{align}\label{eq:regularised_shift}
\left(Z_{\gamma, \epsilon}\circ \phi_c(dt)\right)_{t \in \RR}
\overset{d}{=}
\left(Z_{\gamma, e^c\epsilon}(dt)\right)_{t \in \RR}
\end{align}

\noindent and sending $\epsilon \to 0^+$ on both sides this implies $Z_{\gamma} \circ \phi_c \overset{d}{=} Z_{\gamma}$.
\end{proof}

Let us also recall from \cite[Appendix A]{Won2020} the various probabilistic representations of the reflection coefficient $\overline{C}_{\gamma, d}$.
\begin{lemma}\label{lem:reflection}
Let $d \le 2, \gamma \in (0, \sqrt{2d})$ and $p_c = 2d/\gamma^2$. The reflection coefficient $\overline{C}_{\gamma, d}$ appearing in \Cref{theo:main} is equal to
\begin{align}\label{eq:reflection1}
\overline{C}_{\gamma, d} = \lim_{t \to \infty} t^{p_c - 1} \PP\left(\int_{|x| \le r} |x|^{-\gamma^2} \overline{M}_{\gamma}(dx) > t \right)
\end{align}

\noindent where $\overline{M}_\gamma(dx) = e^{\gamma \overline{X}(x) - \frac{\gamma^2}{2}\EE[\overline{X}(x)^2]}dx$ is the Gaussian multiplicative chaos associated to the exact field $\overline{X}(\cdot)$, and this limit is independent of the value of $r \in (0, 1]$. In particular, it has the alternative probabilistic representation
\begin{align}\label{eq:reflection2}
\overline{C}_{\gamma, d}
= \EE\left[\left(\int_{-\infty}^\infty e^{-\gamma \Ba_t^\mu} Z_\gamma(dt)\right)^{p_c-1}\right]
\end{align}

\noindent with $\mu = \frac{\gamma}{2}(p_c-1) > 0$ where $(\Ba_t^\mu)_{t \ge 0}$ is independent of $Z_{\gamma}(\cdot)$.
\end{lemma}

\begin{remark}\label{rem:ref_indp}
The fact that the limit \eqref{eq:reflection1} is independent of the value $r > 0$ is a consequence of the fact that $\EE\left[\overline{M}_\gamma(B(0, 1))^{p_c-1}\right] < \infty$; see \cite[Theorem 1]{Won2020} for finer results regarding the tail asymptotics of subcritical GMCs.

As for the fact that the probabilistic representations in \Cref{lem:reflection} coincide with the expressions in \eqref{eq:reflection0}, this requires a long detour to the exact integrability of Liouville conformal field theory which is beyond the scope of this paper. The interested readers may consult \cite[Section 1.1 and 4.2]{RV2019} and \cite[Appendix A]{Won2020} for a brief discussion, and \cite{KRV2020, Re2020} for the derivations of these remarkable formulae. 
\end{remark}

\subsection{Some estimates}
In this subsection we collect some basic estimates. The first one concerns the construction of GMC measures via mollifications.
\begin{lemma}\label{lem:GMCbasic}
Using the notations in \Cref{theo:GMCexist}, the following are true:
\begin{itemize}
\item (cf. \cite[Lemma 3.5]{Be2017}) The covariance of $X_\epsilon$ is of the form
\begin{align*}
\EE[X_{\epsilon}(x)X_{\epsilon}(y)] = -\log \left(|x-y| \vee \epsilon \right) + f_\epsilon(x, y)
\end{align*}

\noindent where $\sup_{\epsilon > 0} \sup_{x, y \in D} |f_\epsilon(x, y)| < \infty$.

\item (cf. \cite[Proposition 3.5]{RV2010} or \cite[Section 3.7]{BPNotes}) Let $S \subset D$ be a bounded open subset and $\alpha \in (0, \frac{2d}{\gamma^2})$. Then the collection of random variables
 $(M_{\gamma, \epsilon}(S)^\alpha)_{\epsilon > 0}$ is uniformly integrable. In particular,
\begin{align*}
\EE\left[M_{\gamma, \epsilon}(S)^\alpha\right] \xrightarrow{\epsilon \to 0^+} \EE\left[M_{\gamma}(S)^\alpha\right] < \infty.
\end{align*}
\end{itemize}
\end{lemma}

We now state the multifractality of (regularised)  GMC measures.
\begin{lemma}\label{lem:GMC_multi}
Let $\gamma \in (0, \sqrt{2d})$, and consider a collection of continuous Gaussian fields\footnote{In this particular lemma we do not need $\{X_\epsilon\}_\epsilon$ to be constructed via mollification of one underlying field $X(\cdot)$ as long as the covariance assumption is satisfied.} $\{X_\epsilon(\cdot)\}_{\epsilon > 0}$ on $D$ with covariance of the form
\begin{align*}
\EE\left[X_\epsilon(x) X_\epsilon(y)\right] = -\log \left(|x-y| \vee \epsilon\right) + f_\epsilon(x, y) \qquad \forall x, y \in D
\end{align*}

\noindent where $C_f:= \sup_{\epsilon > 0} \sup_{x, y \in D} |f_\epsilon(x, y)| < \infty$. Suppose $D$ contains the ball $B(0, r)$ centred at origin with radius $r \in [\epsilon, 1/2]$. Then for any $\alpha \in [0, \frac{2d}{\gamma^2})$ and $\kappa \ge 0$, there exists some $C>0$ possibly depending on $C_f, \gamma, \alpha, \kappa$ but independent of $\epsilon, r$ such that the random measure 
\begin{align*}
M_{\gamma, \epsilon}(dx) := e^{\gamma X_\epsilon(x) - \frac{\gamma^2}{2} \EE[X_\epsilon(x)^2]}dx
\end{align*}

\noindent satisfies the moment estimates
\begin{align}
\label{eq:mult_mom}
\EE\left[M_{\gamma, \epsilon}(B(0, r))^\alpha\right]
& \le C r^{\alpha d+ \frac{\gamma^2}{2}\alpha (1-\alpha)},\\
\label{eq:mult_momsing}
\EE\left[\left(\int_{\epsilon \le |x| \le r} |x|^{-\kappa \gamma}M_{\gamma, \epsilon}(dx)\right)^\alpha\right]
&\le C\left(1 \vee \epsilon^{\alpha(d-\kappa \gamma) + \frac{\gamma^2}{2}\alpha(1-\alpha)}\right).
\end{align}
\end{lemma}

\begin{proof}
Both of these estimates are standard in the literature but for completeness we give a sketch of the claim for $\alpha \in (0, \frac{2d}{\gamma^2})$ (the case $\alpha = 0$ is trivial). For convenience, the constant $C > 0$ below may vary from line to line but its uniformity should be evident from the proof.

Let us start with \eqref{eq:mult_mom}: writing
\begin{align*}
M_{\gamma, \epsilon}(B(0, r))
\le M_{\gamma, \epsilon}(B(0, 2\epsilon)) + \sum_{n=0}^{\log \lfloor r/\epsilon \rfloor - 1} M_{\gamma, \epsilon} (A(2^{-n-1}r, 2^{-n}r))
\end{align*}

\noindent where
\begin{align}
\notag
&M_{\gamma, \epsilon} (A(2^{-n-1}r, 2^{-n}r))
:= \int_{2^{-n-1}r \le |x| \le 2^{-n}r} e^{\gamma X_\epsilon(x) - \frac{\gamma^2}{2}\EE[X_\epsilon(x)^2]}dx\\
\label{eq:slice}
&\quad = (2^{-n+1}r)^d \int_{\frac{1}{4}\le |u| \le \frac{1}{2}} e^{\gamma X_\epsilon(2^{-n+1}ru) - \frac{\gamma^2}{2}\EE[X_\epsilon(2^{-n+1}ru)^2]}du.
\end{align}

\noindent For $C>0$ sufficiently large (but independent of $\epsilon, r, n$), we have
\begin{align*}
&\EE\left[X_\epsilon(2^{-n+1}ru)X_\epsilon(2^{-n+1}rv)\right]\\
&\quad = -\log \left( |2^{-n+1}r (u-v)| \vee \epsilon\right) + f_\epsilon(2^{-n+1} ru, 2^{-n+1}rv)\\
&\quad = -\log(2^{-n+1}r) -\log \left( |u-v| \vee (2^{n-1}\epsilon/r)\right) + f_\epsilon(2^{-n+1} ru, 2^{-n+1}rv)\\
& \quad\begin{cases}
\le -\log(2^{-n+1}r)  + \EE\left[\overline{X}_\epsilon(u) \overline{X}_\epsilon(v)\right] + C\\
\ge -\log(2^{-n+1} r) - C
\end{cases} \qquad \forall u, v \in \left\{x: \frac{1}{4} \le |x| \le \frac{1}{2}\right\}
\end{align*}

\noindent where $\overline{X}_\epsilon = \overline{X} \ast \nu_\epsilon$ is a mollified exact field (see the first claim in \Cref{lem:GMCbasic}). If we write $\overline{M}_{\gamma, \epsilon}$ as the GMC measure associated to $\overline{X}_\epsilon$ and also interpret $-\log(2^{-n+1}r)$ as the variance of an independent Gaussian random variable $\Na_{n, r}$, then \Cref{cor:GP} implies that
\begin{align*}
& \EE\left[\left(\int_{\frac{1}{4}\le |u| \le \frac{1}{2}} e^{\gamma X_\epsilon(2^{-n+1}ru) - \frac{\gamma^2}{2}\EE[X_\epsilon(2^{-n+1}ru)^2]}du\right)^{\alpha}\right]\\
& \le \begin{cases}
C \EE\left[\left(\int_{\frac{1}{4}\le |u| \le \frac{1}{2}} e^{\gamma \Na_{n, r} - \frac{\gamma^2}{2}\EE[\Na_{n, r}^2]}du\right)^{\alpha}\right] & \text{for $\alpha \in (0, 1]$}\\
C \EE\left[\left(\int_{\frac{1}{4}\le |u| \le \frac{1}{2}} e^{\gamma (\overline{X}_\epsilon(u) + \Na_{n, r}) - \frac{\gamma^2}{2}\EE[(\overline{X}_\epsilon(u) + \Na_{n, r})^2]}du\right)^{\alpha}\right] &  \text{for $\alpha \in (1, \frac{2d}{\gamma^2})$} 
\end{cases}\\
& \le C \EE\left[\left(e^{\gamma \Na_{n, r} - \frac{\gamma^2}{2}\EE[\Na_{n, r}^2]}\right)^\alpha\right]  \times
\begin{cases}
1 & \text{for $\alpha \in (0, 1]$}\\
\EE\left[\overline{M}_{\gamma, \epsilon}(\overline{B}(0, \frac{1}{2}))^\alpha\right] & \text{for $\alpha \in (1, \frac{2d}{\gamma^2})$} 
\end{cases}\\
& \le C \left(2^{-n} r\right)^{\frac{\gamma^2}{2} \alpha(1-\alpha)}
\end{align*}

\noindent where the last inequality follows from the uniform integrability of $\overline{M}_{\gamma, \epsilon}(\overline{B}(0, \frac{1}{2}))^\alpha$ (see the second claim in \Cref{lem:GMCbasic}). Substituting this back to \eqref{eq:slice}, we obtain
\begin{align*}
\EE\left[M_{\gamma, \epsilon} (A(2^{-n-1}r, 2^{-n}r))^\alpha\right]
\le C \left(2^{-n}r\right)^{\alpha d + \frac{\gamma^2}{2}\alpha(1-\alpha)}.
\end{align*}

\noindent The same principle also shows that
\begin{align*}
\EE\left[M_{\gamma, \epsilon}(B(0, 2\epsilon))^\alpha\right]
\le C \left(2^{-n^*}r\right)^{\alpha d + \frac{\gamma^2}{2}\alpha(1-\alpha)} \qquad \text{where $n^* = \log \lfloor r/\epsilon\rfloor$}.
\end{align*}

Since $\alpha \in (0, \frac{2d}{\gamma^2})$ implies $\alpha d + \frac{\gamma^2}{2} \alpha(1-\alpha)  > 0$, the moment estimates we saw just now are summable: for $\alpha \in (0, 1]$ we have
\begin{align*}
& \EE\left[M_{\gamma, \epsilon}(B(0, r))^\alpha\right]\\
& \le \EE\left[M_{\gamma, \epsilon}(B(0, 2\epsilon))^\alpha\right] + \sum_{n=0}^{\log \lfloor r/\epsilon \rfloor - 1} \EE\left[M_{\gamma, \epsilon} (A(2^{-n-1}r, 2^{-n}r))^\alpha\right]\\
& \le C \sum_{n \ge 0} \left(2^{-n}r\right)^{\alpha d + \frac{\gamma^2}{2}\alpha(1-\alpha)}
\le C r^{\alpha d + \frac{\gamma^2}{2}\alpha(1-\alpha)}
\end{align*}

\noindent and the proof for $\alpha \in (1, \frac{2d}{\gamma^2})$ is similar by considering $ \EE\left[M_{\gamma, \epsilon}(B(0, r))^\alpha\right]^{1/\alpha}$ instead. Thus we have obtained \eqref{eq:mult_mom}.

Next we derive \eqref{eq:mult_momsing} using \eqref{eq:mult_mom}. For $\alpha  \in (0, 1]$, we have
\begin{align*}
\EE\left[\left(\int_{\epsilon \le |x| \le r} |x|^{-\kappa \gamma}M_{\gamma, \epsilon}(dx)\right)^\alpha\right]
&\le \sum_{n = 0}^{\lfloor \log r/\epsilon \rfloor} \EE\left[\left(\int_{2^n\epsilon \le |x| \le 2^{n+1}\epsilon} |x|^{-\kappa \gamma}M_{\gamma, \epsilon}(dx)\right)^\alpha\right]\\
&\le \sum_{n = 0}^{\lfloor \log r/\epsilon \rfloor} (2^n \epsilon)^{-\kappa \gamma \alpha} \EE\left[M_{\gamma, \epsilon}(B(0, 2^{n+1}\epsilon))^\alpha\right]\\
& \le C \sum_{n = 0}^{\lfloor \log r/\epsilon \rfloor} (2^n \epsilon)^{-\kappa \gamma \alpha} (2^{n+1}\epsilon)^{\alpha d + \frac{\gamma^2}{2}\alpha(1-\alpha)}\\
& \le C' \epsilon^{\alpha (d-\kappa \gamma) + \frac{\gamma^2}{2}\alpha(1-\alpha)} \frac{(r / \epsilon)^{\alpha (d-\kappa \gamma) + \frac{\gamma^2}{2}\alpha(1-\alpha)} - 1}{2^{\alpha (d-\kappa \gamma) + \frac{\gamma^2}{2}\alpha(1-\alpha)} - 1}
\end{align*}

\noindent which leads to the bound \eqref{eq:mult_momsing}. As for $\alpha > 1$ the inequality again follows from similar arguments applied to $\EE\left[\left(\int_{\epsilon \le |x| \le r} |x|^{-\kappa \gamma}M_{\gamma, \epsilon}(dx)\right)^\alpha\right]^{1/\alpha}$.
\end{proof}

The next lemma is an elementary estimate that allows us to remove irrelevant mass from our asymptotic analysis in \Cref{subsec:removal}.
\begin{lemma}\label{lem:basic}
Let $k > 0$, and $A, B$ be two non-negative random variables with finite $k$-th moment.
\begin{itemize}
\item[(i)] Suppose there exists some $\eta > 0$ such that $\EE\left[(A+B)^k\right] \le (1+\eta) \EE[B^k]$, then
\begin{align*}
\EE[A^k] \le \max \left(\eta, \left(\eta / k\right)^k\right) \EE[B^k].
\end{align*}
\item[(ii)] Suppose there exists some $\eta > 0$ such that $\EE[A^k] \le \eta \EE[B^k]$, then
\begin{align*}
\EE\left[(A+B)^k\right] \le \left[1 + \max\left(\eta, k2^k \eta^{\frac{1}{k}}\right)\right]\EE[B^k].
\end{align*}
\end{itemize}
\end{lemma}
\begin{proof}
We begin with the first statement. For $k \ge 1$, one can easily check (say by computing the first derivative) that
\begin{align*}
(A+B)^k - B^k - A^k \ge 0 \qquad \forall A, B \ge 0
\end{align*}

\noindent and so $\EE[A^k + B^k] \le (1+\eta) \EE[B^k]$ by assumption and hence $\EE[A^k] \le \eta \EE[B^k]$. As for $k \in (0, 1)$,
\begin{align*}
(A+B)^k - B^k
= k \int_B^{A+B} x^{k-1}dx \ge kAB^{k-1}
\end{align*}

\noindent and thus
\begin{align*}
\eta \EE[B^k]
\ge \EE[(A+B)^k] - \EE[B^k]
\ge k \EE[AB^{k-1}]
\ge k \EE[A^k]^{\frac{1}{k}} \EE[B^k]^{\frac{k-1}{k}}
\end{align*}

\noindent by (reverse) H\"older's inequality. Rearranging terms, we arrive at $\EE[A^k] \le (\eta / k)^k \EE[B^k]$.

Now for the second statement, we start with $k \in (0, 1)$ where
\begin{align*}
\EE[(A+B)^k] \le \EE[A^k] + \EE[B^k] \le (1+\eta) \EE[B^k]
\end{align*}

\noindent by our assumption. As for $k \ge 1$, we have
\begin{align*}
(A+B)^k - B^k 
= k \int_B^{A+B} x^{k-1} dx 
\le kA(A+B)^{k-1}
\le k2^{k-1} \left(A^k + AB^{k-1}\right)
\end{align*}

\noindent from which we obtain
\begin{align*}
\EE[(A+B)^k] 
& \le \EE[B^k] + k2^{k-1}\left(\EE[A^k] + \EE[AB^{k-1}]\right)\\
& \le\EE[B^k] + k2^{k-1}\EE[A^k] + k2^{k-1} \EE[A^k]^{\frac{1}{k}}\EE[B^{k}]^{\frac{k-1}{k}}\\
& \le \left[1 + k2^{k-1} \left(\eta + \eta^{\frac{1}{k}} \right)\right] \EE[B^k]
\end{align*}

\noindent by an application of H\"older's inequality. This concludes the proof.
\end{proof}

\section{Main proofs} \label{sec:mainproof}
For convenience, we shall write $p = p_c = 2d/\gamma^2$ throughout the entire \Cref{sec:mainproof}.
\subsection{Warm-up: a limit from below}
As advertised in \Cref{subsec:idea}, the proof of \Cref{theo:main} largely boils down to the analysis of positive moments of the exponential functional of Brownian motion with negative drift, with the help of Williams' path decomposition. 

To offer a first glance at how some of the results in \Cref{subsec:path_dec} may be applied, we make a digression and study the limit of a much simpler but related quantity. The following result, while not needed for the proof of \Cref{theo:main}, will be used when we sketch a renormalisation result in \Cref{app:lowercont}.

\begin{lemma}\label{lem:lower}
Let $\overline{X}$ be the exact field on the closed unit ball.  If $\gamma \in (0, \sqrt{2d})$ and $p = \frac{2d}{\gamma^2}$, then for any $r \in (0, 1]$,
\begin{align}
\notag
& \lim_{\alpha \to (p-1)^-} 
(p-1 - \alpha)\EE \left[\left(\int_{|x|\le r} e^{\gamma \overline{X}(x) - \frac{\gamma^2}{2} \EE[\overline{X}(x)^2]}\frac{dx}{|x|^{\gamma^2}} \right)^{\alpha} \right]\\
\label{eq:reflection}
& \qquad = (p-1) \EE\left[
\left(\int_{-\infty}^\infty e^{-\gamma \Ba_t^\mu} Z_\gamma(dt)\right)^{p-1}
\right]
\end{align}

\noindent where $Z_{\gamma}(dt)$ is the GMC associated to $\widehat{X}$ as defined in \eqref{eq:exact_cov} and \eqref{eq:lateral}.
\end{lemma}

The fact that \eqref{eq:reflection} is independent of $r \in (0, 1)$ follows from the same reasoning as in \Cref{rem:ref_indp}. Indeed for $0 < \alpha < p-1 \le 1$, we have
\begin{align*}
&\EE \left[\left(\int_{|x|\le r} e^{\gamma \overline{X}(x) - \frac{\gamma^2}{2} \EE[\overline{X}(x)^2]}\frac{dx}{|x|^{\gamma^2}} \right)^{\alpha} \right]\\
& \qquad \le \EE \left[\left(\int_{|x|\le 1} e^{\gamma \overline{X}(x) - \frac{\gamma^2}{2} \EE[\overline{X}(x)^2]}\frac{dx}{|x|^{\gamma^2}} \right)^{\alpha} \right]\\
&\qquad \le \EE \left[\left(\int_{|x|\le r} e^{\gamma \overline{X}(x) - \frac{\gamma^2}{2} \EE[\overline{X}(x)^2]}\frac{dx}{|x|^{\gamma^2}} \right)^{\alpha} \right] + r^{-\gamma^2 \alpha} \mathbb{E}\left[\overline{M}_{\gamma}(\overline{B}(0, 1))^\alpha\right]
\end{align*}

\noindent where $\overline{M}_{\gamma}$ is the GMC measure associated to $\overline{X}$. The last term on the RHS can be bounded uniformly in $\alpha \le p-1$ by the existence of GMC moments (\Cref{lem:GMCbasic}), and hence does not affect the limit in \eqref{eq:reflection}. The case $p - 1 > 1$ is similar. For this reason, we shall assume for simplicity that $r=1$ in the proof below.

\begin{proof}
Starting with the radial decomposition of our exact field
\begin{align*}
\EE\left[\overline{X}(x) \overline{X}(y)\right]
=: \EE[B_{-\log |x|} B_{-\log|y|}] + \EE\left[\widehat{X}(x) \widehat{X}(y)\right],
\end{align*}

\noindent we consider the change of variable $x = e^{-t}$ for $d=1$ and $x = e^{-t} e^{i\theta}$ for $d=2$ and rewrite the $\alpha$-th moment as 
\begin{align*}
\EE \left[\left(\int_{|x|\le 1} e^{\gamma \overline{X}(x) - \frac{\gamma^2}{2} \EE[\overline{X}(x)^2]}\frac{dx}{|x|^{\gamma^2}} \right)^{\alpha} \right]
= \EE\left[\left(\int_0^\infty e^{\gamma \left(B_t -\mu t\right)}Z_\gamma(dt)\right)^{\alpha}\right]
\end{align*}

\noindent with $\mu = \frac{\gamma}{2}(p-1)$. The next step is to re-express $(B_t - \mu t)_{t \ge 0}$ using \Cref{lem:timerev} and \Cref{theo:path_dec}. If $\phi_c: x \mapsto x + c$ is the shift operator, we have
\begin{align*}
\EE\left[\left(\int_0^\infty e^{\gamma \left(B_t -\mu t\right)}Z_\gamma(dt)\right)^{\alpha}\right]
&= \EE\left[\left(\int_{-L_{\Ma}}^\infty e^{\gamma (\Ma-\Ba_t^\mu)}Z_\gamma \circ \phi_{L_{\Ma}} (dt)\right)^{\alpha}\right]\\
&= \EE\left[\left(\int_{-L_{\Ma}}^\infty e^{\gamma (\Ma-\Ba_t^\mu)}Z_\gamma (dt)\right)^{\alpha}\right]
\end{align*}

\noindent where $(\Ba_{\pm t}^\mu)_{t \ge 0}$ are two independent Brownian motions with drift $\mu$ conditioned to stay positive, $\Ma$ is an independent $\mathrm{Exp}(2\mu)$ variable, and $L_{\Ma} := \sup \{t > 0: \Ba_{-t}^\mu = \Ma\}$. Note that the removal of $\phi_{\Ma}$ in the last equality above follows from the translation invariance of $Z_{\gamma}(dt)$, which is inherited from the scale invariance of $\widehat{X}(x)$ (see \Cref{lem:lateralinvariant}).

Now for any $x > 0$, we have
\begin{align*}
\EE\left[\left(\int_{-L_{x}}^\infty e^{\gamma (\Ma-\Ba_t^\mu)}Z_\gamma (dt)\right)^{\alpha} 1_{\{\Ma \ge x\}}\right]
& \le \EE\left[\left(\int_{-L_{\Ma}}^\infty e^{\gamma (\Ma-\Ba_t^\mu)}Z_\gamma (dt)\right)^{\alpha}\right]\\
& \le \EE\left[\left(\int_{-\infty}^\infty e^{\gamma (\Ma-\Ba_t^\mu)}Z_\gamma (dt)\right)^{\alpha}\right].
\end{align*}

\noindent It is enough to focus on the lower bound since along the way we will see that the upper bound is tight as $\alpha \to (p-1)^-$.  By independence, we compute
\begin{align*}
&\EE\left[\left(\int_{-L_{x}}^\infty e^{\gamma (\Ma-\Ba_t^\mu)}Z_\gamma (dt)\right)^{\alpha} 1_{\{\Ma \ge x\}}\right]\\
& \qquad = \EE\left[\left(\int_{-L_{x}}^\infty e^{-\gamma \Ba_t^\mu}Z_\gamma (dt)\right)^{\alpha} \right] \int_x^\infty e^{\gamma \alpha m} 2\mu e^{-2\mu m} dm\\
& \qquad = \EE\left[\left(\int_{-L_{x}}^\infty e^{-\gamma \Ba_t^\mu}Z_\gamma (dt)\right)^{\alpha} \right] \frac{2\mu}{2\mu - \gamma \alpha} e^{-(2\mu - \gamma \alpha) x}.
\end{align*}

Substituting $2 \mu = \gamma (p-1)$ back to our expression, we have
\begin{align*}
&\liminf_{\alpha \to (p-1)^-} (p-1 - \alpha)
\EE\left[\left(\int_{-L_{x}}^\infty e^{\gamma (\Ma-\Ba_t^\mu)}Z_\gamma (dt)\right)^{\alpha} 1_{\{\Ma \ge x\}}\right]\\
& \qquad \ge (p-1) \EE\left[\left(\int_{-L_{x}}^\infty e^{-\gamma \Ba_t^\mu}Z_\gamma (dt)\right)^{p-1} \right] .
\end{align*}

\noindent As $x > 0$ is arbitrary, we let $x \to \infty$ so that $L_{x} \to \infty$ almost surely and the lower bound matches the upper bound. This concludes the proof.
\end{proof}

\medskip

\subsection{Proof of \Cref{theo:main}}
Recall that $X_\epsilon = X \ast \nu_\epsilon$ for some mollifier $\nu \in C_c^\infty(\RR^d)$, and we shall assume without loss of generality that $\mathrm{supp}(\nu) \subset B(0, \frac{1}{4})$. The proof of our main theorem consists of three parts:
\begin{itemize}
\item[1.] We first remove part of the mass in $M_{\gamma, \epsilon}(D)$ which does not contribute to the blow-up.
\item[2.] We then argue that the problem may be reduced to one regarding the exact field \eqref{eq:exact}.
\item[3.] Finally, we prove that the $p$-th moment, after suitable renormalisation, converges to the probabilistic representation \eqref{eq:reflection2} of the reflection coefficient of GMC as $\epsilon \to 0^+$.
\end{itemize}

According to \Cref{lem:GMCbasic}, the covariance structure of $X_\epsilon$ is of the form
\begin{align}\label{eq:reg_cov}
\EE[X_\epsilon(x) X_\epsilon(y)] = -\log\left(|x-y| \vee \epsilon \right) + f_\epsilon(x, y)
\end{align}

\noindent where $\sup_{\epsilon > 0} \sup_{x, y \in D} |f_\epsilon(x, y)| < \infty$. In the following we shall assume without loss of generality that $f, f_\epsilon \ge 0$ everywhere in the domain. Indeed:
\begin{align*}
\EE[X(x) X(y)] &= -\log|x-y| + f(x, y) \\
& = -\log|r^{-1}x- r^{-1}y| - \log r + f(x, y)
\end{align*}

\noindent and since $-\log r + f(x, y) \ge 0$ uniformly for $r$ sufficiently small, we can always rewrite our expression in terms of another log-correlated field (by scaling our coordinate with $r$) with an `$f$'-function that satisfies the non-negativity condition, i.e.
$X^{(r)}(\cdot) := X(r \cdot)$ on the new domain $D/r$ with
\begin{align*}
\EE\left[X^{(r)}(x) X^{(r)}(y)\right] = \EE\left[X(rx) X(ry)\right] = -\log|x-y| + \underbrace{f(rx, ry) - \log r}_{=:f_r(x, y)}.
\end{align*}
\noindent The same argument also applies to $f_\epsilon$, and moreover we have
\begin{align*}
X_\epsilon(rx)
&:= \int X(rx - \epsilon u) \nu(u) du
= \int X^{(r)}(x - \tfrac{\epsilon}{r}u) \nu(u) du
= X^{(r)}_{\epsilon/ r}(u).
\end{align*}

The reduction to the case where $f, f_\epsilon \ge 0$ everywhere is only reasonable if \Cref{theo:main} is invariant under the rescaling of space, which we quickly verify now. If we rewrite
\begin{align*}
&\EE\left[\left(\int_D g(x) M_{\gamma, \epsilon}(dx)\right)^{p}\right]
= \EE\left[\left(\int_D g(x) e^{\gamma X_{\epsilon}(x) - \frac{\gamma^2}{2}\EE[X_\epsilon(x)^2]} dx\right)^{p}\right]\\
& \qquad = \EE\left[\left(\int_{D/r} r^dg(ru) e^{\gamma X_{\epsilon}(ru) - \frac{\gamma^2}{2}\EE[X_\epsilon(ru)^2]} du\right)^{p}\right]\\
& \qquad = \EE\left[\left(\int_{D/r} g_r(u) e^{\gamma X_{\epsilon/r}^{(r)}(u)- \frac{\gamma^2}{2}\EE[X_{\epsilon/r}^{(r)}(u)]} du\right)^{p}\right]
\end{align*}

\noindent where $g_r(u) := r^d g(ru)$ still satisfies the non-negativity and continuity assumptions, then \Cref{theo:main} leads to the asymptotics (as $\epsilon \to 0^+$)
\begin{align*}
&\left(\int_{D/r} e^{d(p-1) f_r(u, u)} g_r(u)^{p} du\right)\frac{\gamma^2}{2}(p-1)^2 \overline{C}_{\gamma, d}\log \frac{1}{\epsilon/r}\\
& =\left(\int_{D} e^{d(p-1) [f(x, x) - \log r]} (r^dg(x))^{p} d(x/r)\right)\frac{\gamma^2}{2}(p-1)^2 \overline{C}_{\gamma, d}\log \frac{1}{\epsilon/r}\\
& =\left(\int_{D} e^{d(p-1) f(x, x)]} g(x)^{p} dx\right)\frac{\gamma^2}{2}(p-1)^2 \overline{C}_{\gamma, d}\log \frac{r}{\epsilon}
\end{align*}

\noindent which is the desired result (since $r > 0$ is fixed and its appearance in the logarithm does not affect the leading order behaviour as $\epsilon \to 0$).

For simplicity our proof below treats the case where $g \equiv 1$ on $\overline{D}$, even though the analysis can easily be adapted to treat continuous $g \ge 0$ (see \Cref{rem:gextend} before \Cref{subsec:genapp}).
\subsubsection{Pre-processing: removal of irrelevant mass}\label{subsec:removal}

Let us recall that $p = 2d/\gamma^2 > 1$. Then
\begin{align*}
\EE\left[M_{\gamma, \epsilon}(D)^p\right]
& = \EE\left[\left(\int_D e^{\gamma X_{\epsilon}(x) - \frac{\gamma^2}{2} \EE[X_\epsilon(x)^2]} dx\right)^p\right]\\
& = \int_D du\EE\left[e^{\gamma X_{\epsilon}(u) - \frac{\gamma^2}{2} \EE[X_\epsilon(u)^2]}\left(\int_D e^{\gamma X_{\epsilon}(x) - \frac{\gamma^2}{2} \EE[X_\epsilon(x)^2]} dx\right)^{p-1}\right]\\
& = \int_D du\EE\left[\left(\int_D \frac{e^{\gamma^2 f_\epsilon(x, u)}}{\left(|x-u|\vee \epsilon\right)^{\gamma^2}} e^{\gamma X_{\epsilon}(x) - \frac{\gamma^2}{2} \EE[X_\epsilon(x)^2]} dx\right)^{p-1}\right]\\
\end{align*}

\noindent where the last equality follows from Cameron-Martin theorem (see \Cref{lem:CMG}) and the fact that
\begin{align*}
e^{\gamma^2 \EE[X_\epsilon(x) X_\epsilon(u)]} = \frac{e^{\gamma^2f_\epsilon(x, u)}}{(|x-u| \vee \epsilon)^{\gamma^2}}
\end{align*}

\noindent based on the expression \eqref{eq:reg_cov}. We shall therefore fix $u \in D$ and consider the integrand
\begin{align}\label{eq:integrand}
\EE\left[\left(\int_D \frac{e^{\gamma^2 f_\epsilon(x, u)}}{\left(|x-u|\vee \epsilon\right)^{\gamma^2}} e^{\gamma X_{\epsilon}(x) - \frac{\gamma^2}{2} \EE[X_\epsilon(x)^2]} dx\right)^{p-1}\right].
\end{align}

To proceed, we use the fact (which is a consequence of \Cref{lem:basic}) that we can remove any mass inside the integral in \eqref{eq:integrand} provided that they do not contribute to the leading order blow-up. The part we want to remove is
\begin{align*}
B_\epsilon := \int_{D \cap B(u, K\epsilon)} \frac{e^{\gamma^2 f_\epsilon(x, u)}}{\left(|x-u|\vee \epsilon\right)^{\gamma^2}} e^{\gamma X_{\epsilon}(x) - \frac{\gamma^2}{2} \EE[X_\epsilon(x)^2]} dx
\end{align*}

\noindent for some $K > 1$. For this we need to show that the $(p-1)$-th moment remains bounded when $\epsilon \to 0$, and this is true because
\begin{align*}
\EE \left[B_\epsilon^{p-1} \right]
&\le C\epsilon^{-\gamma^2(p-1)} \EE\left[M_{\gamma, \epsilon}(D \cap B(u, K \epsilon))^{p-1} \right]\\
&\le C \epsilon^{-\gamma^2(p-1)} (K\epsilon)^{d(p-1) + \frac{\gamma^2}{2}[(p-1)-(p-1)^2]}
\le C K^{\gamma^2(p-1)}
\end{align*}

\noindent where the first inequality follows from the fact that $f_\epsilon$ may be uniformly bounded (see \Cref{lem:GMCbasic}), and the second inequality from the first estimate in \Cref{lem:GMC_multi}.

Actually we will remove more mass from our analysis. Let $r > 0$ be some fixed number, then of course we have
\begin{align*}
&\EE \left[\left(\int_{D \cap \{|x-u| > r\}}  \frac{e^{\gamma^2 f_\epsilon(x, u)}}{\left(|x-u|\vee \epsilon\right)^{\gamma^2}} e^{\gamma X_{\epsilon}(x) - \frac{\gamma^2}{2} \EE[X_\epsilon(x)^2]} dx\right)^{p-1}\right]\\
& \qquad \le C r^{-(p-1)\gamma^2} \EE \left[\left(M_{\gamma, \epsilon}(D)\right)^{p-1}\right]
\end{align*}

\noindent where the RHS stays bounded as $\epsilon \to 0$, by the convergence and existence of $(p-1)$-th moment of Gaussian multiplicative chaos (see \Cref{lem:GMCbasic} and \Cref{lem:GMC_multi}). Summarising the work so far, we have
\begin{lemma} \label{lem:main1}
As $\epsilon \to 0^+$,
\begin{align*}
\EE\left[M_{\gamma, \epsilon}(D)^p\right]
&= \int_D du\EE\left[\left(\int_{D \cap \{K\epsilon \le |x-u| \le r\}} \frac{e^{\gamma^2 f_\epsilon(x, u)}}{|x-u|^{\gamma^2}} e^{\gamma X_{\epsilon}(x) - \frac{\gamma^2}{2} \EE[X_\epsilon(x)^2]} dx\right)^{p-1}\right]\\
& \qquad  + \mathcal{O}(1).
\end{align*}
\end{lemma}

\subsubsection{Reduction to exact fields}
The next step is to rewrite
\begin{align}
\EE\left[\left(\int_{D \cap \{K\epsilon \le |x-u| \le r\}} \frac{e^{\gamma^2 f_\epsilon(x, u)}}{|x-u|^{\gamma^2}} e^{\gamma X_{\epsilon}(x) - \frac{\gamma^2}{2} \EE[X_\epsilon(x)^2]} dx\right)^{p-1}\right]
\end{align}

\noindent in terms of something more analytically tractable. We claim that
\begin{lemma}\label{lem:reduced_exact}
There exist some constants $C_{K, r}(\epsilon) > 0$ independent of $u \in D$ such that
\begin{align*}
\limsup_{r \to 0} \limsup_{K \to \infty} \limsup_{\epsilon \to 0} C_{K, r}(\epsilon) = 0
\end{align*}

\noindent and that
\begin{align}
\notag
& \EE\left[\left(\int_{D \cap \{K\epsilon \le |x-u| \le r\}} \frac{e^{\gamma^2 f_\epsilon(x, u)}}{|x-u|^{\gamma^2}} e^{\gamma X_{\epsilon}(x) - \frac{\gamma^2}{2} \EE[X_\epsilon(x)^2]} dx\right)^{p-1}\right]\\
\notag
& \qquad = \left(1+ \mathcal{O}(C_{K, r}(\epsilon))\right) e^{d(p-1) f(u,u)}\\
\label{eq:reduced_exact}
& \qquad \qquad \times \EE\left[\left(\int_{(D- u)\cap \{K\epsilon \le |x| \le r\}} e^{\gamma \overline{X}_{\epsilon}(x)- \frac{\gamma^2}{2} \EE[\overline{X}_{\epsilon}(x)^2]}\frac{dx}{|x|^{\gamma^2}} \right)^{p-1}\right].
\end{align}
\end{lemma}

\begin{proof}
Recall our earlier reduction to the case where $f, f_\epsilon \ge 0$ everywhere, and without loss of generality suppose $r \le 1/4$. We consider a log-correlated random field $\overline{X}^{(u)}(\cdot)$ on $B(u, 4r)$ with covariance
\begin{align*}
\EE\left[\overline{X}^{(u)}(x) \overline{X}^{(u)}(y)\right]
= - \log |(x-u) - (y-u)| = -\log|x-y| \quad \forall x, y \in B(u, 4r).
\end{align*}

\noindent This field can be easily realised as the exact field in \eqref{eq:exact} up to a shift of domain, i.e. $\overline{X}^{(u)}(\cdot) := \overline{X}(\cdot - u)$. If we write $\overline{X}^{(u)}_{\epsilon} = \overline{X}^{(u)} \ast \nu_\epsilon$, then we have
\begin{align*}
\EE\left[\overline{X}^{(u)}_{\epsilon}(x) \overline{X}^{(u)}_{\epsilon} (y) \right]
&= \int \int \left(- \log|(x-u -\epsilon s) - (y-u - \epsilon t)|\right) \nu(s) \nu(t) ds dt\\
&= \int \int \left(- \log|(x -\epsilon s) - (y - \epsilon t)|\right) \nu(s) \nu(t) ds dt.
\end{align*}

\noindent Similarly, if $X_{\epsilon} = X \ast \nu_\epsilon$, then
\begin{align*}
\EE\left[X_\epsilon(x) X_\epsilon(y)\right]
= \int \int \left[- \log |(x-\epsilon s) - (y-\epsilon t)| + f(x - \epsilon s, y - \epsilon t) \right] \nu(s) \nu(t) ds dt.
\end{align*}
Since $f$ is continuous (and hence uniformly continuous) on the compact set $\overline{D} \times \overline{D}$, it has a modulus of continuity, say $\omega(\cdot)$. If $\Na_{u} \sim \Na(0, f(u,u))$ is independent of everything else, then
\begin{align*}
& \left|\EE\left[(\overline{X}^{(u)}_\epsilon(x)+\Na_u)(\overline{X}^{(u)}_\epsilon(y)+\Na_u)\right]
- \EE\left[X_\epsilon(x) X_\epsilon(y)\right]\right|\\
& =\left|\EE\left[\overline{X}^{(u)}_\epsilon(x)\overline{X}^{(u)}_\epsilon(y)\right] + f(u, u)
- \EE\left[X_\epsilon(x) X_\epsilon(y)\right]\right|\\
& = \left| \int \int [f(u, u) - f(x-\epsilon s, y - \epsilon t)] \nu(s) \nu(t) ds dt \right|
\le 2\omega(\epsilon) + 2\omega(r)
\end{align*}

\noindent for all $x, y \in D \cap \{v: K\epsilon \le |v-u| \le r\}$, and this upper bound goes to $0$ as $\epsilon, r \to 0^+$. Applying \Cref{cor:GP}, we can perform a Gaussian comparison and obtain the estimate
\begin{align*}
&  \frac{\EE\left[\left(\int_{D \cap \{K\epsilon \le |x-u| \le r\}} \frac{e^{\gamma^2 f_\epsilon(x, u)}}{|x-u|^{\gamma^2}} e^{\gamma (\overline{X}_{\epsilon}^{(u)}(x)+\Na_u) - \frac{\gamma^2}{2} \EE[(\overline{X}^{(u)}_{\epsilon}(x)+\Na_u)^2]} dx\right)^{p-1}\right]}{\EE\left[\left(\int_{D \cap \{K\epsilon \le |x-u| \le r\}} \frac{e^{\gamma^2 f_\epsilon(x, u)}}{|x-u|^{\gamma^2}} e^{\gamma X_{\epsilon}(x) - \frac{\gamma^2}{2} \EE[X_\epsilon(x)^2]} dx\right)^{p-1}\right]} \\
& \qquad = 1+ \mathcal{O}(\omega(\epsilon) + \omega(r))
\end{align*}

\noindent uniformly in $u \in D$, and so it suffices to consider the numerator
\begin{align}\label{eq:numerator}
&\EE\left[\left(\int_{D \cap \{K\epsilon \le |x-u| \le r\}} \frac{e^{\gamma^2 f_\epsilon(x, u)}}{|x-u|^{\gamma^2}} e^{\gamma (\overline{X}^{(u)}_{\epsilon}(x)+\Na_u) - \frac{\gamma^2}{2} \EE[(\overline{X}^{(u)}_{\epsilon}(x)+\Na_u)^2]}dx\right)^{p-1}\right].
\end{align}

We note the following observations:
\begin{itemize}
\item $\Na_u$ is independent of  $\overline{X}^{(u)}$ and its expectation may be evaluated directly.
\item We may replace $e^{\gamma^2 f_\epsilon(x, u)}$ by $e^{\gamma^2 f(u,u)}$ and pull it out of the expectation, up to some multiplicative error which vanishes as $\epsilon \to 0$, $K \to \infty$ and $r \to 0$. This is true because
\begin{align*}
|f_\epsilon(x, u) - f(u, u)| \le |f_\epsilon(x, u) - f(x, u)| + \omega(r),
\end{align*}

\noindent and for $|x-u| \ge K\epsilon \ge \epsilon$:
\begin{align*}
&\left|f_\epsilon(x, u) - f(x, u)\right|
=\left| \EE[X_\epsilon(x) X_\epsilon(u)] - \EE[X(x) X(u)]\right|\\
&= \int \int \left[-\log\frac{|(x-\epsilon s) - (u - \epsilon t)|}{|x-u|} +f(x - \epsilon s, u - \epsilon t) - f(x,u)\right]\nu(s) \nu(t) ds dt\\
& \le \int \int\left| \log \left(1 - \epsilon \frac{s-t}{|x-u|}\right) \right| \nu(s)\nu(t) ds dt + 2\omega(\epsilon)
\le \frac{1}{K} + 2 \omega(\epsilon)
\end{align*}

\noindent where the last inequality follows from the assumption that $\mathrm{supp}(\nu) \subset B(0, \frac{1}{4})$.
\end{itemize}

\noindent Taking all these considerations together, \eqref{eq:numerator} can be rewritten (up to a multiplicative error of $C_{K, r}(\epsilon)$ as defined in the statement of the lemma) as
\begin{align*}
\notag & e^{\gamma^2(p-1) f(u,u)}\EE\left[\left(e^{\gamma \Na_u - \frac{\gamma^2}{2} \EE[\Na_u^2]}\right)^{p-1}\right]\\
\notag &\qquad \qquad \times \EE\left[\left(\int_{D \cap \{K\epsilon \le |x-u| \le r\}} \frac{e^{\gamma \overline{X}^{(u)}_{\epsilon}(x)- \frac{\gamma^2}{2} \EE[\overline{X}^{(u)}_{\epsilon}(x)^2]}dx}{|x-u|^{\gamma^2}} \right)^{p-1}\right]\\
& =  e^{d(p-1) f(u,u)}\EE\left[\left(\int_{(D-u) \cap \{K\epsilon \le |x| \le r\}} e^{\gamma \overline{X}_{\epsilon}(x)- \frac{\gamma^2}{2} \EE[\overline{X}_{\epsilon}(x)^2]}\frac{dx}{|x|^{\gamma^2}} \right)^{p-1}\right]
\end{align*}

\noindent since $\overline{X}_\epsilon^{(u)}(\cdot) = \overline{X}^{(u)} \ast \nu_\epsilon = \overline{X}(\cdot - u) \ast \nu_\epsilon = \overline{X}_\epsilon(\cdot - u)$.
\end{proof}

\subsubsection{Convergence to the reflection coefficient}\label{subsec:conref}
The final ingredient we need is the following asymptotic formula.
\begin{lemma}\label{lem:toreflection}
We have
\begin{align} \notag
&\EE\left[\left(\int_{K\epsilon \le |x| \le r} e^{\gamma \overline{X}_{\epsilon}(x)- \frac{\gamma^2}{2} \EE[\overline{X}_{\epsilon}(x)^2]}\frac{dx}{|x|^{\gamma^2}} \right)^{p-1}\right]\\
\label{eq:toreflection}
& \qquad \overset{\epsilon \to 0^+}{\sim} 
\frac{\gamma^2}{2}(p-1)^2 \EE\left[
\left(\int_{-\infty}^\infty e^{-\gamma \Ba_t^\mu} Z_\gamma(dt)\right)^{p-1}\right] \log \frac{1}{\epsilon}.
\end{align}
\end{lemma}

\begin{proof}[Proof of \Cref{theo:main}]
Note that for any $u \in D$, we have the trivial uniform bound
\begin{align*} 
& \EE\left[\left(\int_{(D- u)\cap \{K\epsilon \le |x| \le r\}} e^{\gamma \overline{X}_{\epsilon}(x)- \frac{\gamma^2}{2} \EE[\overline{X}_{\epsilon}(x)^2]}\frac{dx}{|x|^{\gamma^2}} \right)^{p-1}\right]\\
& \qquad \le 
\EE\left[\left(\int_{\{K\epsilon \le |x| \le r\}} e^{\gamma \overline{X}_{\epsilon}(x)- \frac{\gamma^2}{2} \EE[\overline{X}_{\epsilon}(x)^2]}\frac{dx}{|x|^{\gamma^2}} \right)^{p-1}\right]
\end{align*}

\noindent the asymptotics of which is given by \Cref{lem:toreflection}. On the other hand, since $D$ is open, there exists $\tilde{r} = \tilde{r}(u) > 0$ sufficiently small (but not uniform in $u \in D$) such that $(D- u)\cap \{K\epsilon \le |x| \le \tilde{r}\} = \{K\epsilon \le |x| \le \tilde{r}\}$. By dominated convergence, this means that
\begin{align*}
&\lim_{\epsilon \to 0^+} \frac{1}{-\log \epsilon} \int_D du \Bigg\{\EE\left[\left(\int_{\{K\epsilon \le |x| \le r\}} e^{\gamma \overline{X}_{\epsilon}(x)- \frac{\gamma^2}{2} \EE[\overline{X}_{\epsilon}(x)^2]}\frac{dx}{|x|^{\gamma^2}} \right)^{p-1}\right]\\
& \quad -\EE\left[\left(\int_{(D- u)\cap \{K\epsilon \le |x| \le r\}} e^{\gamma \overline{X}_{\epsilon}(x)- \frac{\gamma^2}{2} \EE[\overline{X}_{\epsilon}(x)^2]}\frac{dx}{|x|^{\gamma^2}} \right)^{p-1}\right] \Bigg\} = 0.
\end{align*}

Combining this with \Cref{lem:main1}, \Cref{lem:reduced_exact} and again \Cref{lem:toreflection}, we see that both
\begin{align*}
\liminf_{\epsilon \to 0^+} \frac{1}{-\log \epsilon} \EE\left[M_{\gamma, \epsilon}(D)^p\right]
\qquad \text{and} \qquad
\limsup_{\epsilon \to 0^+} \frac{1}{-\log \epsilon} \EE\left[M_{\gamma, \epsilon}(D)^p\right]
\end{align*}

\noindent are equal to
\begin{align*}
&  \left(1+ \mathcal{O}\left(\limsup_{\epsilon \to 0^+} C_{K, r}(\epsilon)\right)\right) \left(\int_D e^{d(p-1) f(u,u)} du\right)\\
& \qquad \qquad \times \frac{\gamma^2}{2}(p-1)^2 \EE\left[
\left(\int_{-\infty}^\infty e^{-\gamma \Ba_t^\mu} Z_\gamma(dt)\right)^{p-1}\right].
\end{align*}

\noindent Since $K, r > 0$ are arbitrary, we send $K \to \infty$ and $r \to 0$, and conclude our proof by identifying
\begin{align*}
 \EE\left[\left(\int_{-\infty}^\infty e^{-\gamma \Ba_t^\mu} Z_\gamma(dt)\right)^{p-1}\right]
\end{align*}

\noindent as the probabilistic representation of $\overline{C}_{\gamma, d}$ with \Cref{lem:reflection}.
\end{proof}

The rest of the section is to prove \Cref{lem:toreflection}, which is broken into several steps.
\paragraph{Step 1: simplifying the radial part.}
As in the proof of \Cref{lem:lower}, we start by considering the radial decomposition of the exact field $\overline{X}$ and extend it to its regularised version, i.e.
\begin{align*}
\overline{X}_{\epsilon}(x) = B_{\epsilon, -\log |x|} + \widehat{X}_\epsilon(x)
\end{align*}

\noindent where $B_{\epsilon, -\log |x|} =( B_{-\log|\cdot|} \ast \nu_\epsilon)(x)$ and $\widehat{X}_{\epsilon}(x) = \widehat{X} \ast \nu_{\epsilon}(x)$. By performing the substitution $x = e^{-t}$ when $d=1$ or $x=e^{-t} e^{i\theta}$ when $d=2$, we obtain
\begin{align*}
&\EE\left[\left(\int_{K\epsilon \le |x| \le r} e^{\gamma \overline{X}_{\epsilon}(x)- \frac{\gamma^2}{2} \EE[\overline{X}_{\epsilon}(x)^2]}\frac{dx}{|x|^{\gamma^2}} \right)^{p-1}\right]\\
& \qquad = \EE \left[ \left( \int_{-\log r}^T e^{\gamma B_{\epsilon, t} - \frac{\gamma^2}{2} \EE[B_{\epsilon, t}^2]} e^{(d+\gamma^2)t}Z_{\gamma, \epsilon}(dt)\right)^{p-1}\right]
\end{align*}

\noindent for $T = -\log K \epsilon$, and $Z_{\gamma, \epsilon}(dt)$ is defined in the obvious way. If $-\log r \le s \le t \le T$,
\begin{align*}
\EE\left[B_{\epsilon, s} B_{\epsilon, t}\right]
& = \int \int \nu(u) \nu(v) du dv\EE\left[B_{-\log|e^{-s}e_1 + \epsilon v|}B_{-\log|e^{-t}e_1 + \epsilon v|}\right]\\
& \in [s - \log|1 + \epsilon e^{s}|, s - \log |1 - \epsilon e^{s}|]
\end{align*}

\noindent where $\epsilon e^{s} \le K^{-1}$. This means that
\begin{align}\label{eq:BM_approx}
\left|\EE\left[B_s B_t\right] - \EE\left[B_{\epsilon, s} B_{\epsilon, t}\right]\right|
\le \frac{2}{K} \qquad \forall s, t \in [-\log r, T]
\end{align}

\noindent for $K>1$ sufficiently large, and by \Cref{cor:GP} we obtain
\begin{align*}
e^{-\frac{\gamma^2}{K}|(p-1)(p-2)|}
\le \frac{\EE \left[ \left( \int_{-\log r}^T e^{\gamma B_t - \frac{\gamma^2}{2} \EE[B_t^2]} e^{(d+\gamma^2)t}Z_{\gamma, \epsilon}(dt)\right)^{p-1}\right]}{\EE \left[ \left( \int_{-\log r}^T e^{\gamma B_{\epsilon, t} - \frac{\gamma^2}{2} \EE[B_{\epsilon, t}^2]} e^{(d+\gamma^2)t}Z_{\gamma, \epsilon}(dt)\right)^{p-1}\right]}
\le e^{\frac{\gamma^2}{K}|(p-1)(p-2)|}.
\end{align*}
 
Therefore it suffices to consider
\begin{align*}
\EE \left[ \left( \int_{-\log r}^T e^{\gamma B_t - \frac{\gamma^2}{2} \EE[B_t^2]} e^{(d+\gamma^2)t}Z_{\gamma, \epsilon}(dt)\right)^{p-1}\right]
=\EE \left[ \left( \int_{-\log r}^T e^{\gamma (B_t - \mu t)} Z_{\gamma, \epsilon}(dt)\right)^{p-1}\right]
\end{align*}

\noindent where $\mu = \frac{\gamma}{2}(p-1)$. By \Cref{theo:path_dec}, we can decompose the path of $(B_t - \mu t)_{t \ge 0}$ and rewrite our expectation as
\begin{align}
\notag
&\EE \left[ \EE \left[ \left( \int_{-\log r}^T e^{\gamma (B_t - \mu t)} Z_{\gamma, \epsilon}(dt)\right)^{p-1} \Bigg| \max_{t \ge 0} (B_t - \mu t)\right]\right]\\
\notag
& = \EE\left[\EE \left[ e^{\gamma(p-1)\Ma} \left( \int_{-\log r}^{T} e^{\gamma (B_t^\mu -\Ma) 1_{\{t \le \tau_\Ma\}} - \gamma \Ba_{t - \tau_\Ma}^\mu 1_{\{t \ge \tau_\Ma\}}}  Z_{\gamma, \epsilon} (dt)\right)^{p-1} \Bigg|  \Ma \right]\right]\\
\label{eq:expfunc0}
& = 2\mu \int_0^{\infty} dm \EE \left[\left( \int_{-\log r}^{T\wedge \tau_m} e^{\gamma (B_t^\mu -m)} Z_{\gamma, \epsilon} (dt) + \int_0^{(T - \tau_m)_+} e^{-\gamma \Ba_t^\mu} Z_{\gamma, \epsilon} \circ \phi_{\tau_m}(dt)\right)^{p-1}  \right]
\end{align}

\noindent where $\phi_c: x \mapsto x+c$ is the right shift and $
\tau_m:= \inf \{t > 0: B_t^\mu= m\}$. Using \Cref{lem:timerev}, the integrand may also be written as
\begin{align}
\label{eq:expfunc}
 \EE \left[ \left( \int_{-L_{m} - \log r}^{(T-L_{m})_+} e^{-\gamma \Ba_t^\mu} Z_{\gamma, \epsilon}\circ \phi_{L_{m}} (dt)\right)^{p-1}\right]
\end{align}

\noindent with $L_{m} := \sup \{t > 0: \Ba_{-t}^\mu = m\}$.

\paragraph{Step 2: a uniform upper bound.} 
We claim that the integrand in \eqref{eq:expfunc0} is uniformly bounded in $m$ and $\epsilon$ (for any fixed $K>1$ and $r \in (0, \frac{1}{4})$). It is more convenient to deal with the representation \eqref{eq:expfunc}, and we will show that
\begin{align}\label{eq:ub_integrand}
 \EE \left[ \left( \int_{-\infty}^\infty e^{-\gamma \Ba_t^\mu} Z_{\gamma, \epsilon}\circ \phi_{L_{m}} (dt)\right)^{p-1}\right] \le C_p
\end{align}

\noindent for some constants $C_p \in (0, \infty)$ independent of $m$ and $\epsilon$. For future reference, we will actually show a slightly strengthened version of \eqref{eq:ub_integrand}, namely: for any $c \ge 0$, 
\begin{align}\label{eq:ub_integrand2}
 \EE \left[ \left( \int_{|t| \ge c} e^{-\gamma \Ba_t^\mu} Z_{\gamma, \epsilon}\circ \phi_{L_{m}} (dt)\right)^{p-1}\right] \le C_p(c)
\end{align}

\noindent for some $C_p(c) \in (0, \infty)$ independent of $m$ and $\epsilon$, such that $C_p(c)$ is monotonically decreasing to $0$ as $c \to \infty$. There are two cases to consider:
\begin{itemize}
\item $p-1 \ge 1$, i.e. the map $x \mapsto x^{p-1}$ is convex.\\
First observe the distributional equality $Z_{\gamma, \epsilon} \circ \phi_c(dt) \overset{d}{=} Z_{\gamma, \epsilon e^c}(dt)$ for any $c \ge 0$ (possibly random but independent of $\widehat{X}$ (see \eqref{eq:regularised_scale} and \eqref{eq:regularised_shift} in the proof of \Cref{lem:lateralinvariant} for the details). Therefore, our integrand may be bounded by first rewriting it as
\begin{align*}
\EE\left[\EE \left[ \left( \int_{|t| \ge c} e^{-\gamma \Ba_t^\mu} Z_{\gamma, \epsilon \exp(L_{m})} (dt)\right)^{p-1} \Bigg| (\Ba_t^\mu)_{t \in \RR} \right]\right]
\end{align*}

\noindent and then performing Gaussian comparison in the conditional expectation. Since
\begin{align*}
\EE\left[\widehat{X}_{\epsilon}(x) \widehat{X}_\epsilon(y)\right]
\le \EE \left[\widehat{X}(x) \widehat{X}(y)\right] + C
\end{align*}

\noindent for some $C > 0$ independent of $\epsilon > 0$, \Cref{cor:GP} implies that 
\begin{align*}
&\EE\left[\left( \int_{|t| \ge c} e^{-\gamma \Ba_t^\mu} Z_{\gamma, \epsilon \exp(L_{m})} (dt)\right)^{p-1} \right]\\
& \qquad \le  e^{\frac{\gamma^2}{2}p(p-1) C}\EE\left[\left( \int_{|t| \ge c} e^{-\gamma \Ba_t^\mu} Z_{\gamma} (dt)\right)^{p-1} \right].
\end{align*}

\noindent When $c = 0$, the expectation appearing in the last bound is the probabilistic representation of $\overline{C}_{\gamma, d}$ (see \Cref{lem:reflection}) and is thus finite. If we take the above bound as the definition of $C_p(c)$, then $C_p(c)$ satisfies all the desired properties (uniformity in $m, \epsilon$ and monotone convergence to $0$ as $c \to \infty$) in our claim.

\item $p-1 \in (0, 1)$, i.e. the map $x \mapsto x^{p-1}$ is concave.\\
By Jensen's inequality and subadditivity, we have
\begin{align*}
\EE\left[\left( \int_{|t| \ge c} e^{-\gamma \Ba_t^\mu} Z_{\gamma, \epsilon \exp(L_{m})} (dt)\right)^{p-1} \right]
& \le \EE\left[\int_{|t| \ge c} e^{-\gamma \Ba_t^\mu} Z_{\gamma, \epsilon \exp(L_{m})} (dt)\right]^{p-1}\\
& \le  4\pi\EE\left[\int_{c}^\infty e^{-\gamma \Ba_t^\mu} dt\right]^{p-1}.
\end{align*}

\noindent Again we can take the above bound as the definition of $C_p(c)$, which obviously satisfies all the desired properties except the claim that $C_p(0) < \infty$ which we verify now. By \Cref{lem:stocdom}, $(\Ba_t^\mu)_{t \ge 0}$ stochastically dominates a $\mathrm{BES}_0(3)$-process, which we denote by $(\beta_t)_{t \ge 0}$. In particular, using the fact that $\beta_t \overset{d}{=} \sqrt{t}\beta_1$ for any fixed $t > 0$, we have
\begin{align*}
\EE\left[\int_{0}^\infty e^{-\gamma \Ba_t^\mu} dt\right]
& = \int_{0}^\infty \EE\left[e^{-\gamma \Ba_t^\mu}\right] dt\\
& \le \int_{0}^\infty \EE\left[e^{-\gamma \beta_t}\right] dt\\
&= \EE\left[\int_{0}^\infty e^{-\gamma \sqrt{t}\beta_1} dt\right]
\overset{u = \beta_1^2 t}{=} \EE\left[\beta_1^{-2}\right] \int_0^\infty e^{-\gamma \sqrt{u}} du < \infty
\end{align*}

\noindent and so we are done.
\end{itemize}

To analyse the integrand further,  we fix some small $\delta \in (0, 1/4)$ and proceed by splitting the $m$-integral into two regions.
\paragraph{Step 3: main contribution from $m \le \mu(1-2\delta)T$.}
Let us fix some $m_0 > 0$ and ignore all the contributions to the integral from $m < m_0$. We recall from \Cref{lem:timerev} that the law of $L_{m}$ is the same as that of $\tau_m := \inf \{t > 0: B_t^\mu = m\}$ and so the event $A_\delta := \{ L_{m} \le (1-\delta)T\}$ satisfies
\begin{align*}
\PP(A_\delta) 
\ge \PP\left(\max_{t \le T} B_t^\mu \ge \mu(1-\delta)T\right)
\xrightarrow{T \to \infty} 1
\end{align*}

\noindent for any $m \le \mu(1-\delta)T$. Then
\begin{align}
\notag
&  \EE \left[ \left( \int_{-L_{m} - \log r}^{(T-L_{m})_+} e^{-\gamma \Ba_t^\mu} Z_{\gamma, \epsilon}\circ \phi_{L_{m}} (dt)\right)^{p-1}\right]\\
\notag
& \qquad =  \EE \left[ \left( \int_{-L_{m} - \log r}^{(T-L_{m})_+} e^{-\gamma \Ba_t^\mu} Z_{\gamma, \epsilon\exp(L_{m})} (dt)\right)^{p-1}\right]\\
\notag
& \qquad \ge  \EE \left[ \left( \int_{-L_{m_0} - \log r}^{\delta T} e^{-\gamma \Ba_t^\mu} Z_{\gamma, \epsilon\exp(L_{m})} (dt)\right)^{p-1} 1_{A_\delta}\right]\\
\notag
& \qquad \ge  \EE \left[ \left( \int_{-L_{m_0} - \log r}^{M} e^{-\gamma \Ba_t^\mu} Z_{\gamma, \epsilon\exp(L_{m})} (dt)\right)^{p-1} \right]\\
\label{eq:Step3lb1}
& \qquad \qquad - \EE \left[ \left( \int_{-L_{m_0} - \log r}^{M} e^{-\gamma \Ba_t^\mu} Z_{\gamma, \epsilon\exp(L_{m})} (dt)\right)^{p-1} 1_{A_\delta^c}\right]
\end{align}

\noindent where the last inequality just follows from the truncation of the integral at a deterministic (but arbitrary) level $M \le \delta T$ (which is allowed since $T = -\log K \epsilon$ is arbitrarily large as $\epsilon \to 0^+$). 

Since $L_m$ is almost surely finite and only depends on $(\Ba_t^\mu)_{t \in \RR}$ (and in particular independent of $\widehat{X}$), we have $Z_{\gamma, \epsilon\exp(L_m)} \xrightarrow{\epsilon \to 0^+} Z_{\gamma}$ in probability in the weak$^*$-topology of measures on any compact intervals of $\RR$. Combining this with the fact that the $(p-1)$-th power of the total mass of regularised GMC measures are uniformly integrable and that their expectation converges to the $(p-1)$-th moment of limiting GMC measure (see \Cref{lem:GMCbasic}), we see that
\begin{align*}
& \EE \left[ \left( \int_{-L_{m_0} - \log r}^{M} e^{-\gamma \Ba_t^\mu} Z_{\gamma, \epsilon\exp(L_{m})} (dt)\right)^{p-1} \right]\\
& \qquad \xrightarrow{\epsilon \to 0^+}
\EE \left[ \left( \int_{-L_{m_0} - \log r}^{M} e^{-\gamma \Ba_t^\mu} Z_{\gamma} (dt)\right)^{p-1} \right].
\end{align*}

\noindent Similarly, since $\PP(A_\delta^c) \xrightarrow{\epsilon \to 0^+} 0$, the uniform integrability of the $(p-1)$-th moments of the regularised GMCs implies
\begin{align*}
\EE \left[ \left( \int_{-L_{m_0} - \log r}^{M} e^{-\gamma \Ba_t^\mu} Z_{\gamma, \epsilon\exp(L_{m})} (dt)\right)^{p-1} 1_{A_\delta^c}\right] \xrightarrow{\epsilon \to 0^+} 0.
\end{align*}

\noindent Since $M> 0$ is arbitrary, we obtain from \eqref{eq:Step3lb1} the lower bound
\begin{align}
\notag
&\liminf_{\epsilon \to 0^+} \EE \left[ \left( \int_{-L_{m} - \log r}^{(T-L_{m})_+} e^{-\gamma \Ba_t^\mu} Z_{\gamma, \epsilon}\circ \phi_{L_{m}} (dt)\right)^{p-1}\right]\\
\label{eq:Step3lb}
& \qquad \qquad \ge \EE \left[ \left( \int_{-L_{m_0} - \log r}^{\infty} e^{-\gamma \Ba_t^\mu} Z_{\gamma} (dt)\right)^{p-1} \right].
\end{align}

\noindent As for the upper bound, consider
\begin{align}
\notag
& \EE \left[ \left( \int_{-L_{m} - \log r}^{(T-L_{m})_+} e^{-\gamma \Ba_t^\mu} Z_{\gamma, \epsilon}\circ \phi_{L_{m}} (dt)\right)^{p-1}\right]\\
\notag
& \qquad \qquad \le  \EE \left[ \left( \int_{-\infty}^\infty e^{-\gamma \Ba_t^\mu} Z_{\gamma, \epsilon}\circ \phi_{L_{m}} (dt)\right)^{p-1} \right]\\
\label{eq:Step3ub1}
& \qquad \qquad  =  \EE \left[ \left( \int_{-\infty}^\infty e^{-\gamma \Ba_t^\mu} Z_{\gamma, \epsilon \exp(L_{m})} (dt)\right)^{p-1} \right]
\end{align}

\noindent in the following two cases:
\begin{itemize}
\item $p-1 \ge 1$. We can upper bound the $(p-1)$-th root of \eqref{eq:Step3ub1} by
\begin{align*}
& \EE \left[ \left( \int_{|t| \le M} e^{-\gamma \Ba_t^\mu} Z_{\gamma, \epsilon \exp(L_{m})} (dt)\right)^{p-1} \right]^{\frac{1}{p-1}}\\
& \qquad +  \EE \left[ \left( \int_{|t| > M} e^{-\gamma \Ba_t^\mu} Z_{\gamma, \epsilon \exp(L_{m})} (dt)\right)^{p-1} \right]^{\frac{1}{p-1}}
\end{align*}

\noindent for some fixed $M > 0$.
\item $p \in (0, 1)$. We can upper bound \eqref{eq:Step3ub1} by
\begin{align*}
& \EE \left[ \left( \int_{|t| \le M} e^{-\gamma \Ba_t^\mu} Z_{\gamma, \epsilon \exp(L_{m})} (dt)\right)^{p-1} \right]\\
& \qquad +  \EE \left[ \left( \int_{|t| > M} e^{-\gamma \Ba_t^\mu} Z_{\gamma, \epsilon \exp(L_{m})} (dt)\right)^{p-1} \right]
\end{align*}

\noindent for some fixed $M > 0$.

\end{itemize}
Under either scenario, the first term provides the main contribution, and the same argument from the lower bound leads to
\begin{align*}
&\EE \left[ \left( \int_{|t| \le M} e^{-\gamma \Ba_t^\mu} Z_{\gamma, \epsilon \exp(L_{m})} (dt)\right)^{p-1} \right]\\
& \xrightarrow{\epsilon \to 0^+} \EE \left[ \left( \int_{|t| \le M} e^{-\gamma \Ba_t^\mu} Z_{\gamma} (dt)\right)^{p-1} \right]
\le \EE \left[ \left( \int_{-\infty}^{\infty} e^{-\gamma \Ba_t^\mu} Z_{\gamma} (dt)\right)^{p-1} \right].
\end{align*}

\noindent On the other hand, the uniform bound from Step 2 gives us a control on the residual term:
\begin{align*}
\EE \left[ \left( \int_{|t| > M} e^{-\gamma \Ba_t^\mu} Z_{\gamma, \epsilon \exp(L_{m})} (dt)\right)^{p-1} \right]^{\frac{1}{p-1}}
\le C_p(M)^{\frac{1}{p-1}}
\end{align*}

\noindent where $C_p(\cdot) \in (0, \infty)$ is defined in \eqref{eq:ub_integrand2} with the property that $C_p(M) \to 0$ as $M \to \infty$. Since $M > 0$ is arbitrary, we arrive at the upper bound
\begin{align}
\notag
& \limsup_{\epsilon \to 0} \EE \left[ \left( \int_{-L_{m} - \log r}^{(T-L_{m})_+} e^{-\gamma \Ba_t^\mu} Z_{\gamma, \epsilon}\circ \phi_{L_{m}} (dt)\right)^{p-1}\right]\\
\label{eq:Step3ub}
& \qquad\qquad  \le \EE \left[ \left( \int_{-\infty}^{\infty} e^{-\gamma \Ba_t^\mu} Z_{\gamma} (dt)\right)^{p-1} \right].
\end{align}

\noindent Combining \eqref{eq:Step3lb} and \eqref{eq:Step3ub}, we have
\begin{align*}
&  \mu(1-2\delta) \EE \left[ \left( \int_{-L_{m_0} - \log r}^{\infty} e^{-\gamma \Ba_t^\mu} Z_{\gamma} (dt)\right)^{p-1} \right]\\
& \qquad \le \liminf_{\epsilon \to 0^+} \frac{1}{T}  \int_0^{\mu(1-2\delta)T} \EE \left[ \left( \int_{-L_{m} - \log r}^{(T-L_{m})_+} e^{-\gamma \Ba_t^\mu} Z_{\gamma, \epsilon}\circ \phi_{L_{m}} (dt)\right)^{p-1}\right] \\
& \qquad \le\limsup_{\epsilon \to 0^+} \frac{1}{T}  \int_0^{\mu(1-2\delta)T} \EE \left[ \left( \int_{-L_{m} - \log r}^{(T-L_{m})_+} e^{-\gamma \Ba_t^\mu} Z_{\gamma, \epsilon}\circ \phi_{L_{m}} (dt)\right)^{p-1}\right] \\
&\qquad  \le \mu (1-2\delta)\EE \left[ \left( \int_{-\infty}^{\infty} e^{-\gamma \Ba_t^\mu} Z_{\gamma} (dt)\right)^{p-1} \right].
\end{align*}

\noindent The lower bound then matches the upper bound as we send $m_0$ to infinity, and this coincides with $\mu(1-2\delta) \overline{C}_{\gamma, d}$ by \eqref{eq:reflection2} in \Cref{lem:reflection}.

\paragraph{Step 4: remaining contribution.}
We can use the crude bound for the integrand from Step 2 so that
\begin{align*}
\frac{1}{T}\int_{\mu(1-2\delta) T}^{\mu(1+4\delta)T}  \EE \left[ \left( \int_{-\infty}^\infty e^{-\gamma \Ba_t^\mu} Z_{\gamma, \epsilon}\circ \phi_{L_{m}} (dt)\right)^{p-1}\right]
\le 6\delta \mu C_{p}
\end{align*}

\noindent where $C_p$ is defined in \eqref{eq:ub_integrand}. It suffices to control the remaining contribution from $m \ge \mu(1+4\delta)T$. We consider a different event here, namely 
\begin{align*}
\widetilde{A}_\delta 
= \{\tau_{(1-\delta)m} > (1+\delta)T\}
= \left\{ \max_{t \le (1+\delta)T} B_t^\mu < (1-\delta)m\right\}.
\end{align*}

\noindent On this event, it is actually more convenient to work with the original representation of the integrand, i.e. \eqref{eq:expfunc0}. Since $\tau_m > \tau_{(1-\delta)m} > T$ and $\max_{t \le T} B_t^\mu \le (1-\delta)m$, we see that 
\begin{align*}
&\EE \left[\left( \int_{-\log r}^{T\wedge \tau_m} e^{\gamma (B_t^\mu -m)} Z_{\gamma, \epsilon} (dt) + \int_0^{(T - \tau_m)_+} e^{-\gamma \Ba_t^\mu} Z_{\gamma, \epsilon} \circ \phi_{\tau_m}(dt)\right)^{p-1}  1_{\widetilde{A}_\delta}\right]\\
& \le e^{-(p-1)\delta \gamma m} \EE \left[\left( \int_{-\log r}^{\tau_{(1-\delta)m}} e^{\gamma (B_t^\mu -(1-\delta) m)} Z_{\gamma, \epsilon} (dt) \right)^{p-1}  1_{\widetilde{A}_\delta}\right]\\
& \le e^{-(p-1)\delta \gamma m} \EE \left[\left( \int_{-\infty}^0 e^{-\gamma \Ba_t^\mu} Z_{\gamma, \epsilon} \circ \phi_{L_{(1-\delta)m}} (dt) \right)^{p-1}  1_{\widetilde{A}_\delta}\right]\\
&\le C_p e^{-(p-1)\delta \gamma m}
\end{align*}

\noindent where the second inequality follows from time reversal (\Cref{lem:timerev}), and the constant $C_p < \infty$ in the last inequality is again the bound in \eqref{eq:ub_integrand}.

On the complementary event, we have 
\begin{align*}
& \EE \left[\left( \int_{-\log r}^{T\wedge \tau_m} e^{\gamma (B_t^\mu -m)} Z_{\gamma, \epsilon} (dt) + \int_0^{(T - \tau_m)_+} e^{-\gamma \Ba_t^\mu} Z_{\gamma, \epsilon} \circ \phi_{\tau_m}(dt)\right)^{p-1}  1_{\widetilde{A}_\delta^c}\right] \\
&\le 2^{p-1} \Bigg\{\EE \left[\left( \int_{-L_m}^{0} e^{-\gamma \Ba_t^\mu} Z_{\gamma, \epsilon} \circ \phi_{L_m}(dt)\right)^{p-1}  1_{\widetilde{A}_\delta^c}\right]\\
& \qquad \qquad + \EE \left[\left( \int_{0}^{(T - L_m)_+} e^{-\gamma \Ba_t^\mu} Z_{\gamma, \epsilon} \circ \phi_{L_m}(dt)\right)^{p-1}  1_{\widetilde{A}_\delta^c}\right]
\Bigg\}\\
& =: 2^{p-1}(I + II)
\end{align*}

\noindent where we abused the notation in the first inequality and considered the dual interpretation $\widetilde{A}_\delta^c = \left\{L_{(1-\delta)m} \le (1+\delta)T\right\}$ when we applied time reversal. Using the fact that $\PP(\Na(0,1) > t) \le e^{-t^2/2}$ and $m \ge \mu(1+4\delta) T$, we have
\begin{align*}
\PP\left(\widetilde{A}_\delta^c\right)
& = \PP\left(\tau_{(1-\delta)m} \le (1+\delta)T\right)\\
& \le \PP\left(\max_{t \le (1+\delta)T} B_t \ge (1-\delta)m - \mu T\right)\\
& = 2\PP\left(B_{(1+\delta)T}\ge (1-\delta)m - \mu T\right)\\
& \le2 \exp\left(-\frac{1}{2} \left(\frac{(1-\delta)m - \mu T}{\sqrt{(1+\delta)T}}\right)^2\right)
\le 2 \exp\left(-\frac{\mu}{2} \frac{3\delta - 4\delta^2}{1+\delta}  \left[(1-\delta)m - \mu T \right]\right)
\end{align*}

\noindent so that
\begin{align*}
II 
& \le \EE \left[\left( \int_{0}^{(T - L_m)_+} e^{-\gamma \Ba_t^\mu} Z_{\gamma, \epsilon} \circ \phi_{L_m}(dt)\right)^{p-1} \right] \PP\left(\widetilde{A}_\delta^c\right)\\
& \le 2 C_p \exp\left(-\frac{\mu}{2} \frac{3\delta - 4\delta^2}{1+\delta}  \left[(1-\delta)m - \mu T \right]\right)
\end{align*}

\noindent by independence (as $\widetilde{A}_\delta^c$ only depends on $(\Ba_t^\mu)_{t \le 0}$). As for the first term, pick $\alpha > 1$ small enough so that $1 \le \alpha(p-1) < p$, and we obtain
\begin{align*}
I  & \le \EE \left[\left( \int_{-L_m}^{0} e^{-\gamma \Ba_t^\mu} Z_{\gamma, \epsilon} \circ \phi_{L_m}(dt)\right)^{\alpha(p-1)}\right]^{1/\alpha}  \PP\left({\widetilde{A}_\delta^c}\right)^{1 - \frac{1}{\alpha}}\\
& \le C\EE \left[\left( \int_{-L_m}^{0} e^{-\gamma \Ba_t^\mu} Z_{\gamma}(dt)\right)^{\alpha(p-1)}\right]^{1/\alpha}  \PP\left({\widetilde{A}_\delta^c}\right)^{1 - \frac{1}{\alpha}}\\
& \le C'\EE \left[\left( \int_{-\infty}^0 e^{-\gamma \Ba_t^\mu} Z_{\gamma}(dt)\right)^{\alpha(p-1)}\right]^{1/\alpha}  \exp\left(-\frac{\mu}{2} \frac{3\delta - 4\delta^2}{1+\delta}  (1-\alpha^{-1})\left[(1-\delta)m - \mu T \right]\right)
\end{align*}

\noindent by H\"older's inequality and Gaussian comparison (to replace $Z_{\gamma, \epsilon} \circ \phi_{L_m}$ with $Z_{\gamma}$ up to a multiplicative error $C>0$). Our task now is to show that
\begin{align}\label{eq:Step4bound}
\EE \left[\left( \int_{-\infty}^0 e^{-\gamma \Ba_t^\mu} Z_{\gamma}(dt)\right)^{\alpha(p-1)}\right]^{\frac{1}{\alpha(p-1)}} < \infty.
\end{align}

For this, consider
\begin{align}
\notag
&\EE \left[\left( \int_{-\infty}^0 e^{-\gamma \Ba_t^\mu} Z_{\gamma}(dt)\right)^{\alpha(p-1)}\right]^{\frac{1}{\alpha(p-1)}} \\
\notag
&\qquad \le \sum_{n \ge 1} e^{-\gamma n}\EE \left[\left( \int_{L_{n-1}}^{L_{n}} e^{-\gamma (\Ba_{-t}^\mu - \Ba_{-L_{n-1}}^\mu)} Z_{\gamma}(dt)\right)^{\alpha(p-1)}\right]^{\frac{1}{\alpha(p-1)}}\\
\label{eq:remaining}
&\qquad \le \left(\sum_{n \ge 1} e^{-\gamma n}\right)\EE \left[\left( \int_{0}^{L_1} e^{-\gamma \Ba_{-t}^\mu} Z_{\gamma}(dt)\right)^{\alpha(p-1)}\right]^{\frac{1}{\alpha(p-1)}}
\end{align}

\noindent where the last equality follows from \Cref{cor:indep}. To analyse the expectation, we first split it into
\begin{align}\label{eq:remainsum}
\EE \left[\left( \int_{0}^{L_1} e^{-\gamma \Ba_{-t}^\mu} Z_{\gamma}(dt)\right)^{\alpha(p-1)}\right]^{\frac{1}{\alpha(p-1)}}
& \le \sum_{n \ge 0} \EE \left[\left( \int_{n}^{n+1} Z_{\gamma}(dt) \right)^{\alpha(p-1)} 1_{\{L_1 \ge n\}}\right]^{\frac{1}{\alpha(p-1)}} 
\end{align}

\noindent where
\begin{align*}
\EE \left[\left( \int_{n}^{n+1} Z_{\gamma}(dt) \right)^{\alpha(p-1)}\right]^{\frac{1}{\alpha(p-1)}}
= \EE \left[\left( \int_{0}^{1} Z_{\gamma}(dt) \right)^{\alpha(p-1)}\right]^{\frac{1}{\alpha(p-1)}} 
\end{align*}

\noindent by translation invariance (\Cref{lem:lateralinvariant}), and this quantity is finite due to the existence of GMC moments. Meanwhile, from \Cref{theo:2m-b}, we have
\begin{align*}
\PP\left(L_1 \ge n\right)
= \PP\left(\inf_{s \ge n} \Ba_s^\mu \le 1 \right)
&= \PP\left(\max_{s \le n} B_s^\mu \le 1\right)\\
& \le \PP\left(B_n + \mu n \le 1\right)
\le \exp\left(-\frac{1}{2n} (1-\mu n)^2\right).
\end{align*}

\noindent Applying H\"older's inequality, \eqref{eq:remainsum} can be further bounded by
\begin{align*}
&\sum_{n \ge 0} \EE \left[\left( \int_{n}^{n+1} Z_{\gamma}(dt) \right)^{\alpha q(p-1)}\right]^{\frac{1}{\alpha q(p-1)}}  \PP\left(L_1 \ge n\right)^{1 - \frac{1}{q}}
\end{align*}

\noindent which is summable provided that $q > 1$ is chosen such that $\alpha q(p-1) < p$ so that
\begin{align*}
\EE \left[\left( \int_{n}^{n+1} Z_{\gamma}(dt) \right)^{\alpha q(p-1)}\right]
=\EE \left[\left( \int_{0}^{1} Z_{\gamma}(dt) \right)^{\alpha q(p-1)}\right] < \infty
\end{align*}

\noindent by the existence of GMC moments again. Therefore, \eqref{eq:remaining} and hence \eqref{eq:Step4bound} are all finite, as claimed.

Gathering all terms, we see that there exists some constant $C>0$ independent of $T$ (or equivalently $\epsilon$) and $m \ge \mu (1+4\delta)T$ (but possibly depends on other parameters) such that
\begin{align*}
& \EE \left[\left( \int_{-\log r}^{T\wedge \tau_m} e^{\gamma (B_t^\mu -m)} Z_{\gamma, \epsilon} (dt) + \int_0^{(T - \tau_m)_+} e^{-\gamma \Ba_t^\mu} Z_{\gamma, \epsilon} \circ \phi_{\tau_m}(dt)\right)^{p-1}  \right]\\
& \qquad \le C \exp\left(-C^{-1} \left((1-\delta)m - \mu T\right)\right)
\end{align*}

\noindent and hence
\begin{align*}
& \int_{m \ge \mu(1+4\delta)T} dm \EE \left[\left( \int_{-\log r}^{T\wedge \tau_m} e^{\gamma (B_t^\mu -m)} Z_{\gamma, \epsilon} (dt) + \int_0^{(T - \tau_m)_+} e^{-\gamma \Ba_t^\mu} Z_{\gamma, \epsilon} \circ \phi_{\tau_m}(dt)\right)^{p-1}  \right]\\
& \qquad \le \frac{1}{1-\delta} e^{-\mu \left[(1-\delta)(1+4\delta) - 1\right] T} \le 2
\end{align*}

\noindent for $\delta > 0$ sufficiently small.

\paragraph{Step 5: conclusion.} Summarising all the work we have done so far, we have
\begin{align*}
& \EE\left[\left(\int_{K\epsilon \le |x| \le r} e^{\gamma \overline{X}_{\epsilon}(x)- \frac{\gamma^2}{2} \EE[\overline{X}_{\epsilon}(x)^2]}\frac{dx}{|x|^{\gamma^2}} \right)^{p-1}\right]\\
& = \left(1 + \mathcal{O}(K^{-1})\right) 2\mu \int_0^{\infty} dm \EE \Bigg[\Bigg( \int_{-\log r}^{T\wedge \tau_m} e^{\gamma (B_t^\mu -m)} Z_{\gamma, \epsilon} (dt) \\
& \qquad \qquad \qquad \qquad \qquad \qquad \qquad  
+ \int_0^{(T - \tau_m)_+} e^{-\gamma \Ba_t^\mu} Z_{\gamma, \epsilon} \circ \phi_{\tau_m}(dt)\Bigg)^{p-1}  \Bigg]\\
& =  \left(1 + \mathcal{O}(K^{-1})\right) 2\mu  \left[\mu \overline{C}_{\gamma, d} T + \mathcal{O}(1+\delta T)\right]
\end{align*}

\noindent where $T = -\log K\epsilon$. Dividing both sides by $-\log \epsilon$ and taking the limit $\epsilon \to 0^+$, the RHS becomes $\left(1 + \mathcal{O}(K^{-1})\right) 2\mu^2 \overline{C}_{\gamma, d} + \Oa(\delta)$. Since $K > 1$ and $\delta > 0$ are arbitrary, we send $K$ to infinity and $\delta$ to zero and obtain the desired result.

\begin{remark}\label{rem:gextend}
To treat the case where $g \ge 0$ is continuous on $\overline{D}$, one would first obtain the analogue of \Cref{lem:main1} by following the same argument, i.e.
\begin{align*}
& \EE\left[\left(\int_D g(x)M_{\gamma, \epsilon}(dx)\right)^p\right]\\
&= \int_D g(u) du\EE\left[\left(\int_{D \cap \{K\epsilon \le |x-u| \le r\}} g(x) \frac{e^{\gamma^2 f_\epsilon(x, u)}}{|x-u|^{\gamma^2}} M_{\gamma, \epsilon}(dx)\right)^{p-1}\right] + \mathcal{O}(1).
\end{align*}

When $r$ is small, $g(x) \approx g(u)$ for $|x-u| \le r$, and so the integral above can be further approximated by
\begin{align*}
&= \int_D g(u)^p du\EE\left[\left(\int_{D \cap \{K\epsilon \le |x-u| \le r\}} \frac{e^{\gamma^2 f_\epsilon(x, u)}}{|x-u|^{\gamma^2}} M_{\gamma, \epsilon}(dx)\right)^{p-1}\right]
\end{align*}

\noindent and no further changes are needed for the rest of the analysis. This explains the appearance of $g(u)^p$ in the leading order coefficient in \Cref{theo:main}.
\end{remark}

\subsection{General approximation: proof of \Cref{cor:main}}\label{subsec:genapp}
Not surprisingly, the generalisation of \Cref{theo:main} to \Cref{cor:main} relies on Gaussian comparison. To achieve this, we need the following technical estimate.
\begin{lemma}\label{lem:interpolate_help}
Let $\epsilon, \delta \in (0, 1)$ and $r = \epsilon^{1-\delta}$. Then for $\gamma \in (0, \sqrt{2d})$, $p = p_c = 2d/\gamma^2$, we have
\begin{align*}
\EE\left[M_{\gamma, \epsilon}(B(x_0, r)\cap D)^2 M_{\gamma, \epsilon}(D)^{p-2}\right]
\le  C r^d \log \frac{r}{\epsilon}
\end{align*}

\noindent where the constant $C>0$ is uniform in $\epsilon, \delta$, $x_0 \in D$, and can be further chosen to be uniformly for all Gaussian fields $X_\epsilon$ with covariance of the form
\begin{align*}
\EE\left[X_\epsilon(x) X_\epsilon(y)\right] = -\log \left(|x- y| \vee \epsilon \right) + f_\epsilon(x, y) \qquad \forall x, y \in D
\end{align*}

\noindent where $\sup_{\epsilon > 0} \sup_{x, y \in D} |f_\epsilon(x, y)| \le K$ for some $K > 0$.
\end{lemma}

\begin{proof}
To simplify the notations let us verify the inequality for $x_0 = 0$, assuming that the domain contains the origin and that $\mathrm{diam}(D) := \sup \{|x-y| : x, y \in D\} \le 1$. The constant $C>0$ below will vary from line to line but its uniformity in $\epsilon, \delta$ and $x_0 \in D$ should be evident from the proof.

To begin with, consider
\begin{align}
\notag
& \EE\left[M_{\gamma, \epsilon}(B(x_0, r))^2 M_{\gamma, \epsilon}(D)^{p-2}\right]\\
\notag
& = \int_{|u| \le r} \int_{|v| \le r} du dv\EE\left[ e^{\gamma (X_\epsilon(u) + X_\epsilon(v)) - \frac{\gamma^2}{2} \EE[X_\epsilon(u)^2 + X_\epsilon(v)^2]}M_{\gamma, \epsilon}(D)^{p-2}\right]\\ 
\label{eq:ip_estimate}
&
\le C \int_{|u|\le r} du \int_{|v| \le r} \frac{dv}{(|u-v| \vee \epsilon)^{\gamma^2}}
\EE\left[\left( \int_D \frac{M_{\gamma, \epsilon}(dx)}{(|x-u|\vee \epsilon)^{\gamma^2}(|x-v|\vee \epsilon)^{\gamma^2}} \right)^{p-2}\right].
\end{align}

\noindent Before we continue with our bounds, let us mention that the uniformity of $C>0$ over the class of fields $X_\epsilon$ stated in the lemma is a simple consequence of the fact that the expectation that appears in the integrand in \eqref{eq:ip_estimate} can be bounded uniformly in this class via Gaussian comparison (see \Cref{cor:GP}).

Now, let us split the $v$-integral in \eqref{eq:ip_estimate} into two parts, namely when $|u-v| \le 2\epsilon$ or $2\epsilon \le |u-v| \le r$. We proceed by considering two separate cases.
\begin{itemize}
\item Case 1: $p \ge 2$.  For the first part where $|u-v| \le 2\epsilon$, we have
\begin{align*}
&\int_{|u-v| \le 2\epsilon} \frac{dv}{(|u-v| \vee \epsilon )^{\gamma^2}}\EE\left[\left( \int_D \frac{M_{\gamma, \epsilon}(dx)}{(|x-u|\vee \epsilon)^{\gamma^2}(|x-v|\vee \epsilon)^{\gamma^2}} \right)^{p-2}\right]\\
& \qquad \le C \epsilon^{d-\gamma^2}\EE\left[\left( \int_D \frac{M_{\gamma, \epsilon}(dx)}{(|x-u|\vee \epsilon)^{2\gamma^2}} \right)^{p-2}\right]\\
& \qquad \le C \epsilon^{d-\gamma^2}
\Bigg\{\EE\left[\left( \epsilon^{-2\gamma^2} M_{\gamma, \epsilon}(B(u, \epsilon)) \right)^{p-2}\right] \\
& \qquad \qquad \qquad \qquad+
\EE\left[\left( \int_{\{|x-u| \ge \epsilon\} \cap D} \frac{M_{\gamma, \epsilon}(dx)}{(|x-u|\vee \epsilon)^{2\gamma^2}} \right)^{p-2}\right]\Bigg\}\\
& \qquad \le C \epsilon^{d-\gamma^2} \epsilon^{-\frac{\gamma^2}{2}(p-2)} \le C
\end{align*}

\noindent where the second last inequality follows from \Cref{lem:GMC_multi}.

\noindent As for the part where $|u-v| \in [2\epsilon, r]$, we first focus on the expected value in \eqref{eq:ip_estimate} and split it into two terms, depending on whether $x$ is closer to $u$ or $v$. By symmetry, it suffices to deal with
\begin{align*}
& \EE\left[\left( \int_{\{ |x-u| \le |x-v| \} \cap D} \frac{M_{\gamma, \epsilon}(dx)}{(|x-u|\vee \epsilon)^{\gamma^2}(|x-v|\vee \epsilon)^{\gamma^2}} \right)^{p-2}\right]\\
& \qquad \le C \Bigg\{
\EE\left[\left( \int_{\{ |x-u| \ge 2|u-v| \}} \frac{M_{\gamma, \epsilon}(dx)}{|x-u|^{2\gamma^2}} \right)^{p-2}\right]\\
& \qquad \qquad \qquad \qquad 
+ \EE\left[\left( \int_{\{ |x-u| \le 2|u-v| \}} \frac{M_{\gamma, \epsilon}(dx)}{(|x-u|\vee \epsilon)^{\gamma^2}|u-v|^{\gamma^2}} \right)^{p-2}\right]
\Bigg\}
\end{align*}

\noindent where 
\begin{align*}
\EE\left[\left( \int_{\{ |x-u| \ge 2|u-v| \}} \frac{M_{\gamma, \epsilon}(dx)}{|x-u|^{2\gamma^2}} \right)^{p-2}\right]
&\le C |u-v|^{-\frac{\gamma^2}{2}(p-2)}
\end{align*}

\noindent and
\begin{align*}
& \EE\left[\left( \int_{\{ |x-u| \le 2|u-v| \}} \frac{M_{\gamma, \epsilon}(dx)}{(|x-u|\vee \epsilon)^{\gamma^2} |u-v|^{\gamma^2}} \right)^{p-2}\right] \\
& \qquad \le 
C |u-v|^{-\gamma^2(p-2)} \EE\left[\left( \int_{\{|x-u| \le 2|u-v| \}} \frac{M_{\gamma, \epsilon}(dx)}{|x-u|^{\gamma^2}} \right)^{p-2}\right]\\
& \qquad \le C |u-v|^{-\gamma^2(p-2)} |u-v|^{\frac{\gamma^2}{2}(p-2)}
\le C |u-v|^{-\frac{\gamma^2}{2}(p-2)}
\end{align*}

\noindent again by \Cref{lem:GMC_multi}. Therefore,
\begin{align*}
& \int_{2\epsilon \le |u-v| \le r} \frac{dv}{(|u-v| \vee \epsilon )^{\gamma^2}}\EE\left[\left( \int_D \frac{M_{\gamma, \epsilon}(dx)}{(|x-u|\vee \epsilon)^{\gamma^2}(|x-v|\vee \epsilon)^{\gamma^2}} \right)^{p-2}\right]\\
& \qquad \le C\int_{2\epsilon \le |u-v| \le r} \frac{dv}{|u-v|^{\gamma^2 + \frac{\gamma^2}{2}(p-2)}}\\
& \qquad =  C \int_{2\epsilon \le |u-v| \le r} \frac{dv}{|u-v|^d} 
\le C \log \frac{r}{\epsilon}
\end{align*}

\noindent which implies the desired inequality.

\item Case 2: $p < 2$. For the first part where $|u-v| \le 2\epsilon$, we compare the Gaussian field $X_\epsilon$ with a Gaussian random variable $\Na_\epsilon$ with variance $-\log \epsilon$ on $|x-u| \le \epsilon$, and obtain
\begin{align}
\notag
&\int_{|u-v| \le 2\epsilon} \frac{dv}{(|u-v| \vee \epsilon )^{\gamma^2}}\EE\left[\left( \int_D \frac{M_{\gamma, \epsilon}(dx)}{(|x-u|\vee \epsilon)^{\gamma^2}(|x-v|\vee \epsilon)^{\gamma^2}} \right)^{p-2}\right]\\
\notag
& \qquad \le C \int_{|u-v| \le 2\epsilon} \frac{dv}{(|u-v| \vee \epsilon )^{\gamma^2}}\EE\left[\left( \int_{\{|x-u| \le \epsilon\} \cap D} \frac{e^{\gamma \Na_\epsilon - \frac{\gamma^2}{2} \EE[\Na_\epsilon^2]}dx}{\epsilon^{\gamma^2}(3\epsilon)^{\gamma^2}} \right)^{p-2}\right]\\
\notag
& \qquad \le C \epsilon^{d-\gamma^2} \epsilon^{-\frac{\gamma^2}{2}(p-2)} \left(\int_{\{|x-u| \le \epsilon \} \cap D} \epsilon^{-d} dx \right)^{p-2}\\
\label{eq:coruniform1}
& \qquad = C \left(\int_{\{|x-u| \le \epsilon \} \cap D} \epsilon^{-d} dx \right)^{p-2}
\end{align}

\noindent by \Cref{cor:GP}. This can be uniformly bounded in $u \in \overline{D}$ by 
\begin{align*}
C \epsilon^{-d} \int_{\overline{B}(u, \epsilon)} dx \le C \int_{\overline{B}(0,1)} dx \le C.
\end{align*}

The second part  may be treated in a similar way. If $|x-u| \le \frac{1}{2}|u-v|$, it is easy to check by triangle inequality that $\frac{1}{2} |u-v| \le |x-v| \le \frac{3}{2} |u-v|. $ Let us write $c = c(u, v):= \frac{1}{2}|u-v|$ If $\Na_{c} \sim \Na(0, -\log c)$ is an independent variable, we have
\begin{align*}
\left|\EE\left[X_\epsilon(x-u) X_\epsilon(y-u)\right] -
\EE\left[X_{\epsilon/c}(c^{-1}(x-u)) + \Na_c)( X_{\epsilon/c}(c^{-1}(y-u)) + \Na_c)\right]\right|
\end{align*}

\noindent for any $x, y \in \overline{B}(u, c)$ and hence
\begin{align*}
& \int_{B(0,r)}du \int_{2\epsilon \le |u-v| \le r} \frac{dv}{|u-v|^{\gamma^2}} \EE\left[\left( \int_{D} \frac{M_{\gamma, \epsilon}(dx)}{(|x-u|\vee \epsilon)^{\gamma^2}(|x-v|\vee \epsilon)^{\gamma^2}} \right)^{p-2}\right]\\
& \le \int_{B(0,r)}du \int_{2\epsilon \le |u-v| \le r} \frac{dv}{|u-v|^{\gamma^2}}\EE\left[\left( \int_{\{|x-u| \le \frac{1}{2}|u-v|\} \cap D} \frac{M_{\gamma, \epsilon}(dx)}{(|x-u|\vee \epsilon)^{\gamma^2}(|x-v|\vee \epsilon)^{\gamma^2}} \right)^{p-2}\right]\\
& \le C\int_{B(0,r)}du \int_{2\epsilon \le |u-v| \le r} \frac{dv}{|u-v|^{\gamma^2}} \EE\left[\left( \int_{\overline{B}(0, c) \cap (D-u)} \frac{e^{\gamma \Na_{c} - \frac{\gamma^2}{2} \EE[\Na_{c}^2]} e^{\gamma X_{\epsilon/c}(x/c) - \frac{\gamma^2}{2} \EE[X_{\epsilon/c}(x/c)^2]}dx}{|u-v|^{2\gamma^2}} \right)^{p-2}\right]\\
& = C\int_{B(0,r)}du \int_{2\epsilon \le |u-v| \le r} \frac{dv}{|u-v|^{d}} \EE\left[\left(\int_{\overline{B}(0, 1) \cap \left(c^{-1}(D-u)\right)} M_{\gamma, \epsilon/c}(dx)\right)^{p-2} \right].
\end{align*}

\noindent Arguing as before (and combined with the convergence and existence of GMC moments, see \Cref{lem:GMCbasic}) we can obtain a bound for the last expectation that is uniform for all $u, v \in D$ and $\epsilon > 0$. Hence the above expression is bounded by $C r^d \log \frac{r}{\epsilon}$ and the proof is complete.
\end{itemize}
\end{proof}
\begin{proof}[Proof of \Cref{cor:main}]
Again for simplicity we just treat the case $g \equiv 1$ here. Let us construct on the same probability space a Gaussian field $X$ with covariance \eqref{eq:LGF}, independent of $\{X_\epsilon\}_\epsilon$, and write $\widetilde{X}_\epsilon(x) = X \ast \nu_\epsilon$ for some mollifier $\nu \in C_c^\infty(\RR^d)$. We further introduce the notations
\begin{align*}
X_\epsilon^t(x) 
= \sqrt{t} \widetilde{X}_\epsilon(x) + \sqrt{1-t} X_\epsilon(x), 
\qquad M_{\gamma, \epsilon}^t(dx) 
= e^{\gamma X_\epsilon^t(x) - \frac{\gamma^2}{2}\EE[X_\epsilon^t(x)^2]}dx
\end{align*}

\noindent for all $t \in [0,1]$. Using the interpolation formula from \Cref{lem:Kahane}, we have
\begin{align}
\notag
&\left|\EE\left[M_{\gamma, \epsilon}^1(D)^p\right]
- \EE\left[M_{\gamma, \epsilon}^0(D)^p\right]\right|\\
\notag
& \qquad \le \frac{p(p-1) }{2} \int_0^1 dt\int_D \int_D \left| \EE[X_\epsilon(x) X_\epsilon(y)]-\EE[ \widetilde{X}_\epsilon(x) \widetilde{X}_\epsilon(y)]\right|\\
\label{eq:interpolate_error}
& \qquad \qquad \qquad \qquad \qquad \qquad \times\EE\left[M_{\gamma, \epsilon}^t(dx) M_{\gamma, \epsilon}^t(dy)M_{\gamma, \epsilon}^t(D)^{p-2}\right].
\end{align}

Now fix $\delta > 0$ and write $r = \epsilon^{1-\delta}$. By assumption, we have for $\epsilon > 0$ sufficiently small
\begin{align*}
\sup_{|x-y| \ge r} \left|\EE[X_\epsilon(x) X_\epsilon(y)] - \EE[\widetilde{X}_\epsilon(x) \widetilde{X}_\epsilon(y)]\right| \le \delta.
\end{align*}

\noindent Together with the other assumption on the covariance of $X_\epsilon$, we can bound \eqref{eq:interpolate_error} by
\begin{align}
\notag
&\frac{p(p-1)\delta }{2}  \sup_{t \in [0, 1]} \EE\left[M_{\gamma, \epsilon}^t(D)^{p}\right]\\
\label{eq:interpolate_error2}
&+ \frac{Cp(p-1)}{2} \int_0^1 dt\int_D \int_{\{|x-y| \le r\}\cap D} \EE\left[M_{\gamma, \epsilon}^t(dx) M_{\gamma, \epsilon}^t(dy) M_{\gamma, \epsilon}^t(D)^{p-2}\right].
\end{align}

Using Gaussian comparison, $\EE\left[M_{\gamma, \epsilon}^t(D)^p\right]$ is, up to a multiplicative factor, upper bounded by $\EE\left[M_{\gamma, \epsilon}^{1}(D)^p\right]$. It then follows from \Cref{theo:main} that the first term in \eqref{eq:interpolate_error2} is of order $\Oa(\delta \log \frac{1}{\epsilon})$. As for the second term, we can bound it by
\begin{align*}
& \frac{Cp(p-1)}{2}\int_0^1 dt \int_D \EE\left[M_{\gamma, \epsilon}^t(dx) M_{\gamma, \epsilon}^t(B(x, r) \cap D) M_{\gamma, \epsilon}^t(D)^{p-2}\right]\\
& \qquad \le \frac{Cp(p-1)}{2} \sum_{x \in D_r}\sup_{t \in [0, 1]} \EE\left[M_{\gamma, \epsilon}^t(B(x, 2r) \cap D)^2 M_{\gamma, \epsilon}^t(D)^{p-2}\right]
\end{align*}

\noindent where $D_r$ is a collection of points $x\in D$ such that $D \subset \bigcup_{x \in D} B(x, r)$. We can of course choose $D_r$ such that its cardinality is at most $\mathcal{O}(r^{-d})$. Combining this with \Cref{lem:interpolate_help} (which holds uniformly in $t \in [0, 1]$) we see that the second term is also $\Oa(\delta \log \frac{1}{\epsilon})$, and therefore
\begin{align*}
\limsup_{\epsilon \to 0^+} \frac{1}{-\log \epsilon} \left|\EE\left[M_{\gamma, \epsilon}^1(D)^p\right]
- \EE\left[M_{\gamma, \epsilon}^0(D)^p\right]\right|
= \Oa(\delta).
\end{align*}

\noindent Since $\delta > 0$ is arbitrary, we let $\delta \to 0^+$ and conclude that $\EE\left[M_{\gamma, \epsilon}^1(D)^{p}\right]$ has the same leading order asymptotics as $\EE\left[M_{\gamma, \epsilon}^0(D)^{p}\right]$, which is described by \Cref{theo:main}.
\end{proof}

\subsection{On circular unitary ensemble: proof of \Cref{theo:CUEinteger}} \label{subsec:CUEproof}
The $k=2$ case was first proved in  \cite[Theorem 1.15]{CK2015}, relying on the full Painlev\'e V description of the strong asymptotics of Toeplitz determinant with 2 merging singularities, but it is not clear how this may be extended to deal with $k > 2$ singularities with arbitrary merging behaviour. We shall show, however, that the asymptotics for critical-subcritical moments can be derived using only the following ingredients.
\begin{lemma}\label{lem:FH}
The following estimates holds for any integers $k \ge 1$.
\begin{itemize}
\item[(i)] (cf. \cite[Theorem 1.1]{DIK2014}) For any $s \ge 0$ and any $\delta > 0$,
\begin{align*}
\EE_{\mathrm{U(N)}}\left[\prod_{j=1}^k |P_N(\theta_j)|^{2s}\right] 
= (1+o(1)) N^{ks^2}\left[ \frac{G(1+s)^2}{G(1+2s)}\right]^k  \prod_{1 \le u < v \le k}|e^{i\theta_u} - e^{i\theta_v}|^{-2s^2}
\end{align*}

uniformly in $\theta_1, \dots, \theta_k \in [0, 2\pi)$ satisfying $|e^{i\theta_u} -e^{i \theta_v}| \ge N^{-(1-\delta)}$ for any distinct $u, v \in 1, \dots, k$.
\item[(ii)] (cf. \cite[Theorem 1.1]{Fa2019}) For any $s \ge 0$, there exists some constant $C>0$ such that
\begin{align*}
\EE_{\mathrm{U(N)}}\left[\prod_{j=1}^k |P_N(\theta_j)|^{2s}\right] 
\le C N^{ks^2} \prod_{1 \le u < v \le k}\left(|e^{i\theta_u} - e^{i\theta_v}|\vee \frac{1}{N}\right)^{-2s^2}
\end{align*}

\noindent uniformly in $\theta_1, \dots, \theta_k \in [0, 2\pi)$ and for all $N$ sufficiently large.
\end{itemize}
\end{lemma}

\begin{remark}
The asymptotics in \Cref{lem:FH} (i) was originally stated for fixed singularities $\theta_1, \dots, \theta_k \in [0, 2\pi)$ in \cite{DIK2014}, but can be shown to hold in the current setting. We explain this briefly in \Cref{app:FH}.
\end{remark}

\begin{proof}[Proof of \Cref{theo:CUEinteger}]
Let $X(e^{i\theta})$ be the log-correlated Gaussian field in \eqref{eq:hcue1} and $\{X_\epsilon(e^{i\theta})\}_{\epsilon > 0}$ a sequence of approximate fields in the sense of \Cref{cor:main}. If we set $\epsilon = \epsilon(N) = N^{-1}$, $s = 1/\sqrt{k}$ and $\gamma = \sqrt{2}s$, then \Cref{lem:FH} may be re-interpreted as follows:
\begin{itemize}
\item For any $\delta \in (0, 1)$,
\begin{align*}
\frac{\EE_{\mathrm{U}(N)} \left[\prod_{j=1}^k|P_N(\theta_j)|^{2s}\right]}{\prod_{j=1}^k\EE_{\mathrm{U}(N)} \left[|P_N(\theta_j)|^{2s}\right]}
\overset{N \to \infty}{=} (1+o(1)) \frac{\EE\left[\prod_{j=1}^k e^{\gamma X_{\epsilon}(e^{i\theta_j})} \right]}{\prod_{j=1}^k\EE\left[ e^{\gamma X_{\epsilon}(e^{i\theta_j})} \right]}
\end{align*}

\noindent for some $o(1)$ depending on $\delta$ but uniformly in $\theta_1, \dots, \theta_k \in [0, 2\pi)$ satisfying $|e^{i\theta_u} -e^{i \theta_v}| \ge N^{-(1-\delta)}$ for any distinct $u, v \in \{1, \dots, k\}$.

\item There exists some constant $C>0$ such that
\begin{align*}
\frac{\EE_{\mathrm{U}(N)} \left[\prod_{j=1}^k|P_N(\theta_j)|^{2s}\right]}{\prod_{j=1}^k\EE_{\mathrm{U}(N)} \left[|P_N(\theta_j)|^{2s}\right]}
\le C \frac{\EE\left[\prod_{j=1}^k e^{\gamma X_{\epsilon}(e^{i\theta_j})} \right]}{\prod_{j=1}^k\EE\left[ e^{\gamma X_{\epsilon}(e^{i\theta_j})} \right]}
\end{align*}

\noindent uniformly in $\theta_1, \dots, \theta_k \in [0,2\pi)$.
\end{itemize}

For any $\delta_1, \dots, \delta_{k-1} > 0$, let us define 
\begin{align*}
D_{(\delta_1, \dots, \delta_{k-1})}
& = D_{(\delta_1, \dots, \delta_{k-1})}(N)\\
& = \{ (\theta_1, \dots, \theta_k) \in [0, 2\pi)^k: |e^{i\theta_u} - e^{i\theta_v}| \ge \epsilon^{1-\delta_u} \quad \text{for any } v > u, u \in [k]\}
\end{align*}

\noindent and $D_{(\delta_1, \dots, \delta_{k-1})}^c = [0, 2\pi)^k \setminus D_{(\delta_1, \dots, \delta_{k-1})}$. Then we may expand the integer moments and write
\begin{align*}
& \mathrm{MoM}_{\mathrm{U}(N)}(k, s)
= \frac{1}{(2\pi)^k} \int_{[0, 2\pi)^k}\EE_{\mathrm{U}(N)} \left[\prod_{j=1}^k|P_N(\theta_j)|^{2s}\right] d\theta_1 \dots d\theta_k\\
& \qquad = \frac{\EE_{\mathrm{U(N)}}[|P_N(0)|^{2s}]^k}{(2\pi)^k}  \int_{[0, 2\pi)^k}\frac{\EE_{\mathrm{U}(N)} \left[\prod_{j=1}^k|P_N(\theta_j)|^{2s}\right]}{\prod_{j=1}^k\EE_{\mathrm{U}(N)} \left[|P_N(\theta_j)|^{2s}\right]} d\theta_1 \dots d\theta_k\\
& \qquad = \frac{\EE_{\mathrm{U(N)}}[|P_N(0)|^{2s}]^k}{(2\pi)^k}
\Bigg\{  
\int_{D_{(\delta_1, \dots, \delta_{k-1})}}\EE\left[\prod_{j=1}^k M_{\gamma, \epsilon}(d\theta_j) \right]
+ \Oa\left(\int_{D_{(\delta_1, \dots, \delta_{k-1})}^c}\EE\left[\prod_{j=1}^k M_{\gamma, \epsilon}(d\theta_j) \right]\right)\Bigg\}\\
& \qquad = \frac{\EE_{\mathrm{U(N)}}[|P_N(0)|^{2s}]^k}{(2\pi)^k}
\Bigg\{  
\EE\left[M_{\gamma, \epsilon}([0, 2\pi))^k \right]
+ \Oa\left(\int_{D_{(\delta_1, \dots, \delta_{k-1})}^c}\EE\left[\prod_{j=1}^k M_{\gamma, \epsilon}(d\theta_j) \right]\right)\Bigg\}
\end{align*}

\noindent where $M_{\gamma, \epsilon}(d\theta) = e^{\gamma X_\epsilon(e^{i\theta}) - \frac{\gamma^2}{2} \EE[X_\epsilon(e^{i\theta})^2]}d\theta$. All we have to show now is that for suitably chosen $\delta_1, \dots, \delta_{k-1}$, the estimate
\begin{align}\label{eq:FH_goal}
\EE\left[\int_{D_{(\delta_1, \dots, \delta_{k-1})}^c}\prod_{j=1}^k M_{\gamma, \epsilon}(d\theta_j)\right]
= \Oa\left(c(\delta_1, \dots, \delta_{k-1}; \epsilon) \log \frac{1}{\epsilon}\right).
\end{align}

\noindent holds for some $c(\delta_1, \dots, \delta_{k-1}; \epsilon) > 0$ such that $c(\delta_1, \dots, \delta_{k-1}; \epsilon)$ tends to zero when we first let $\epsilon \to 0$ and then $\delta_1, \dots, \delta_{k-1} \to 0$.

We proceed in the following way. Let us fix $\delta_1 > 0$, and consider\footnote{We suppress the fact that all our dummy variables lie in $[0, 2\pi)$ to simplify our notations in the proof.}
\begin{align}\label{eq:CM_FH1}
&\EE\left[\left(\int_{|e^{iu} - e^{i\theta_1}| \le \epsilon^{1-\delta_1}} \frac{M_{\gamma, \epsilon}(du)}{(|e^{iu} - e^{i\theta_1}|\vee \epsilon)^{\gamma^2}}\right)^{k-1} \right].
\end{align}

\noindent Using the fact that $|e^{iu} - e^{i\theta_1}| = \Omega(|u-\theta_1|)$ where  $|u-\theta_1|$ should be interpreted as the distance on $[0, 2\pi)$ with the end-points identified, we have the following multifractal estimates for (regularised) GMC measures (see \Cref{lem:GMC_multi}):
\begin{align*}
\EE\left[\left(\int_{|e^{iu} - e^{i\theta_1}| \le \epsilon} \frac{M_{\gamma, \epsilon}(du)}{(|e^{iu} - e^{i\theta_1}|\vee \epsilon)^{\gamma^2}}\right)^{k-1} \right] & \le C,\\
\EE\left[\left(\int_{\epsilon \le |e^{iu} - e^{i\theta_1}| \le \epsilon^{1-\delta_1}} \frac{M_{\gamma, \epsilon}(du)}{(|e^{iu} - e^{i\theta_1}|\vee \epsilon)^{\gamma^2}}\right)^{k-1} \right] & \le C \delta_1 \log \frac{1}{\epsilon}.
\end{align*}

\noindent (In this proof we use $C > 0$ to denote some constant which varies from line to line below but is uniform in every relevant parameter.) In particular, \eqref{eq:CM_FH1} is of order $\Oa(\delta_1 \log \frac{1}{\epsilon})$, and
\begin{align*}
& \EE\left[M_{\gamma, \epsilon}([0, 2\pi))^k \right]
= \int_{[0, 2\pi)} d\theta_1 \EE\left[e^{\gamma X_\epsilon(e^{i\theta_1}) - \frac{\gamma^2}{2}\EE[X_\epsilon(e^{i\theta})^2]} \left(\int_{[0, 2\pi)} M_{\gamma, \epsilon}(du)\right)^{k-1} \right] \\
& \qquad = \int_{[0, 2\pi)} d\theta_1 \EE\left[\left(\int_{[0, 2\pi)} e^{\gamma^2 \EE[X_\epsilon(e^{i\theta_1}) X_\epsilon(e^{iu})]} M_{\gamma, \epsilon}(du)\right)^{k-1} \right] \\
& \qquad = \Oa\left(\delta_1^{\frac{1}{k-1}} \log \frac{1}{\epsilon}\right) +  (1+o(1)) \int_{[0, 2\pi)} d\theta_1 \EE\left[\left(\int_{|e^{iu}-e^{i\theta_1}| \ge \epsilon^{1-\delta_1}}  \frac{M_{\gamma, \epsilon}(du)}{|e^{i u} - e^{i\theta_1}|^{\gamma^2}}\right)^{k-1} \right]
\end{align*}

\noindent where the second equality follows from Cameron-Martin (\Cref{lem:CMG}), and in the last line the error term $\Oa\left(\delta_1^{\frac{1}{k-1}} \log \frac{1}{\epsilon}\right)$ comes from the removal of contribution from $|e^{iu} - e^{i\theta_1}| \le \epsilon^{1-\delta}$, combined with the fact that the leading order is $\Omega(\log \frac{1}{\epsilon})$ (\Cref{cor:main}) and the basic estimate in \Cref{lem:basic}. If $k=2$ then the previous identity reads as
\begin{align*}
& \EE\left[M_{\gamma, \epsilon}([0, 2\pi))^2 \right]
 = \Oa\left(\delta_1 \log \frac{1}{\epsilon}\right) +  (1+o(1)) \int_{[0, 2\pi)} d\theta_1 \int_{|e^{iu}-e^{i\theta_1}| \ge \epsilon^{1-\delta_1}} \frac{du}{|e^{i u} - e^{i\theta_1}|^{\gamma^2}}
\end{align*}

\noindent and we have verified \eqref{eq:FH_goal} successfully.

If $k \ge 3$ we continue to `bring down the powers':
\begin{align*}
& \int d\theta_1 \EE\left[\left(\int_{|e^{iu}-e^{i\theta_1}| \ge \epsilon^{1-\delta_1}}  \frac{M_{\gamma, \epsilon}(du)}{|e^{i u} - e^{i\theta_1}|^{\gamma^2}}\right)^{k-1} \right]\\
& = \int d\theta_1 \int_{|e^{i\theta_2} - e^{i\theta_1}| \ge \epsilon^{1-\delta_1}} \frac{d\theta_2}{|e^{i\theta_2}-e^{i\theta_1}|^{\gamma^2}}
\EE\left[\left(\int_{|e^{iu}-e^{i\theta_1}| \ge \epsilon^{1-\delta_1}}  \frac{e^{\gamma^2 \EE[X_\epsilon(e^{i\theta_2}) X_\epsilon(e^{iu})]}M_{\gamma, \epsilon}(du)}{|e^{i u} - e^{i\theta_1}|^{\gamma^2}}\right)^{k-2} \right].
\end{align*}

\noindent Let us now pick $\delta_2 > 0$ such that $2\delta_2 < \delta_1$, and consider
\begin{align*}
&\EE\left[\left(\int_{|e^{iu}-e^{i\theta_1}| \ge \epsilon^{1-\delta_1}, |e^{iu} - e^{i\theta_2}| \le \epsilon^{1-\delta_2}}  \frac{e^{\gamma^2 \EE[X_\epsilon(e^{i\theta_2}) X_\epsilon(e^{iu})]}M_{\gamma, \epsilon}(du)}{|e^{i u} - e^{i\theta_1}|^{\gamma^2}}\right)^{k-2} \right]\\
& \le C \EE\left[\left(\int_{|e^{iu}-e^{i\theta_1}| \ge \epsilon^{1-\delta_1}, |e^{iu} - e^{i\theta_2}| \le \epsilon^{1-\delta_2}}  \frac{M_{\gamma, \epsilon}(du)}{|e^{i u} - e^{i\theta_1}|^{\gamma^2} (|e^{i u} - e^{i\theta_2}| \vee \epsilon)^{\gamma^2}}\right)^{k-2} \right].
\end{align*}

\noindent For $\epsilon > 0$ sufficiently small, we may assume without loss of generality that $\epsilon^{1-\delta_1} \ge 10 \epsilon^{1-2\delta_2}$  so that
\begin{align*}
|e^{iu} - e^{i\theta_1}| 
&\ge |e^{i\theta_2} - e^{i\theta_1}|  - |e^{iu} - e^{i\theta_2}|
\ge \frac{1}{2} |e^{i\theta_2} - e^{i\theta_1}|, \\
\text{and}\qquad  |e^{i\theta_2} - e^{i\theta_1}|
&\ge \frac{1}{2} \left[|e^{i u} - e^{i\theta_2}|+ |e^{i\theta_2} - e^{i\theta_1}|\right]
\ge \frac{1}{2}|e^{iu} - e^{i\theta_1}|
\end{align*}

\noindent when $|e^{iu} - e^{i\theta_2}| \le \epsilon^{1-2\delta_2} \le \epsilon^{1-\delta_1} \le |e^{i\theta_2} - e^{i\theta_1}|$. This means that
\begin{align*}
& \EE\left[\left(\int_{|e^{iu}-e^{i\theta_1}| \ge \epsilon^{1-\delta_1}, |e^{iu} - e^{i\theta_2}| \le \epsilon^{1-\delta_2}}  \frac{M_{\gamma, \epsilon}(du)}{|e^{i u} - e^{i\theta_1}|^{\gamma^2} (|e^{i u} - e^{i\theta_2}| \vee \epsilon)^{\gamma^2}}\right)^{k-2} \right]\\
& \qquad \le C |e^{i\theta_2} - e^{i\theta_1}|^{-(k-2) \gamma^2}\EE\left[\left(\int_{|e^{iu} - e^{i\theta_2}| \le \epsilon^{1-\delta_2}}  \frac{M_{\gamma, \epsilon}(du)}{ (|e^{i u} - e^{i\theta_2}| \vee \epsilon)^{\gamma^2}}\right)^{k-2} \right]\\
& \qquad \le C |e^{i\theta_2} - e^{i\theta_1}|^{-(k-2) \gamma^2} \epsilon^{(1-\delta_2)\frac{\gamma^2}{2}(k-2)}\\
& \qquad \le C |e^{i\theta_2} - e^{i\theta_1}|^{-(k-2) \gamma^2} \epsilon^{\delta_2 \frac{\gamma^2}{2}(k-2)}\EE\left[\left(\int_{|e^{iu} - e^{i\theta_2}| \le \epsilon^{1-2\delta_2}}  \frac{M_{\gamma, \epsilon}(du)}{ (|e^{i u} - e^{i\theta_2}| \vee \epsilon)^{\gamma^2}}\right)^{k-2} \right]\\
& \qquad \le C  \epsilon^{\delta_2 \frac{\gamma^2}{2}(k-2)}\EE\left[\left(\int_{|e^{iu} - e^{i\theta_2}| \le \epsilon^{1-2\delta_2}}  \frac{M_{\gamma, \epsilon}(du)}{(|e^{i u} - e^{i\theta_1}|\vee \epsilon)^{\gamma^2} (|e^{i u} - e^{i\theta_2}| \vee \epsilon)^{\gamma^2}}\right)^{k-2} \right]\\
\end{align*}

\noindent where the second and third inequalities follow from the multifractal estimates (\Cref{lem:GMC_multi}) again. In particular,
\begin{align*}
& \int d\theta_1 \int_{|e^{i\theta_2} - e^{i\theta_1}| \ge \epsilon^{1-\delta_1}} \frac{d\theta_2}{|e^{i\theta_2}-e^{i\theta_1}|^{\gamma^2}}
 \EE\left[\left(\int_{\substack{|e^{iu}-e^{i\theta_1}| \ge \epsilon^{1-\delta_1} \\|e^{iu} - e^{i\theta_2}| \le \epsilon^{1-\delta_2}}}  \frac{M_{\gamma, \epsilon}(du)}{|e^{i u} - e^{i\theta_1}|^{\gamma^2} (|e^{i u} - e^{i\theta_2}| \vee \epsilon)^{\gamma^2}}\right)^{k-2} \right]\\
& \qquad \le C \epsilon^{\delta_2\frac{\gamma^2}{2}(k-2)}\int d\theta_1 \int_{|e^{i\theta_2} - e^{i\theta_1}| \ge \epsilon^{1-\delta_1}} \frac{d\theta_2}{|e^{i\theta_2}-e^{i\theta_1}|^{\gamma^2}}\\
& \qquad \qquad \times 
\EE\left[\left(\int_{|e^{iu} - e^{i\theta_2}| \le \epsilon^{1-2\delta_2}}  \frac{M_{\gamma, \epsilon}(du)}{(|e^{i u} - e^{i\theta_1}|\vee \epsilon)^{\gamma^2} (|e^{i u} - e^{i\theta_2}| \vee \epsilon)^{\gamma^2}}\right)^{k-2} \right]\\
& \qquad \le C \epsilon^{\delta_2\frac{\gamma^2}{2}(k-2)} \EE\left[M_{\gamma, \epsilon}([0, 2\pi)^k\right] \le C \epsilon^{\delta_2\frac{\gamma^2}{2}(k-2)}  \log \frac{1}{\epsilon},
\end{align*}

\noindent and using the basic estimate from \Cref{lem:basic} we have
\begin{align*}
& \int d\theta_1 \int_{|e^{i\theta_2} - e^{i\theta_1}| \ge \epsilon^{1-\delta_1}} \frac{d\theta_2}{|e^{i\theta_2}-e^{i\theta_1}|^{\gamma^2}}
\EE\left[\left(\int_{|e^{iu}-e^{i\theta_1}| \ge \epsilon^{1-\delta_1}}  \frac{e^{\gamma^2 \EE[X_\epsilon(e^{i\theta_2}) X_\epsilon(e^{iu})]}M_{\gamma, \epsilon}(du)}{|e^{i u} - e^{i\theta_1}|^{\gamma^2}}\right)^{k-2} \right]\\
& = \int d\theta_1 \int_{|e^{i\theta_2} - e^{i\theta_1}| \ge \epsilon^{1-\delta_1}} \frac{d\theta_2}{|e^{i\theta_2}-e^{i\theta_1}|^{\gamma^2}}
 \EE\left[\left(\int_{\substack{|e^{iu}-e^{i\theta_1}| \ge \epsilon^{1-\delta_1} \\|e^{iu} - e^{i\theta_2}| \ge \epsilon^{1-\delta_2}}}  \frac{M_{\gamma, \epsilon}(du)}{|e^{i u} - e^{i\theta_1}|^{\gamma^2} |e^{i u} - e^{i\theta_2}|^{\gamma^2}}\right)^{k-2} \right]\\
& \qquad + \Oa\left(c(\delta_2; \epsilon) \log \frac{1}{\epsilon}\right)
\end{align*}

\noindent for some suitable $c(\delta_2; \epsilon)$. This procedure can be repeated until we bring down all the powers, leading to the identity
\begin{align*}
\EE\left[M_{\gamma, \epsilon}([0, 2\pi))^k\right]
& = \int_{D_{(\delta_1, \dots, \delta_{k-1})}} \frac{d\theta_1 \dots d\theta_k}{\prod_{1 \le u < v \le k} |e^{i\theta_u} - e^{i\theta_v}|^{\gamma^2}} + \Oa\left(c(\delta_1, \dots, \delta_{k-1}; \epsilon) \log \frac{1}{\epsilon}\right)\\
&= \EE\left[\int_{D_{(\delta_1, \dots, \delta_{k-1})}} \prod_{j=1}^k M_{\gamma, \epsilon}(d\theta_j)\right] + \Oa\left(c(\delta_1, \dots, \delta_{k-1}; \epsilon) \log \frac{1}{\epsilon}\right)
\end{align*}

\noindent for some suitable $c(\delta_1, \dots, \delta_{k-1}; \epsilon)$ and this necessarily implies the validity of \eqref{eq:FH_goal}. Since all these $\delta_j$'s are arbitrary,
we conclude that
\begin{align*}
\mathrm{MoM}_{\mathrm{U}(N)}(k, s)
& = (1+o(1)) \frac{\EE_{\mathrm{U(N)}}[|P_N(0)|^{2s}]^k}{(2\pi)^k}
\EE\left[M_{\gamma, \epsilon}([0, 2\pi))^k \right]
\end{align*}

\noindent as $\epsilon \to 0^+$ and our proof is complete.
\end{proof}

\appendix

\section{Circular unitary ensemble and subcritical GMCs} \label{app:subsub}
Earlier in the introduction, we mentioned that the asymptotics \eqref{eq:cue_subsub} for $\mathrm{MoM}_{\mathrm{U}(N)}(k, s)$ is essentially resolved in the subcritical regime where $k < 1/s^2$ and $0 < s < 1$. By essentially resolved we meant that the GMC argument could be made rigorous provided that one could establish the convergence of moments 
\begin{align}\label{eq:momconv}
\lim_{N \to \infty} \EE_{\mathrm{U}(N)} \left[ \left(\int_0^{2\pi} \frac{|P_N(\theta)|^{2s}}{\EE_{\mathrm{U}(N)}|P_N(\theta)|^{2s}} d\theta\right)^{k}\right]
= \EE\left[ \left(\int_0^{2\pi} M_{\gamma}(d\theta)\right)^k\right],
\end{align}

\noindent where $\gamma = \sqrt{2}s$ and $M_{\gamma}(d\theta)  =e^{\gamma X(e^{i\theta}) - \frac{\gamma^2}{2}\EE[X(e^{i\theta})^2]} d\theta$ is the GMC measure associated to the Gaussian free field \eqref{eq:cuegff}. In many situations, the condition \eqref{eq:momconv} can be verified. For instance, if $k \in \NN$, we can consider the expanded moments \eqref{eq:expandedmom} and show that the main contribution to the multiple integral comes from the set where the angles $\theta_1, \dots, \theta_k$ are at a macroscopic distance apart from each other, in which case existing asymptotics for Toeplitz determinant with singularities may be applied. This was done as part of \cite[Corollary 2.1]{Fa2019}.

For non-integer values of $k$, the expanded moment approach no longer works. One may want to rely on the weak convergence
\begin{align}
\notag
\frac{|P_N(\theta)|^{2s}}{\EE_{\mathrm{U}(N)}|P_N(\theta)|^{2s}} d\theta 
& =: \mu_{N,s}(d\theta)\\
\label{eq:cuegmccon}
&\xrightarrow[N \to \infty]{d} M_{\gamma}(d\theta)  =e^{\gamma X(e^{i\theta}) - \frac{\gamma^2}{2}\EE[X(e^{i\theta})^2]} d\theta, \qquad \theta \in [0, 2\pi),
\end{align}

\noindent but this distributional result per se does not immediately imply \eqref{eq:momconv}. Fortunately, the proof of \eqref{eq:cuegmccon} is based on a moment method \cite{NSW2018, We2015}, and stronger conclusions may be obtained. For example, in the $L^2$-regime \cite{We2015} where $s^2 < 1/2$ , there exists a further sequence of random measures $\mu_{N, s}^{(M)}(d\theta)$ such that
\begin{align*}
\lim_{M \to \infty} \lim_{N \to \infty}\EE_{\mathrm{U}(N)}\left[ \left(\mu_{N, s}([0, 2\pi))-\mu_{N, s}^{(M)}([0, 2\pi))\right)^2\right] = 0
\end{align*}

\noindent and
\begin{align*}
\lim_{M \to \infty} \lim_{N \to \infty}\EE_{\mathrm{U}(N)}\left[ \left(\mu_{N, s}^{(M)}([0, 2\pi))\right)^k\right] = \EE\left[ \left(\int_0^{2\pi} M_{\gamma}(d\theta)\right)^k\right] \qquad \text{for any }k \in (0, 1/s^2).
\end{align*}

\noindent Using these ingredients, we may conclude that for any $1 \le k \le 2$,
\begin{align*}
& \EE_{\mathrm{U}(N)} \left[ \left(\int_0^{2\pi} \frac{|P_N(\theta)|^{2s}}{\EE_{\mathrm{U}(N)}|P_N(\theta)|^{2s}} d\theta\right)^{k}\right]^{1/k}\\
& =\EE_{\mathrm{U}(N)}\left[ \left(\mu_{N, s}^{(M)}([0, 2\pi))\right)^k\right]  + \mathcal{O}\left( \EE_{\mathrm{U}(N)}\left[ \left(\mu_{N, s}([0, 2\pi))-\mu_{N, s}^{(M)}([0, 2\pi))\right)^2\right]^{1/2}\right)\\
& \xrightarrow{N \to \infty, M \to\infty} \EE\left[ \left(\int_0^{2\pi} M_{\gamma}(d\theta)\right)^k\right]
\end{align*}

\noindent and the case where $k \in (0, 1)$ is similar. This approach is generalised in \cite{NSW2018} to the $L^1$-regime where $s^2 \in [1/2, 1)$, from which one could obtain \eqref{eq:momconv} at least for all $k \le 1$. Of course, these results have not covered the entire subcritical regime $k < 1/s^2$ in the asymptotic analysis of moments of moments \eqref{eq:CUEmom} of CUE matrices, but we expect that more refined analysis along a similar line of argument may make this possible and that the problem is conceptually solved.

\section{GMC moments: renormalisation from below} \label{app:lowercont}
Here we briefly mention the `lower continuity' of GMC moments up to renormalisation at critical exponent in the following sense.
\begin{proposition}
Let $d \le 2$, $\gamma \in (0, \sqrt{2d})$, $p_c = 2d / \gamma^2$ and $M_{\gamma}(dx)$ be the GMC measure associated to the log-correlated field \eqref{eq:LGF}. If $g \ge 0$ is a continuous function on $\overline{D}$, then
\begin{align*}
\lim_{\alpha \to p_c^-} (p_c-1-\alpha)
\EE\left[\left(\int_D g(x) M_{\gamma}(dx)\right)^{\alpha}\right]
= \left(\int_D  e^{d(p_c-1)f(u, u)} g(u)^{p_c} du\right) (p_c-1) \overline{C}_{\gamma, d}
\end{align*}

\noindent and
\begin{align*}
\lim_{\widetilde{\gamma} \to \gamma^-} \left(\frac{2d}{\widetilde{\gamma}^2} - p_c\right)
\EE\left[\left(\int_D g(x) M_{\widetilde{\gamma}}(dx)\right)^{p_c}\right]
= \left(\int_D  e^{d(p_c-1)f(u, u)} g(u)^{p_c} du\right) (p-1) \overline{C}_{\gamma, d}.
\end{align*}
\end{proposition}

\begin{proof}
We only sketch the proof of the first limit for $g \equiv 1$ and $f \ge 0$. Let us fix $r \in (0, 1/4)$ as before. By 
considerations similar to those in \Cref{subsec:removal}, we have
\begin{align}
\notag
\EE\left[\left(\int_D M_{\gamma}(dx)\right)^{\alpha}\right]
&= \int_D du \EE\left[\left(\int_D \frac{e^{\gamma^2 f(x, u)}}{|x-u|^{\gamma^2}} M_{\gamma}(dx)\right)^{\alpha-1}\right]\\
\label{eq:lower_pf1}
& = \mathcal{O}(1) + \int_D du \EE\left[\left(\int_{D\cap \{|x-u| \le r\}} \frac{e^{\gamma^2 f(x, u)}}{|x-u|^{\gamma^2}} M_{\gamma}(dx)\right)^{\alpha-1}\right].
\end{align}

Given the uniform continuity of $f$ on $\overline{D} \times \overline{D}$, it has a modulus of continuity $\omega(\cdot)$ so that
\begin{align*}
|f(x, u) - f(u,u)| \le \omega(r) \qquad \forall x \in B(u, r) \quad \forall u \in D
\end{align*}

\noindent and $\omega(r) \to 0$ as $r \to 0^+$. Based on the same argument of Gaussian comparison we saw in \Cref{lem:reduced_exact}, \eqref{eq:lower_pf1} can be rewritten as
\begin{align*}
& \mathcal{O}(1) + \left(1+\mathcal{O}(\omega(r))\right)\left(\int_D du  \EE\left[\left(e^{\gamma^2 f(u, u)} e^{\gamma \Na_{u} - \frac{\gamma^2}{2} \EE[\Na_{u}^2]}\right)^{\alpha-1}\right]\right)
\EE \left[\left(\int_{|x|\le r} e^{\gamma \overline{X}(x) - \frac{\gamma^2}{2} \EE[\overline{X}(x)^2]}\frac{dx}{|x|^{\gamma^2}} \right)^{\alpha} \right]\\
& = \mathcal{O}(1) + \left(1+\mathcal{O}(\omega(r))\right)\left(\int_D  e^{\frac{\gamma^2}{2}\alpha(\alpha-1) f(u, u)} du\right)
\EE \left[\left(\int_{|x|\le r} e^{\gamma \overline{X}(x) - \frac{\gamma^2}{2} \EE[\overline{X}(x)^2]}\frac{dx}{|x|^{\gamma^2}} \right)^{\alpha} \right].
\end{align*}

\noindent Our claim follows immediately by first applying \Cref{lem:lower} and then sending $r$ to $0$.
\end{proof}

\section{Supercritical moments of approximate GMCs} \label{app:supercritical}
Let us again consider $X_\epsilon := X \ast \nu_\epsilon$ for some mollifier $\nu \in C_c^\infty(\RR^d)$ and study the integrand on the RHS of
\begin{align*}
\EE\left[M_{\gamma, \epsilon}(D)^{p}\right]
= \int_D du\EE\left[\left(\int_D \frac{e^{\gamma^2 f_\epsilon(x, u)}}{\left(|x-u|\vee \epsilon\right)^{\gamma^2}} e^{\gamma X_{\epsilon}(x) - \frac{\gamma^2}{2} \EE[X_\epsilon(x)^2]} dx\right)^{p-1}\right].
\end{align*}

Suppose we are in the `supercritical' case where $p > 2d/\gamma^2$ for any $\gamma \in \RR$. Looking at (part of) the microscopic contribution, one can show that
\begin{align*}
\EE\left[\left(\int_{|x-u| \le \epsilon} \frac{e^{\gamma^2 f_\epsilon(x, u)}}{\left(|x-u|\vee \epsilon\right)^{\gamma^2}} e^{\gamma X_{\epsilon}(x) - \frac{\gamma^2}{2} \EE[X_\epsilon(x)^2]} dx\right)^{p-1}\right]
= \Omega\left(\epsilon^{\frac{\gamma^2}{2}(p-1)\left(\frac{2d}{\gamma^2} - p\right)}\right)
\end{align*}

\noindent captures the size of the leading order of $\EE[M_{\gamma, \epsilon}(D)^p]$. On the other hand, for any $\delta > 0$, it is possible to show that the contribution from $\{|x-u| \ge \epsilon^{1-\delta}\}$ does not yield the leading order, so that
\begin{align*}
&\EE\left[\left(\int_D \frac{e^{\gamma^2 f_\epsilon(x, u)}}{\left(|x-u|\vee \epsilon\right)^{\gamma^2}} e^{\gamma X_{\epsilon}(x) - \frac{\gamma^2}{2} \EE[X_\epsilon(x)^2]} dx\right)^{p-1}\right]\\
& \qquad \overset{\epsilon \to 0^+}{\sim }
\EE\left[\left(\int_{D\cap \{|x-u| \le \epsilon^{1-\delta}\}} \frac{e^{\gamma^2 f_\epsilon(x, u)}}{\left(|x-u|\vee \epsilon\right)^{\gamma^2}} e^{\gamma X_{\epsilon}(x) - \frac{\gamma^2}{2} \EE[X_\epsilon(x)^2]} dx\right)^{p-1}\right]
\end{align*}

\noindent by \Cref{lem:basic}. Since $\delta > 0$ is arbitrary, a finer analysis will arrive at the conclusion that the leading coefficient depends only on the microscopic regime where $|x-u| = \mathcal{O}(\epsilon)$.

To proceed further, let us define $X_\epsilon(x; u) = X_\epsilon(x) - X_\epsilon(u)$ and rewrite the last expression as
\begin{align}
\notag
& \EE\left[\left(\int_{D\cap\{ |x-u| \le \epsilon^{1-\delta}\}} e^{\gamma^2 \EE[X_\epsilon(x) X_\epsilon(u)]}e^{\gamma X_{\epsilon}(x) - \frac{\gamma^2}{2} \EE[X_\epsilon(x)^2]} dx\right)^{p-1}\right]\\
\notag
& = e^{\frac{\gamma^2}{2}p(p-1) \EE[X_\epsilon(u)^2]} \\
\notag
& \qquad \times \EE\Bigg[e^{\gamma(p-1) X_\epsilon(u) - \frac{\gamma^2(p-1)^2}{2} \EE[X_\epsilon(u)^2]}\left(\int_{D\cap\{ |x-u| \le \epsilon^{1-\delta}\}} e^{\gamma X_{\epsilon}(x;u) - \frac{\gamma^2}{2} \EE[X_\epsilon(x;u)^2]} dx\right)^{p-1}\Bigg]\\
\label{eq:supercritical1}
& = e^{\frac{\gamma^2}{2}p(p-1) \EE[X_\epsilon(u)^2]}
 \EE\Bigg[\left(\int_{D\cap \{|x-u| \le \epsilon^{1-\delta}\}} e^{\gamma^2(p-1) \EE[X_\epsilon(u) X_\epsilon(x; u)]}e^{\gamma X_{\epsilon}(x;u) - \frac{\gamma^2}{2} \EE[X_\epsilon(x;u)^2]} dx\right)^{p-1}\Bigg]
\end{align}

\noindent by Cameron-Martin (\Cref{lem:CMG}). We have
\begin{align*}
\EE[X_\epsilon(u) X_\epsilon(x; u)]
& = - \log\left(\frac{|u-x|}{\epsilon} \vee 1\right) + \left(f_\epsilon(u, x) - f_\epsilon(u, u)\right), \\
\EE[X_\epsilon(u+\epsilon s; u) X_\epsilon(u + \epsilon t; u)]
& = \log \left(|s| \vee 1\right) + \log \left(|t| \vee 1\right) - \log \left(|s-t| \vee 1\right)\\
& \qquad + f_\epsilon(u, u) + f_\epsilon(u+\epsilon s, u+ \epsilon t)
- f_\epsilon(u, u+ \epsilon s) - f_\epsilon(u, u+ \epsilon t).
\end{align*}

\noindent In particular, the Gaussian field $\mathbf{X}_\epsilon(s; u) := X_\epsilon(u + \epsilon s; u)$ does not have log-singularity anymore as $\epsilon \to 0^+$, and with a simple change of variable we see that \eqref{eq:supercritical1} is asymptotically equal to
\begin{align*}
& e^{\frac{\gamma^2}{2}p(p-1) \EE[X_\epsilon(u)^2]} \epsilon^{d(p-1)}\\
&\quad \times
 \EE\Bigg[\left(\int_{\RR^d}\frac{e^{\gamma^2(p-1) \left[f_\epsilon(u, u + \epsilon s) - f_\epsilon(u, u)\right]}}{\left(|s| \vee 1\right)^{\gamma^2 (p-1)}} 
\frac{e^{-\frac{\gamma^2}{2}\left[f_\epsilon(u,u) + f_\epsilon(u+\epsilon s, u+\epsilon s) - 2f_\epsilon(u, u+\epsilon s)\right]}}{\left(|s| \vee 1\right)^{\gamma^2}}
e^{\gamma \mathbf{X}_{\epsilon}(s;u)} ds\right)^{p-1}\Bigg]\\
& = \epsilon^{\frac{\gamma^2}{2}(p-1)\left(\frac{2d}{\gamma^2} - p\right)} e^{-\frac{\gamma^2}{2}(p-1)^2f_\epsilon(u, u)}
 \EE\Bigg[\left(\int_{\RR^d}\frac{e^{p\gamma^2f_\epsilon(u, u + \epsilon s) - \frac{\gamma^2}{2} f_\epsilon(u+\epsilon s, u+\epsilon s)}}{\left(|s| \vee 1\right)^{p\gamma^2}} e^{\gamma \mathbf{X}_{\epsilon}(s;u)}ds\right)^{p-1}\Bigg].
\end{align*}

\noindent The leading order coefficient remains $\Omega(1)$ as $\epsilon \to 0^+$, but it is very sensitive to the choice of approximate fields $X_\epsilon$ (or equivalently the mollifiers $\nu$). In addition, the Gaussian quantity here is not able to accurately capture the properties of e.g. random matrix models, the microscopic statistics of which are described by the Sine process instead.

Let us also quickly comment on the `critical-supercritical' case, i.e. $p = 2d/\gamma^2 \in (0, 1)$. If we look at the toy model of exponential functional of Brownian motion, i.e.
\begin{align*}
\EE\left[\left(\int_0^T e^{\gamma (B_t - \mu t)}dt\right)^{p-1}\right]
\end{align*}

\noindent we now have a Brownian motion with drift $-\mu = -\frac{\gamma}{2}(p-1) > 0$, so that the leading order contribution to the above integral should come from near $T = -(1+o(1))\log \epsilon$. In terms of the GMC asymptotics, it again suggests that the leading order depends heavily on the microscopic region.

\section{Slowly merging Fisher-Hartwig singularities}\label{app:FH}
We now discuss why the strong asymptotics in \Cref{lem:FH} (i) continue to hold when we allow singularities $\theta_1, \dots, \theta_k$ to possibly merge at a a rate not faster than $N^{-(1-\delta)}$. Since this requires an inspection of the Riemann-Hilbert analysis that is largely independent of the main text of the current article, we shall follow the notations in \cite{DIK2011, DIK2014} as closely as possible in this appendix unless otherwise specified. Much of the following content is rather standard in this kind of study, and we shall only briefly explain the arguments with heavy reference to existing literature where interested readers can find more details.

To begin with, it is helpful to have a quick recap of the proof of the strong asymptotics of Toeplitz determinant with Fisher-Hartwig singularities.\footnote{We do not include the jump singularities here as they are not needed in our discussion.} Letting\footnote{Unlike \cite{DIK2011,DIK2014} we do not include the singularity $z_0 = 1$ a priori, as we do not need the connections to Hankel determinants.}
\begin{align*}
f(z) = \prod_{j=1}^m |z-z_j|^{2\alpha_j}
\qquad \text{where } z_j = e^{i\theta_j}, \quad j = 1, \dots, m, \quad 0 \le \theta_1 < \dots < \theta_m < 2\pi
\end{align*}

\noindent where $\Re(\alpha_j) > -\frac{1}{2}$ for $j = 1, \dots, m$, and
\begin{align*}
f_j = \frac{1}{2\pi}\int_0^{2\pi} f(e^{i\theta}) e^{-ij\theta}d\theta, \qquad j \in \ZZ,
\end{align*}

\noindent we are interested in the asymptotics for $n$-dimensional Toeplitz determinant with symbol $f$, i.e. $D_n(f) = D_n(\alpha_1, \dots, \alpha_m) := \det\left(f_{j-k}\right)_{j,k=0}^{n-1}$, as $n \to \infty$. This is closely related to Haar-distributed unitary matrices since the Heine-Szeg\"o identity (see e.g. \cite[Equation (2.10)]{DIK2014}) implies that
\begin{align*}
\EE_{\mathrm{U(N)}}\left[ \prod_{j=1}^k |\det(I - U_Ne^{-i\theta_j})|^{2s}\right]
= D_{N-1}(\underbrace{s, \dots, s}_{\text{$k$ terms}}).
\end{align*}

\noindent Furthermore, the Toeplitz determinants are connected to orthogonal polynomials on the unit circle in the following way: there exists a dual pair of orthogonal polynomials\footnote{This is true when all the $\alpha_j$'s are real, which is all we need in our case; for complex-valued $\alpha_j$'s see e.g. the discussion after \cite[Equation (2.10)]{DIK2014}.}
\begin{align*}
\phi_k(z) = \chi_k z^k + \dots
\qquad \text{and} \qquad
\widehat{\phi}_k(z) = \chi_k z^k + \dots
\end{align*}

\noindent of degree $k$ with positive leading coefficient $\chi_k$  for each $k \in \NN$ such that (writing $z=e^{i\theta}$)
\begin{align*}
\frac{1}{2\pi}\int_0^{2\pi} \phi_k(z) z^{-j} f(z) d\theta = \chi_k^{-1} \delta_{jk},
\qquad \frac{1}{2\pi}\int_0^{2\pi} \widehat{\phi}_k(z) z^{j} f(z) d\theta = \chi_k^{-1} \delta_{jk}, \qquad j = 0, 1, \dots, k
\end{align*}

\noindent and it can be shown that (see e.g. \cite[Section 2]{DIK2014})
\begin{align*}
\chi_k = \sqrt{\frac{D_k(f)}{D_{k+1}(f)}}
\qquad \text{so that}
\qquad D_n(f) := \prod_{j=1}^{n-1} \chi_j^{-2}.
\end{align*}

The Riemann-Hilbert analysis kicks in when we try to encode these objects in some matrix-valued functions: for each $k \in \NN$, one can define
\begin{align*}
Y(z) = Y^{(k)}(z) = \begin{pmatrix}
\chi_k^{-1} \phi_k(z) & \chi_k^{-1} \int_C \frac{\phi_k(\xi)}{\xi - z} \frac{f(\xi)d\xi}{2\pi \xi^k}\\
-\chi_{k-1} z^{k-1}\widehat{\phi}_{k-1}(z^{-1}) & -\chi_{k-1} \int_C \frac{\widehat{\phi}_{k-1}(\xi^{-1})}{\xi - z} \frac{f(\xi)d\xi}{2\pi \xi}
\end{pmatrix}
\qquad C = \text{unit circle}
\end{align*}

\noindent and show that it is the unique solution to the so-called $Y$ Riemann-Hilbert problem (RHP), see \cite[Equations (2.5-2.9)]{DIK2014}. One can then show that (see \cite[Equation (5.3)]{DIK2014} and the references therein) for $\gamma \in \cup_{j=1}^m \{\alpha_j\}$, we have
\begin{align}
\notag
\frac{\partial}{\partial \gamma} \log D_n(f)
& = -2\frac{\frac{\partial \chi_n}{\partial \gamma}}{\chi_n}
\left(n + \sum_{j=1}^m \alpha_j \left\{1 - \widetilde{Y}_{12}^{(n)}(z_j) Y_{21}^{(n+1)}(z_j) - z_j Y_{11}^{(n)}(z_j) \widetilde{Y}_{22}^{(n+1)}(z_j)\right\}\right)\\
\label{eq:diffid}
& \qquad - 2\sum_{j=1}^m \alpha_j \left(\widetilde{Y}_{12}^{(n)}(z_j) \frac{\partial}{\partial \gamma}Y_{21}^{(n+1)}(z_j) - z_j \frac{\partial}{\partial \gamma}Y_{11}^{(n)}(z_j) \widetilde{Y}_{22}^{(n+1)}(z_j)\right)
\end{align}

\noindent where $Y^{(k)}_{ij}$ is the $ij$-th entry of the solution to the $Y^{(k)}$-RHP, and $\widetilde{Y}^{(k)}_{ij}$ are the regularised versions.\footnote{See the discussion between Equations (3.16) and (3.17) in \cite{DIK2014}.} In particular, the RHS of \eqref{eq:diffid} only depends on the `last two' RHPs $Y^{(n)}$ and $Y^{(n+1)}$. Provided that we can establish their asymptotics (including any derivatives that have appeared in the differential identity) as $n \to \infty$, we can obtain our desired formula by integrating both sides with respect to each variable $\gamma = \alpha_j$, $j = 1, \dots, m$, using the fact that $D_n(0, \dots, 0) = 1$.

So far everything we have discussed applies regardless of the behaviour of the singularities, and it all boils down to the analysis of the asymptotics for $Y$- RHP. This relies on the method of non-linear steepest descent, first applied in the random matrix context in \cite{DKMVZ1999}, which involves the following steps:
\begin{itemize}
\item The $Y$-RHP is transformed to the so-called $T$-RHP which `normalises the behaviour at infinity', see e.g. \cite[Equation (4.1)]{DIK2014}.
\item The $T$-RHP is further transformed to a new $S$-RHP, a procedure known as opening of lenses, see e.g.  \cite[Equation (4.2)]{DIK2014} and the discussion that follows it.
\item The $S$-RHP is approximated by two sets of functions:
\begin{itemize}
\item  the global parametrix $N(z)$ in \cite[Equations (4.7-4.10)]{DIK2014} when $z$ is away from the unit circle and not close to any of the singularities;
\item a collection of functions $P_{z_j}(z)$, known as the the local parametrices (see \cite[Equation (4.16)]{DIK2014} and other related definitions), each of which is defined on a small neighbourhood $U_{z_j}$ of $z_j = e^{i\theta_j}$, say $U_{z_j} := \{z \in \CC: |z - z_j| < \epsilon\}$ for some $\epsilon > 0$ chosen such that $U_{z_j}$'s are non-overlapping.
\end{itemize}
\item The $R$-RHP \cite[Equation (4.26)]{DIK2014}, which measures the `multiplicative error' when the $S$-RHP is approximated by the global/local parametrices, is shown to be close to the identity matrix with other `nice properties' as $n \to \infty$.
\end{itemize}

When the singularities $z_1, \dots, z_m$ are fixed, the parameter $\epsilon > 0$ in the definition of $U_{z_j}$ may be chosen as some fixed number, and one can easily verify the `jump condition' for the $R$-RHP in \cite[Equation (4.27-4.29)]{DIK2014} is of the form $R_+(z) = R_-(z)\Delta(z)$ with $\Delta(z) = I + \mathcal{O}(n^{-1})$. Standard analysis then shows that $R$ has a unique solution and that $R-I$ has an asymptotic expansion in powers of $n^{-1}$, see more details in \cite[Section 4.1]{DIK2014}.

When the singularities are allowed to vary but remain $n^{-(1-\kappa)}$ (with $\kappa> 0$) separated from each other\footnote{The separation distance $n^{-(1-\kappa)}$ here should be seen as our $N^{-(1-\delta)}$ in the main text of the current paper.}, the first three steps of the steepest descent will still go through except that the local parametrices $P_{z_j}$ will now be defined on $U_{z_j}$ with e.g. radius $\epsilon = \epsilon_n = \frac{1}{10} n^{-(1-\kappa)}$ in order to satisfy the non-overlapping condition. The analysis of $R$-RHP with shrinking contours already appeared in \cite{DIK2011,DIK2014} for the `extension to smooth $V(z)$', even though in the context of merging singularities they were perhaps first analysed in \cite{CK2015} and more recently in \cite{CFLW2019,Fa2019} The consequence of shrinking contours (in the notations of \cite{DIK2011,DIK2014}) are as follows:
\begin{itemize}
\item The jump conditions \cite[Equations (4.27-4.28)]{DIK2014} remain true, but now the jump matrix
\begin{align*}
\Delta(z) := I + N(z) \begin{pmatrix} 0 & 0 \\ f(z)^{-1} z^{\mp n} & 0\end{pmatrix} N(z)^{-1}
\end{align*}

\noindent (where the exponent for $z$ is negative if $z \in \Sigma_j^{\mathrm{out}}$ and positive if $z \in \Sigma_j^{''\mathrm{out}}$ for $j=1, \dots, m$) is no longer $I+\mathcal{O}(e^{-c n})$ for some small $c > 0$, but rather $I + \mathcal{O}(e^{-c n^c})$ where $c>0$ may possibly depend on $\kappa$ but can be made uniform in compact sets of $\alpha_1, \dots, \alpha_m$.
\item The jump condition \cite[Equation (4.29)]{DIK2014} is also true, but now the matching condition for the global/local parametrices is weakened to
\begin{align*}
\Delta(z) := P_{z_j}(z) N(z)^{-1} = I + \mathcal{O}(1/n\epsilon_n) = I + \mathcal{O}(n^{-\kappa}), \qquad \forall z \in \partial U_{z_j} \setminus \{ \text{intersection points} \}
\end{align*}

\noindent for each $j=1, \dots, m$.
\end{itemize}

While $R(z) - I$ can no longer be expanded in an asymptotic series in powers of $n^{-1}$, the $R$-RHP is still a `small-norm problem' and so for $n$ sufficiently large it may still be solved using the Neumann series construction.\footnote{The arguments for shrinking contours are explained in greater detail in \cite[Section 6.4]{Fa2019} for the problem of uniform asymptotics of Toeplitz determinants with arbitrary singularities, as well as \cite[Section 8.1]{CFLW2019} for a problem on Hankel determinants.} In particular, one can still verify that
\begin{align*}
R(z) = I + o(1) \qquad \text{ and } \qquad \frac{\partial}{\partial \gamma} R(z) = o(1)
\end{align*}

\noindent and so all the leading order asymptotics for $\chi_n$, $Y^{(n)}$ and $Y^{(n+1)}$ and its regularised versions (as well as their derivatives) appearing in \eqref{eq:diffid} will remain identical to what was observed in the fixed singularity case, and hence the strong asymptotics for Toeplitz determinants with singularities may be extended to cover the slightly more general situation with a separation distance of $n^{-(1-\kappa)}$ for any $\kappa > 0$.



\begin{thebibliography}{99}
\bibitem{ABK2020}
T. Assiotis, E. C. Bailey and J. P. Keating.
On the moments of the moments of the characteristic polynomials of Haar distributed symplectic and orthogonal matrices.
Preprint arXiv:1910.12576.

\bibitem{AK2020}
T. Assiotis and J.P. Keating.
Moments of moments of characteristic polynomials of random unitary matrices and lattice point counts.
Random Matrices: Theory and Applications. https://doi.org/10.1142/S2010326321500192.

\bibitem{As2020}
T. Assiotis.
On the moments of the partition function of the C$\beta$E field.
Preprint arXiv:2011.10323.

\bibitem{Be2017}
N. Berestycki.
An elementary approach to Gaussian multiplicative chaos.
Electron. Commun. Probab. Volume 22 (2017), paper no. 27, 12 pp.


\bibitem{BHR2019}
M. Dal Borgo, E. Hovhannisyan and A. Rouault:
Mod-Gaussian Convergence for Random Determinants.
Ann. Henri Poincar\'e 20 (2019), 259.298.


\bibitem{BK2019}
E. C. Bailey and J. P. Keating.
On the moments of the moments of the characteristic polynomials of random unitary matrices.
Communications in Mathematical Physics, 371(2):689--726, 2019.

\bibitem{BK2020a}
E. C. Bailey and J. P. Keating.
On the moments of the moments of $\zeta(1/2+it)$.
Preprint arXiv:2006.04503.


\bibitem{BK2020b}
E. C. Bailey and J. P. Keating.
Moments of Moments and Branching Random Walks.
Preprint arXiv:2008.09536.


\bibitem{BW2018}
G. Baverez and M. D. Wong.
Fusion asymptotics for Liouville correlation functions.
Preprint arXiv:1807.10207.

\bibitem{BPNotes}
N. Berestycki and E. Powell.
Lecture notes on \emph{Gaussian free field, Liouville quantum gravity and Gaussian multiplicative chaos}.
Available at \url{https://homepage.univie.ac.at/nathanael.berestycki/Articles/master.pdf}.

\bibitem{BWW2018}
N. Berestycki, C. Webb and M. D. Wong.
Random Hermitian matrices and Gaussian multiplicative chaos.
Probab. Theory Relat. Fields 172, 103-189 (2018). https://doi.org/10.1007/s00440-017-0806-9.



\bibitem{CFLW2019}
T. Claeys, B. Fahs, G. Lambert and C. Webb.
How much can the eigenvalues of a random Hermitian matrix fluctuate?
Preprint arXiv:1906.01561.

\bibitem{CGMY2020}
T. Claeys, G. Glesner, A. Minakov, and M. Yang.
Asymptotics for averages over classical orthogonal ensembles, 
to appear in Int. Math. Res. Notices .

\bibitem{CK2015}
T. Claeys and I. Krasovsky.
Toeplitz determinants with merging singularities.
Duke Mathematical Journal, 164(15):2897--2987, 2015.


\bibitem{CN2019}
R. Chhaibi and J. Najnudel.
On the circle, $\displaystyle \mathrm{GMC}^\gamma = \lim_{\leftarrow} \mathrm{C\beta E}_n$ for $\gamma = \sqrt{\frac{2}{\beta}}, (\gamma \le 1)$.
Preprint arXiv:1904.00578.

\bibitem{DIK2011}
P. Deift, A. Its, and I. Krasovsky. 
Asymptotics of Toeplitz, Hankel, and Toeplitz$+$Hankel determinants with Fisher-Hartwig singularities.
Annals of mathematics (2011): 1243-1299.

\bibitem{DIK2014}
P. Deift, A. Its and I. Krasovsky.
On the asymptotics of a Toeplitz determinant with singularities. Random matrix theory, interacting particle systems and integrable systems, 65, p.93, 2014.

\bibitem{DKMVZ1999}
P. Deift, T. Kriecherbauer, K. T-R McLaughlin, S. Venakides and X. Zhou.
Strong asymptotics of orthogonal polynomials with respect to exponential weights.
Communications on Pure and Applied Mathematics, no. 12 (1999): 1491-1552.

\bibitem{DKRV2016}
F. David, A. Kupiainen, R. Rhodes and V. Vargas.
Liouville quantum gravity on the Riemann sphere.
Comm. Math. Phys. 342 (2016), no. 3, 869-907.

\bibitem{DKRV2017}
F. David, A. Kupiainen, R. Rhodes, and V. Vargas.
Renormalizability of Liouville Quantum Gravity at the Seiberg bound.
Electron. J. Probab. Volume 22 (2017), paper no. 93, 26 pp.

\bibitem{DMS2014}
B. Duplantier, J. Miller and S. Sheffield. 
Liouville quantum gravity as a mating of trees. 
Preprint arXiv:1409.7055.

\bibitem{DS2011}
B. Duplantier and S. Sheffield.
Liouville quantum gravity and KPZ.
Invent. Math. 185 (2011), no. 2, 333-393.

\bibitem{Fa2019}
B. Fahs. Uniform asymptotics of Toeplitz determinants with Fisher--Hartwig singularities. 
Preprint arXiv:1909.07362.


\bibitem{FB2008}
F.V. Fyodorov and  J.P. Bouchaud.
Freezing and extreme value statistics in a Random Energy Model with logarithmically correlated potential.
Journal of Physics A: Mathematical and Theoretical, Volume 41, Number 37, 372001 (2008)


\bibitem{FK2014}
Y. V. Fyodorov and J. P. Keating. 
Freezing transitions and extreme values: random matrix theory, and disordered landscapes. 
Philosophical Transactions of the Royal Society A: Mathematical, Physical and Engineering Sciences, 372(2007):20120503, 2014.


\bibitem{FK2020}
J. Forkel and J. P. Keating.
The Classical Compact Groups and Gaussian Multiplicative Chaos.
Preprint arXiv:2008.07825.

\bibitem{JS2017}
J. Junnila and E. Saksman.
Uniqueness of critical Gaussian chaos.
Electron. J. Probab. 22: 1-31 (2017). DOI: 10.1214/17-EJP28.

\bibitem{JSW2020}
J. Junnila, E. Saksman and C. Webb.
Imaginary multiplicative chaos: Moments, regularity and connections to the Ising model.
Ann. Appl. Probab. 30(5): 2099-2164 (October 2020). DOI: 10.1214/19-AAP1553.

\bibitem{Kah1985}
J.-P. Kahane.
Sur le chaos multiplicatif.
Ann. Sci. Math. Qu\'ebec, 9 no.2 (1985), 105-150.

\bibitem{KN2004}
R. Killip and L. Nenciu. 
Matrix models for circular ensembles. 
Int Math Res Not, 50: 2665-2710 (2004).

\bibitem{KRV2020}
A. Kupiainen, R. Rhodes and V. Vargas (2020). 
Integrability of Liouville theory: proof of the DOZZ Formula. 
Annals of Mathematics, 191(1), 81-166.


\bibitem{KS2000}
J.P. Keating and N.C. Snaith.
Random matrix theory and $\zeta(1/2 + it)$.
Commun. Math. Phys., 214, 57-89, 2000.

\bibitem{La2019}
G. Lambert.
Mesoscopic central limit theorem for the circular beta-ensembles and applications.
Preprint arXiv:1902.06611.

\bibitem{LOS2018}
G. Lambert, D. Ostrovsky and N. Simm. 
Subcritical Multiplicative Chaos for Regularized Counting Statistics from Random Matrix Theory. Commun. Math. Phys. 360, 1-54 (2018). https://doi.org/10.1007/s00220-018-3130-z.

\bibitem{NSW2018}
M. Nikula, E. Saksman and Christian Webb.
Multiplicative chaos and the characteristic polynomial of the CUE: the $L^1$-phase. Preprint arXiv:1806.01831.

\bibitem{Re2020}
G. Remy. 
The Fyodorov-Bouchaud formula and Liouville conformal field theory. Duke Mathematical Journal, 169(1), 177-211. (2020)

\bibitem{RP1981}
L. C. G. Rogers and J. W. Pitman.
Markov functions, Ann. Prob., 9 (1981), 573-582.

\bibitem{RV2010}
R. Robert and V. Vargas.
Gaussian multiplicative chaos revisited.
Ann. Probab. 38(2): 605-631 (March 2010). DOI: 10.1214/09-AOP490.


\bibitem{RV2019}
R. Rhodes and V. Vagras. The tail expansion of Gaussian multiplicative chaos and the Liouville reflection coefficient. 
Ann. Probab. Volume 47, Number 5 (2019), 3082-3107.

\bibitem{RY1998}
D. Revuz and M. Yor. 
Continuous Martingales and Brownian Motion. 3rd edition. 
Grundlehren der mathematischen Wissenschaften Vol 293, Springer Berlin Heidelberg, 1998.


\bibitem{RZ2020a}
The distribution of Gaussian multiplicative chaos on the unit interval.
Ann. Probab. Volume 48, Number 2 (2020), 872-915.


\bibitem{RZ2020b}
G. Remy and T. Zhu.
Integrability of boundary Liouville conformal field theory
Preprint arXiv:2002.05625.


\bibitem{Sha2016}
A. Shamov.
On Gaussian multiplicative chaos.
J. Funct. Anal. 270 (2016), 3224--3261

\bibitem{Sh2016}
S. Sheffield.
Conformal weldings of random surfaces: SLE and the quantum gravity zipper.
Ann. Probab. Volume 44, Number 5 (2016), 3474-3545.


\bibitem{We2015}
C. Webb. 
The characteristic polynomial of a random unitary matrix and Gaussian multiplicative chaos - The $L^2$-phase.
Electron. J. Probab. 20 (2015), no. 104, 1-21.

\bibitem{Wil1974}
D. Williams.
Path decomposition and continuity of local times for one-dimensional diffusions, I. 
Proceed. of London Math. Soc. (3), 28, 738-768, (1974).

\bibitem{Won2020}
M. D. Wong.
Universal tail profile of Gaussian multiplicative chaos.
Probability Theory and Related Fields volume 177, pages 711-746 (2020).

\bibitem{WW2019}
C. Webb and M. D. Wong.
On the moments of the characteristic polynomial of a Ginibre random matrix.
Proc. London Math. Soc., 118: 1017-1056. https://doi.org/10.1112/plms.12225.

\end{thebibliography}
\end{document}